\documentclass[12pt,oneside]{amsart}
\usepackage{amssymb,amsmath,verbatim,amstext,eucal,amscd}

\newtheorem{lemma}{Lemma}[section] 
\newtheorem{proposition}[lemma]{Proposition}
\newtheorem{theorem}[lemma]{Theorem}
\newtheorem{corollary}[lemma]{Corollary}

\newtheorem{prop}[lemma]{Proposition}
\newtheorem{thm}[lemma]{Theorem}
\newtheorem{cor}[lemma]{Corollary}

\newenvironment{remark}[1]{\refstepcounter{lemma}%
\vskip 5pt \par\noindent {\bf #1\ \thelemma .}}{\vskip 5pt \par}

\newenvironment{remark*}[1]{\par \vskip 5pt \noindent 
{\bf #1.}}{\vskip 5pt \par}

\newcommand{\alg}{{\rm Alg}\,}
\newcommand{\Aut}{\operatorname{Aut}}
\newcommand{\Ad}{\operatorname{Ad}}
 
\newcommand{\bh}{\ensuremath{{\mathcal B}({\mathcal H})}}

\newcommand{\cstaralg}{$C^*$-algebra}
\newcommand{\cstar}{\hbox{$C^*$}}
\providecommand{\dual}[1]{\ensuremath{#1^{\#}}}
\providecommand{\ddual}[1]{\ensuremath{#1^{\#\#}}}

\newcommand{\ds}{\displaystyle}
\newcommand{\eps}{\mbox{$\varepsilon$}}

\newcommand{\innerprod}[1]{\left\langle #1\right\rangle}

\newcommand{\indlim}{\varinjlim}
\newcommand{\indlimit}{\varinjlim}

\newcommand{\lat}{{\rm Lat}\,}

\newcommand{\norm}[1]{\left\|{#1}\right\|}

\providecommand{\qed}%
{\hfill \vrule height5pt width4pt depth1pt \vspace{+2.00ex}}

\newcommand{\spn}{\operatorname{span}}

\newcommand{\bbA}{{\mathbb{A}}}

\newcommand{\bbC}{{\mathbb{C}}}

\newcommand{\bbN}{{\mathbb{N}}}

\newcommand{\bbQ}{{\mathbb{Q}}}
\newcommand{\bbR}{{\mathbb{R}}}

\newcommand{\bbT}{{\mathbb{T}}}

\newcommand{\bbZ}{{\mathbb{Z}}}
  \newcommand{\A}{{\mathcal{A}}}
  \providecommand{\sA}{\A}
  \newcommand{\B}{{\mathcal{B}}}
  \newcommand{\C}{{\mathcal{C}}}
  
  \newcommand{\D}{{\mathcal{D}}}
  
  \newcommand{\E}{{\mathcal{E}}}

\renewcommand{\H}{{\mathcal{H}}}
  \newcommand{\I}{{\mathcal{I}}}
  
\renewcommand{\L}{{\mathcal{L}}}
  \newcommand{\M}{{\mathcal{M}}}
  \newcommand{\N}{{\mathcal{N}}}
\renewcommand{\O}{{\mathcal{O}}}

  \newcommand{\R}{{\mathcal{R}}}

  \newcommand{\U}{{\mathcal{U}}}
  
  \newcommand{\W}{{\mathcal{W}}}
  \newcommand{\X}{{\mathcal{X}}}
  
  \newcommand{\Z}{{\mathcal{Z}}}
\newcommand{\fB}{{\mathfrak{B}}}

\newcommand{\fM}{{\mathfrak{M}}}
\newcommand{\fN}{{\mathfrak{N}}}

\newcommand{\fT}{{\mathfrak{T}}}

\providecommand{\sA}{\A}

\theoremstyle{definition}
\newtheorem{defn}[lemma]{Definition}
\theoremstyle{remark}

\providecommand{\Eig}{\E}
\providecommand{\Eigone}{\E^1}

\providecommand{\chull}{\overline{\text{co}}^\sigma}

\providecommand{\cstardiag}{\cstar-diagonal}
\providecommand{\Ad}{\hbox{Ad}\,}
\providecommand{\indlimit}{\indlim}

\def\symbdown#1{\Big\downarrow\rlap{$\vcenter{\hbox{$\scriptstyle #1$}}$}}
\def\symbup#1{\Big\uparrow\rlap{$\vcenter{\hbox{$\scriptstyle #1$}}$}}

\def\symbse#1{\!\!\!\!\searrow\rlap{$\vcenter{\hbox{$\scriptstyle #1$}}$}}
\providecommand{\id}{\operatorname{id}}

\begin{document}
\title[Coordinates and Isomorphisms ]{Coordinate Systems and 
  Bounded Isomorphisms for Triangular Algebras}
\author[A.P. Donsig]{Allan P. Donsig}
\email{adonsig@math.unl.edu}
\author[D.R. Pitts]{David R. Pitts}
\email{dpitts@math.unl.edu}
\address{Dept. of Mathematics\\
University of Nebraska-Lincoln\\ Lincoln, NE\\ 68588-0323}

\date{June 29, 2005}

\maketitle

\begin{abstract}
For a Banach $\D$-bimodule $\M$ over an abelian unital \cstaralg\ $\D$, 
we define $\Eigone(\M)$ as the collection of norm-one eigenvectors for the
dual action of $\D$ on the Banach space dual $\dual{\M}$.
Equip $\Eigone(\M)$ with the weak-$*$ topology. 
We develop general properties of $\Eigone(\M)$.
It is properly viewed as a coordinate system for $\M$ when 
$\M\subseteq \C$, where $\C$ is a
  unital \cstaralg\ containing $\D$ as a regular MASA with the
  extension property; moreover, $\Eigone(\C)$ coincides with Kumjian's twist in
  the context of \cstardiag s.  We identify the \cstar-envelope of a
  subalgebra $\A$ of a \cstardiag\ when $\D\subseteq\A\subseteq \C$.
For triangular subalgebras, each containing the MASA, a bounded isomorphism 
induces an algebraic isomorphism of the coordinate systems 
which can be shown to be continuous in certain cases.  
For subalgebras, each containing the MASA, a bounded isomorphism that maps 
one MASA to the other MASA induces an isomorphism of the coordinate systems.
  We show that the weak operator closure of the image of a
  triangular algebra in an appropriate representation is a CSL algebra
  and that bounded isomorphism of triangular algebras extends to an
  isomorphism of these CSL algebras.  We prove that for 
  triangular algebras in our context, any bounded isomorphism is
  completely bounded.  Our methods simplify and extend various
  known results; for example, isometric isomorphisms of the triangular
  algebras extend to isometric isomorphisms of the \cstar-envelopes,
  and the conditional expectation $E:\C\to\D$ is multiplicative when
  restricted to a triangular subalgebra.  Also, we use our methods to
  prove that the inductive limit of
  \cstardiag s with regular connecting maps is again a \cstardiag. 

\end{abstract}

\vfill\eject
\tableofcontents

\section{Introduction}\label{S:Intro}
This paper presents the results of our study of bounded isomorphisms
of coordinatized (nonselfadjoint) operator algebras.  Isometric
isomorphisms have been extensively studied~(see, for example,
\cite{MR94k:47068,MR95a:46080,PowerClTrSuAFC*Al}) and are quite
natural, as they include restrictions of $*$-isomorphisms.  Isometric
isomorphisms of \cstaralg s preserve adjoints.  Bounded isomorphisms,
in contrast, need not preserve adjoints or map $*$-subalgebras to
$*$-subalgebras.  Nonetheless, we obtain structural
results, most notably, that bounded isomorphisms are completely
bounded and that they factor into diagonal-fixing, spatial, and
isometric parts, analogous to Arveson-Josephson's factorization of
bounded isomorphisms of analytic crossed products.

Coordinates have been used in the categories of $W^*$-algebras and
$C^*$-algebras with $*$-homomorphisms for decades, going back at least
to the work of Feldman and Moore,
\cite{FeldmanMooreErEqReI,FeldmanMooreErEqReII} for von Neumann
algebras and Renault's construction of \cstar-algebras for a wide
range of topological groupoids \cite{MR82h:46075}.  Of particular
interest to us is Kumjian's construction, in~\cite{MR88a:46060}, of a
certain $\bbT$-groupoid which he called a \textit{twist}, which he
showed is a classifying invariant for a \textit{diagonal pair}, a
separable \cstar-algebra with a distinguished MASA satisfying various
properties.  Renault describes twists in terms of the dual
groupoid~\cite{RenaultTwApDuGrC*Al} and this perspective is often
helpful, see for example the work of Thomsen in
~\cite{ThomsenOnFrTrGpC*Al}.  Twists have been used by various
authors, most notably Muhly, Qiu, and Solel, and Muhly and Solel, 
to study varied categories of subalgebras and submodules of groupoid
\cstaralg s with isometric morphisms \cite{MR94i:46075,MR1271693}.

To apply coordinate methods to bounded isomorphisms, we found it
necessary to revisit these coordinate constructions, eliminating as
much as possible the use of adjoints and clarifying the role of the
extension property.  We define coordinates for bimodules over an
abelian \cstar-algebra which are intrinsic to the bimodule structure
and not \textit{a priori} closely tied to the $*$-structure.  These
definitions allow us to simplify and extend some of the structural
results in the literature.  In particular, we obtain a number of
results for algebras containing the abelian \cstaralg: e.g.\
isomorphism of the coordinates is equivalent to diagonal-preserving
isomorphism of the algebras.  This analysis of coordinate systems will
be useful, we expect, in applying coordinate constructions to more
general settings.

Bounded isomorphisms play a role in the study of norm-closed operator
 algebras which is parallel to similarity transforms in the study of
 weakly-closed subalgebras of \bh.  During the mid-1980's and early
 1990's, there was considerable interest in the structural analysis of
 such algebras via their similarity theory.  This was particularly
 successful with the class of nest algebras (see, for example,
 \cite{DavidsonSiCoPeNeAl,DavidsonOrrPittsCoInCeNeAl,%
 LarsonNeAlSiTr,LarsonPittsIdNeAl,OrrStIdNeAl}) and, to a lesser
 degree, the CSL algebras.  Interestingly, by using certain
 faithful representations of \cstar-diagonals, we can employ
 similarity theory for atomic CSL algebras to obtain structural
 results for bounded isomorphisms between triangular subalgebras of
 \cstar-diagonals. 

We turn now to a more detailed outline of the paper.  Throughout the
paper, we consider bimodules over an abelian unital \cstaralg\
$\D$. Our view is that the set of coordinates for such a bimodule $\M$
is the collection $\Eigone(\M)$ of norm one elements of the Banach
space dual $\dual{\M}$ which are eigenvectors for the bimodule action
of $\D$ on $\dual{\M}$.  We use \textit{eigenfunctional} for
such elements, and we topologize $\Eigone(\M)$ using the weak-$*$
topology.  With this structure, we call $\Eigone(\M)$ the coordinate
system for $\M$.  In Section~\ref{eigeninter}, we establish some very
general, but useful, properties.  For example, Theorem~\ref{extension}
shows that when $\M_1$ is a submodule of $\M$, then an element of
$\Eigone(\M_1)$ can be extended (but not necessarily uniquely) to an 
element of $\Eigone(\M)$.  

A key step in minimizing the use of adjoints is replacing normalizers
with intertwiners, that is, elements $m \in \M$ so that $m\D = \D m$.
Section~\ref{S:normal} starts by showing that intertwiners and
normalizers are closely related, at least when $\D$ is a MASA in a
unital \cstaralg\ $\C$ containing $\M$
(Propositions~\ref{interTOnormal} and~\ref{normalTOinter}).  The
extension property is used to construct intertwiners and to slightly
strengthen a key technical result of Kumjian.  In this generality,
$\Eigone(\M)$ need not separate points of $\M$.

In Section~\ref{S:inclusions} we work in the context of $\D$-bimodules
$\M\subseteq \C$, where $\C$ is a unital \cstaralg\ and $\D\subseteq
\C$ is a regular MASA with the extension property.  Such a pair
$(\C,\D)$ we term a \textit{regular \cstar-inclusion}.  Here the
coordinates are better behaved: eigenfunctionals on submodules of $\M$
extend uniquely to $\M$.  Nevertheless, the coordinates for regular
\cstar-inclusions again are not sufficiently rich to separate points.
However, the failure to separate points is intimately related to a
certain ideal of $\fN\subseteq \C$, and Theorem~\ref{Dinclusion} shows
that the quotient of $\C$ by $\fN$ is a \cstar-diagonal, which is a
mild generalization of Kumjian's notion of diagonal pair due to
Renault~\cite{RenaultTwApDuGrC*Al}.  Essentially, a \cstardiag\ is a
regular \cstar-inclusion where the conditional expectation
$E:\C\rightarrow \D$ arising from the extension property is faithful.
There are an abundance of \cstardiag s: crossed products of abelian
\cstar\ algebras by freely acting amenable groups are \cstardiag s,
and Theorem~\ref{indlimit} shows that inductive limits of \cstardiag s
are again \cstardiag s when the connecting maps satisfy a regularity
condition.  In particular, AF-algebras and circle algebras can be
viewed as \cstardiag s.  While our primary interest in this paper is
the use of coordinate methods to study nonselfadjoint algebras,
Theorem~\ref{Dinclusion} and Theorem~\ref{indlimit} are examples of
results in the theory of \cstaralg s obtained using our perspective.

When $\M$ is a $\D$-bimodule contained in a \cstar-diagonal, the
elements of $\Eigone(\M)$ do separate points, and when the bimodule is
an algebra, also have a continuous product.  In
Section~\ref{S:inclusions}, we use the extension property to show that
the coordinate system $\Eigone(\C)$ for a \cstar-diagonal agrees with
Kumjian's twist.  Our methods provide some simplifications and
generalizations of Kumjian's results.  One of the interesting features
of our approach is that it allows us to show that the coordinate
systems for bimodules $\M\subseteq \C$ are intrinsic to the bimodule
alone, and not dependent on the choice of the embedding into the
particular \cstar-diagonal.  This, and the agreement of our
construction with Kumjian's is achieved in Theorems~\ref{partA},
Corollary~\ref{extensionunique}, and Proposition~\ref{loccmpt}.  An
interesting application of our coordinate methods is
Theorem~\ref{envelope}, which shows that if the pair $(\C,\D)$ is a
\cstardiag, and $\A$ is a norm closed algebra with $\D\subseteq
\A\subseteq \C$, then the \cstar-envelope of $\A$ coincides with the
\cstar-subalgebra of $\C$ generated by $\A$.

We study the representation theory of \cstar-diagonals in
Section~\ref{repns}, obtaining a faithful atomic representation compatible 
with the \cstar-diagonal structure.  
Our methods are reminiscent of Gardner's work on isomorphisms of \cstaralg s 
in~\cite{GardnerOnIsC*Al}.
The significance to us of these representations is
that they carry subalgebras containing the diagonal to algebras weakly
dense in a CSL-algebra, see Theorem~\ref{faithfulcompatable}.
This theorem enables us to prove, with fewer hypotheses than previously needed,
that the conditional expectation is multiplicative when restricted to a 
triangular subalgebra, Theorem~\ref{ExpectHomomorph}.

Section~\ref{S:Iso} considers diagonal-preserving bounded
isomorphisms, those that map the diagonal of one algebra onto the diagonal of
the other.
Theorem~\ref{bimodiso} shows that in this case, there is an
isomorphism between coordinate systems arising naturally from the
algebra isomorphism.  Consequently, we are able to prove several
results, such as Theorem~\ref{EigoneUnique}, which shows an
automorphism of a triangular algebra which fixes the diagonal
pointwise arises from a cocycle, and Theorems~\ref{trialgcor}
and~\ref{bigSteve}, which show that coordinates are invariant under 
diagonal-preserving bounded isomorphisms, extending previous results
for isometric isomorphisms.

We then turn to bounded isomorphisms of triangular algebras which do
not preserve the diagonal.  A main result, Theorem~\ref{algisom} shows
that a given bounded isomorphism of triangular subalgebras induces an
\textit{algebraic} isomorphism $\gamma$ of their coordinate systems,
but our methods are not strong enough to show that this isomorphism is
continuous everywhere.  Nevertheless, this does show that the
algebraic structure of coordinate systems for triangular algebras is
invariant under bounded isomorphism, a fact we believe would be
difficult to show using previously existing methods.

While we are not able to prove continuity of the map $\gamma$ on coordinate
systems in general, we can prove it in various special cases.
A bounded isomorphism of triangular algebras induces a canonical
$*$-isomorphism of the diagonals.  If this $*$-isomorphism extends to
a $*$-isomorphism of the \cstar-envelopes, then
Corollary~\ref{envelopeiso} shows that the product of $\gamma$ with an
appropriate cocycle yields a continuous isomorphism.
Also, Theorem~\ref{algpresgamma} shows that bounded
isomorphism of triangular algebras implies an isomorphism of their
coordinate systems when the triangular algebras are generated by their
\textit{algebra-preserving normalizers}.  This is a new class of
algebras which contains a variety of known classes, including those
generated by order-preserving normalizers or those generated by
monotone $G$-sets.

Another of our main results, Theorem~\ref{CSLextend}, shows that if 
boundedly isomorphic triangular subalgebras are represented in the faithful
representation constructed in Theorem~\ref{faithfulcompatable}, then the 
isomorphism extends to an isomorphism of the weak closures of the 
triangular algebras.  
A crucial ingredient in proving this is Theorem~\ref{algisom}.
Theorem~\ref{CSLextend} allows us to
use known results about isomorphisms of CSL algebras to prove another
main result, Theorem~\ref{CB}, which asserts that every bounded
isomorphism of triangular subalgebras of \cstardiag s is completely
bounded.  Consequently, we are able to prove that every isometric
isomorphism of triangular algebras extends to an isometric isomorphism
of the corresponding \cstar-envelopes.
Another application of Theorem~\ref{CSLextend} is Theorem~\ref{isostructure},
which gives the factorization into diagonal-fixing, spatial, and isometric 
parts mentioned earlier.

Together with Power, we asserted \cite[Theorem~4.1]{MR2002k:47148}
that two limit algebras are isomorphic if and only if a certain type
of coordinates for the algebras, namely their \textit{spectra}, are
isomorphic.  Unfortunately, there is a serious gap in the proof, and
another of our motivations for the work in this paper was an attempt
to provide a correct proof.  While we have not yet done this, our
results provide evidence that \cite[Theorem~4.1]{MR2002k:47148} is
true.  The main result of Section~\ref{S:IndLim},
Theorem~\ref{trienvelopeiso}, shows that an (algebraic) isomorphism of
a limit algebra $\A_1$ onto another limit algebra $\A_2$ implies the
existence of a $*$-isomorphism $\tau$ of their \cstar-envelopes.  If
we could choose $\tau$ so that $\tau(\A_1)=\A_2$, then
\cite[Theorem~4.1]{MR2002k:47148} would follow easily, but we do not
know this.  However, any isomorphism of $\A_1$ onto $\A_2$ induces a
$*$-isomorphism $\alpha$ of $\A_1\cap \A_1^*$ onto $\A_2\cap \A_2^*$.
If $\tau$ can be chosen so that $\tau$ extends $\alpha$,
Corollary~\ref{envelopeiso} shows that $\tau$ carries $\A_1$ onto
$\A_2$, and moreover, there is an isomorphism of the corresponding
coordinate systems.  It is somewhat encouraging that
Theorem~\ref{trienvelopeiso} shows that there is no $K$-theoretic
obstruction to the existence of a $\tau$ which extends $\alpha$.  This
is as close as we have been able to come in our efforts to provide a
correct proof of \cite[Theorem~4.1]{MR2002k:47148}.

\section{Intertwiners and Eigenfunctionals}\label{eigeninter}

In this section we provide a very general discussion of coordinates.
Although our focus in this paper is on \cstar-diagonals and regular
\cstar-inclusions, defined in Section~\ref{S:inclusions}, we start
in a more general framework, with a view to extending coordinate
methods beyond our focus here.
Indeed, there are several useful general results, most notably
Theorem~\ref{extension}, which shows that eigenfunctionals can be extended 
from one bimodule to another bimodule containing the first.  
{F}rom this result, we characterize when an eigenfunctional with a given range 
and source exists, using a suitable seminorm.

\begin{remark*}{Notational Convention}
Given a Banach space $X$, we denote its Banach space dual by $\dual{X}$,
to minimize confusion with adjoints.
\end{remark*}

Throughout this section, $\D$ will be a unital, abelian
\cstar-algebra, and $\M$ will be a Banach space which is also a
bounded $\D$-bimodule, that is,
there exists a constant
$K>0$ such that for every $d, f\in\D$ and $m\in \M$, 
$$\norm{d\cdot m\cdot f}\leq K \norm{d}\, \norm{m}\, \norm{f}.$$
As usual, $\dual{\M}$ becomes a Banach $\D$-module with the action,
$$\innerprod{m,f\cdot \phi\cdot d}= \innerprod{d\cdot m\cdot f,\phi}
\qquad d,\, f\in \D, \, m\in \M, \text{ and } \phi\in \dual{\M}.$$

\begin{defn}\label{Intertwiners}
A nonzero element $m\in \M$ is a \textit{$\D$-intertwiner}, or
more simply, an \textit{intertwiner} if
$$m\cdot \D=\D\cdot m.$$ If $m\in\M$ is an intertwiner such that
for every $d\in\D$, $d\cdot m\in\bbC m$,  we call $m$ a \textit{minimal
intertwiner}.  

A minimal intertwiner of  $\M^\#$ will be called a \textit{eigenfunctional}; 
when necessary for clarity, we use \textit{$\D$-eigenfunctional}.  
That is, an eigenfunctional is a nonzero linear functional $\phi : \M \to \bbC$
so that, for all $d \in \D$,
$x \mapsto \phi(dx)$, $x \mapsto \phi(xd)$ are multiples of $\phi$.

Denote the set of all $\D$-eigenfunctionals by $\Eig_\D(\M)$
(or $\Eig(\M)$, if the context is clear).
We equip $\Eig(\M)$ with the relative weak$^*$ topology (i.e.\ the
relative $\sigma(\M^\#,\M)$-topology).

Denote the set of all \textsl{norm-one} $\D$-eigenfunctionals by
$\Eigone_\D(\M)$ or $\Eigone(\M)$.

\end{defn}

Given an eigenfunctional $\phi\in\Eig_\D(\M)$, 
the associativity of the maps $d\in\D\mapsto d\cdot \phi$ 
and $d\in\D\mapsto \phi\cdot d$ yields the existence of unique 
multiplicative linear functionals $s(\phi)$ and $r(\phi)$ on
$\D$ satisfying 
$s(\phi)(d) \phi = d \cdot \phi$ and
$r(\phi)(d) \phi = \phi \cdot d$, that is,
\begin{equation}\label{sr} \phi(xd) = \phi(x) \bigl[ s(\phi)(d) \bigr], \qquad
   \phi(dx) = \bigl[ r(\phi)(d)  \bigr] \phi(x).
\end{equation}
\begin{defn}
We call $s(\phi)$ and $r(\phi)$ the \textit{source} and \textit{range}
of $\phi$ respectively.
\end{defn}

There is a natural action of the nonzero complex numbers on $\Eig(\M)$,
sending $(\lambda,\phi)$ to the functional $m \mapsto \lambda
\phi(m)$; clearly 
 $s(\lambda \phi)= s(\phi)$ and $r(\lambda \phi)= r(\phi)$.

We next record a few basic properties of eigenfunctionals.

\begin{proposition}\label{rangecontinuity}
With the weak$^*$-topology, $\Eig(\M) \cup \{0\}$ is closed.
Further, $r:\Eig(\M)\rightarrow\hat{\D}$ and
$s:\Eig(\M)\rightarrow \hat{\D}$ are continuous.
\end{proposition}

\begin{proof}
Suppose $(\phi_\lambda)_{\lambda\in\Lambda}$ is a net in
$\Eig(\M)\cup\{0\}$ and $\phi_\lambda\stackrel{\text{w}*}{\rightarrow}
\phi\in \dual{\M}$.

If $\phi=0$, there is nothing to do, so we assume that $\phi\neq 0$.
Choose $m\in\M$ such that $\phi(m)\neq 0$.
Then $\phi_\lambda(m)\rightarrow \phi(m)$, so there
exists $\lambda_0\in\Lambda$ such that $\phi_\lambda(m)\neq 0$ for
every $\lambda\succeq \lambda_0$.  For any $d\in\D$, and
$\lambda\succeq \lambda_0$,  we have
$$s(\phi_\lambda)(d)=
\frac{\phi_\lambda(md)}{\phi_\lambda(m)}\rightarrow
\frac{\phi(md)}{\phi(m)}$$
Thus, $s(\phi_\lambda)$ converges weak-$*$ to the functional
$\rho$ in $\hat{\D}$ given by $\rho(d)=\phi(md)/\phi(m)$.
Similarly, $r(\phi_\lambda)$ converges weak-$*$ to $\tau\in\hat{\D}$
given by $\tau(d)=\phi(dm)/\phi(m).$
It now follows that $\phi$ is a eigenfunctional with
$r(\phi)=\tau$ and $s(\phi)=\rho$.
Thus, $\Eig(\M)\cup\{0\}$ is closed.

In particular, if $\phi_\lambda \to \phi \in \Eig(\M)$,
then $s(\phi_\lambda) \to s(\phi)$ and
similarly for the ranges, so $s$ and $r$ are continuous.
\end{proof}

To be useful, there should be many eigenfunctionals.  This need not
occur for arbitrary bimodules, as Example~\ref{noeigenfunctionals}
shows.  However, for bimodules of \cstar-diagonals, which are the
bimodules of principal interest in the present paper,
Proposition~\ref{Mdensity} below will show that eigenfunctionals exist
in abundance.

We need the following seminorm to extend eigenfunctionals and to
characterize when eigenfunctionals exist.

\begin{defn}
For $\sigma, \rho\in\hat{D}$,  define
$B_{\sigma,\rho}:\M\rightarrow \bbR$ by
\[
B_{\sigma,\rho}(m):= \inf\{\norm{dmf}: d,f\in \D, \sigma(d)=\rho(f)=1\}.
\]
\end{defn}

These infima do not increase if we restrict to elements $d$ or $f$ of
norm one.  Indeed, for any elements $d$ and $f$ as above, since
$|\sigma(d)|\leq \norm{d}$ and $|\rho(f)|\leq \norm{f}$, we can
replace them with $d/\norm{d}$ and $f/\norm{f}$ and this will only
decrease the norm of $\norm{dmf}$.  Thus,
\[ B_{\sigma,\rho}(m)=\inf\{\norm{dmf}: d,f\in \D,
    \sigma(d)=\rho(f)=1=\norm{d}=\norm{f}\}. \]
In particular, we have $B_{\sigma,\rho}(m)\leq \norm{m}$.

A variant of this seminorm was used by Steve Power in~\cite{MR93a:46114}
to distinguish families of limit algebras associated to singular MASAs.

\begin{proposition}\label{seminormprop}
For $\sigma,\rho\in\hat{\D}$, we have the following:
\begin{enumerate}
\item $B_{\sigma,\rho}$ is a seminorm.
\item For $m\in\M$ and $\phi \in \Eig(\M)$,
$|\phi(m)| \le \|\phi\| B_{r(\phi),s(\phi)}(m).$
\item For $m\in\M$, $d,f \in \D$, $B_{\sigma,\rho}(dmf)
    = |\sigma(d)| B_{\sigma,\rho}(m) |\rho(f)|$.
\item If $f \in \dual{\M}\setminus \{0\}$ satisfies $|f(m)|\leq
B_{\sigma,\rho}(m)$ for all $m\in\M$, then $f\in \Eig(\M)$ with
$s(f)=\rho$, $r(f)=\sigma$, and $\norm{f}\leq 1$.
\end{enumerate}
\end{proposition}

\begin{proof}
For (1), it is immediate that
$B_{\sigma,\rho}(\lambda m)=|\lambda| B_{\sigma,\rho}(m)$,
for $\lambda\in\bbC$ and $m\in\M$.

To show subadditivity, let $a, b\in\M$ and choose $\eps >0$.
Pick norm one elements $d_1$, $d_2$, $f_1$, $f_2$ of $\D$,
satisfying $\sigma(d_1)=\sigma(d_2)=1=\rho(f_1)=\rho(f_2)$,
 $\norm{d_1af_1}<B_{\sigma,\rho}(a)+\eps$ and
$\norm{d_2bf_2}< B_{\sigma,\rho}(b)+\eps$.
Then
\begin{align*}
\norm{d_1d_2(a+b)f_1f_2}&\leq \norm{d_2d_1af_1f_2}+\norm{d_1d_2bf_2f_1}\\
& \leq \norm{d_1af_1}+\norm{d_2bf_2}\\
& < B_{\sigma,\rho}(a)+B_{\sigma,\rho}(b) +2\eps,
\end{align*}
whence $B_{\sigma,\rho}(a+b)\leq B_{\sigma,\rho}(a)+B_{\sigma,\rho}(b)$.

For (2), suppose $d, f\in\D$ with $s(\phi)(d)=r(\phi)(f)=1$.
Then $|\phi(m)|=|\phi(dmf)|\leq \norm{\phi}\, \norm{dmf}$.
Taking the infimum over all such $d$ and $f$ gives the inequality.

For (3), we show that
$B_{\sigma,\rho}(dm) = |\sigma(d)| B_{\sigma,\rho}(m)$;
the proof for $\rho$ and $f$ is similar.

For any $a,b\in\D$ with $\norm{a}=\norm{b}=1=\sigma(a)=\rho(b)$ we
have,
\begin{align*}\norm{admb}&\leq
  \norm{\sigma(d)amb}+\norm{(ad-\sigma(d)a)mb}\\
&\leq |\sigma(d)|\norm{amb}+\norm{ad-\sigma(d)a}\,\norm{m}.
\end{align*}
Hence $B_{\sigma,\rho}(dm)\leq
|\sigma(d)|\norm{amb}+\norm{ad-\sigma(d)a}\,\norm{m}.$
The definition of $B_{\sigma,\rho}$ and the fact that 
 $\inf\{\norm{ad-\sigma(d)a}: a\in\D,\,  \norm{a}=1=\sigma(d)\}=0$, imply that
given $\eps>0$, we may find norm one elements $a_1,b_1$ and
$a_2$ of $\D$ such that $\sigma(a_1)=\sigma(a_2)=\rho(b_1)=1$ and 
$$\norm{a_1mb_1}<B_{\sigma,\rho}(m)+\eps\quad \text{and}\quad
\norm{a_2d-\sigma(d)a_2}\, \norm{m} <\eps.$$
Then
\begin{align*}
B_{\sigma,\rho}(dm)&\leq
|\sigma(d)|\norm{a_1a_2mb_1}+\norm{a_1a_2d-\sigma(d)a_1a_2}\,
\norm{m}\\
&\leq |\sigma(d)|\norm{a_1mb_1}+\norm{a_2d-\sigma(d)a_2}\,\norm{m}\\
&\leq |\sigma(d)|(B_{\sigma,\rho}(m)+\eps) +\eps,
\end{align*} whence $B_{\sigma,\rho}(dm)\leq
|\sigma(d)|B_{\sigma,\rho}(m).$

To obtain $|\sigma(d)|B_{\sigma,\rho}(m)\leq B_{\sigma,\rho}(dm)$,
observe that for $a,b\in \D$ and
$\norm{a}=\norm{b}=1=\sigma(a)=\rho(b)$, 
$$|\sigma(d)|B_{\sigma,\rho}(m)\leq\norm{\sigma(d)amb}
\leq \norm{\sigma(d)a-ad}\, \norm{m}
+\norm{admb},$$ and argue as above.

Finally, suppose $f$ is a nonzero linear functional on $\M$ satisfying
 $|f(m)|\leq B_{\sigma,\rho}(m)$ for every $m\in\M$.  
 As $B_{\sigma,\rho}(x)\leq \norm{x}$, we see that
 $f$ is bounded and $\norm{f}\leq 1$.

Suppose $d\in\D$, and let $k=d-\rho(d)I$.  Clearly $\rho(k)=0$, and
so, by for $x\in \M$,
$|f(xk)|\leq B_{\rho,\sigma}(xk)=0$.  Therefore
$f(xd)=f(x)\rho(d)+ f(xk)=f(x)\rho(d)$.  Similarly,
$f(dx)=\sigma(d)f(x)$.  Thus, $f$ is an eigenfunctional with range
$\sigma$ and source $\rho$.
\end{proof}

\begin{theorem}\label{extension}
Suppose $\M$ is a  norm-closed $\D$-bimodule and $\N\subseteq \M$ is a
norm-closed sub-bimodule.
Given $\phi \in \Eig(\N)$, there is $\psi \in \Eig(\M)$
with $\psi |_{\N} = \phi$ and $\norm{\phi}=\norm{\psi}$.
\end{theorem}

Necessarily, $\psi$ has the same range and source as $\phi$.

\begin{proof}
Let $\rho=s(\phi)$ and $\sigma=r(\phi)$.  From
 Proposition~\ref{seminormprop}~(2), $|\phi(n)| \le
 \|\phi\|B_{\sigma,\rho}(n)$ for all $n \in \N$.  By the Hahn-Banach
 Theorem, there exists an extension of $\phi$ to a linear functional
 $\psi$ on $\M$ satisfying $|\psi(x)|\leq \norm{\phi}B_{\sigma,\rho}(x)$ for
 all $x\in\M$.  Now apply the last part of
 Proposition~\ref{seminormprop}. 
\end{proof}

We would like to be able to say that the extension in
Theorem~\ref{extension} is unique, but this need not be true.  For
example, if $\M_1=\D$ and $\M_2=\C$, then we are considering
extensions of pure states, which need not be unique
(see~\cite{ArchboldBunceGregsonExStC*AlII,MR0123922}, for example).
However, in the context of regular \cstar-inclusions the extension is
unique, as we show in Section~\ref{S:inclusions}.

We now characterize the existence of eigenfunctionals in terms of the
$B_{\sigma, \rho}$ seminorms. 

\begin{thm}\label{nonzeroeig}
Suppose $\sigma, \rho\in\hat{\D}$.
There is $\phi \in \Eig(\M)$ with $r(\phi)=\sigma$ and $s(\phi)=\rho$
if and only if $B_{\sigma,\rho} |_{\M} \neq 0.$
\end{thm}

\begin{proof}
If $\phi \in \Eig(\M)$ with $s(\phi)=\rho$ and $r(\phi)=\sigma$, then 
Proposition~\ref{seminormprop}~(2) implies $B_{\sigma,\rho} |_{\M} \neq 0$.

Conversely, suppose $B_{\sigma,\rho}(m)\neq 0$.  Define a linear
functional $f$ on $\bbC m$ by $f(\lambda m )=\lambda
B_{\sigma,\rho}(m)$.  Now use the Hahn-Banach Theorem to extend 
 $f$ to a linear functional on all of $\M$ which satisfies $|f(m)|\leq
B_{\sigma,\rho}(m)$, and apply Proposition~\ref{seminormprop}.
\end{proof}

\begin{remark}{Example}\label{noeigenfunctionals}
Here is an example  where $\Eig(\M)$ is empty. 
In $\B(L^2[0,1])$, let $\D$ be the operators of multiplication
by elements of $C[0,1]$, and let $\M$ be the compact operators.  If
$\sigma\in [0,1]$ and $\Lambda=\{d\in\D: \hat{d}(\sigma)=1 \text{ and
} 0\leq d\leq I\}$, then $\Lambda$ becomes a directed set with the direction
$d\preceq e$ if and only if $d-e\geq 0$. Viewing $\Lambda$ as a net
indexed by itself, then $\Lambda$ is a bounded net converging
 strongly to zero.
Hence given any compact operator
$K$, the net $\{dK\}_{d\in \Lambda}$  converges
to zero in norm.  It follows that $B_{\sigma,\rho}(K)=0$ for all
$\sigma,\rho\in\hat{\D}$, so by Theorem~\ref{nonzeroeig}, the set of
eigenfunctionals is $\{0\}$.
\end{remark}

Not surprisingly, eigenfunctionals behave appropriately under bimodule maps.

For $i=1,2$, let $\M_i$ be $\D$-bimodules and let
$\theta:\M_1\rightarrow \M_2$ be a bounded $\D$-bimodule map.
Recall the Banach adjoint map\label{adjointdefn}
$\theta^\# : \M_2^\# \to \M_1^\#$,
given by $\theta^\#(\phi) = \phi \circ \theta$.
If  $\phi \in \Eig(\M_2)$, then $\theta^\# \phi \in \Eig(\M_1)$, and 
we have $s(\phi \circ \theta) = s(\phi)$,
 and $r(\phi \circ \theta) = r(\phi)$.

We include the following simple result for reference purposes; the
proof is left to the reader.

\begin{prop}\label{adjointmap}
For $i=1,2$, let $\M_i$ be $\D$-bimodules and suppose $\theta : \M_1
\to \M_2$ is a bounded linear map which is also a
$\D$-bimodule map.  Then $\theta^\# |_{\Eig(\M_2)}$ is a
continuous map of $\Eig(\M_2)$ into $\Eig(\M_1)\cup\{0\}$.  

If $\theta$ is bijective, then $\theta^\#|_{\Eig(\M_2)}$ is
a homeomorphism of $\Eig(\M_2)$ onto $\Eig(\M_1)$.
If $\theta$ is isometric, $\theta^\#|_{\Eigone(\M_2)}$ is
a homeomorphism of $\Eigone(\M_2)$ onto $\Eigone(\M_1)$.
\end{prop}

\section{Normalizers and the Extension Property}
\label{S:normal}

We begin this section with a discussion of the relationship between
normalizers and intertwiners.  Together,
Propositions~\ref{interTOnormal} and~\ref{normalTOinter} show that
all intertwiners of a maximal abelian \cstaralg\ are normalizers
and every normalizer can be approximated as closely as desired by
intertwiners.  This shows that for our purposes, there is no
disadvantage in using intertwiners instead of normalizers; moreover,
the fact that intertwiners behave well under bounded isomorphism is a
considerable advantage.

In the second part of the section, we consider a \cstaralg\ $\C$
containing a MASA $\D$ which has the extension property (see
Definition~\ref{D:extension}) and a $\D$-bimodule $\M\subseteq \C$.
We use the extension property and a Theorem from
~\cite{ArchboldBunceGregsonExStC*AlII} to construct intertwiners in
$\M$ from intertwiners in $\C$,  and strengthen a key technical
result of Kumjian, Proposition~\ref{Kumjian}.  We see these results as
a step towards extending coordinates from \cstar-diagonals to more
general settings.

\begin{remark*}{Context for Section~\ref{S:normal}} 
Throughout this section, $\C$ will be a unital \cstar-algebra and 
$\D\subseteq\C$ will be an abelian $\cstar$ subalgebra containing the 
unit of $\C$.  
Bimodules considered in this section will be closed subspaces $\M$ of $\C$ 
which are $\D$-bimodules under multiplication in $\C$.
\end{remark*}

\begin{defn}
An element $v\in\M$ is a \textit{normalizer} of $\D$ if $v\D
v^*\cup v^*\D v\subseteq \D$.  
The set of all such elements is denoted $\N_{\D}(\M)$ or,
if $\D$ is clear, $\N(\M)$. 
Recall that $\M$ is said to be \textit{regular} if 
the closed span of $\N_{\D}(\M)$ equals $\M$.  
\end{defn}

Typically, normalizers play a major role in constructing coordinates
for operator algebras.
Since normalizers depend on the involution, it can be difficult to
determine if isomorphisms that are not $*$-extendible or isometric 
preserve normalizers.
Intertwiners (Definition~\ref{Intertwiners}) are not defined in terms
of the involution and so it can be easier to decide if they are
preserved by such isomorphisms.
We begin with a comparison of normalizers and intertwiners.

It is easy to find examples of intertwiners of abelian \cstaralg s
which are not normalizers:  for a simple example, observe
that every element of $M_2(\bbC)$ is an intertwiner for $\bbC I_2$.
However, the next theorem shows that when the abelian algebra is a MASA, 
intertwiners are normalizers.

\begin{proposition}\label{interTOnormal}
If $v\in\C$ is an intertwiner for $\D$, then $v^*v,\, vv^* \in \D'\cap
\C$.  If $\D$ is maximal abelian in $\C$, then $v$ is a normalizer of
$\D$.
\end{proposition}

\begin{proof}
Let $v$ be an intertwiner.  
Let $J_s(v):=\{d\in\D: vd=0\}$ and $J_r(v):=\{d\in\D: dv=0\}$.  
Then $J_s$ and $J_r$ are norm-closed ideals in $\D$.  
Define a  mapping $\alpha_v$ between $\D/J_s$ and
$\D/J_r$  by $\alpha_v(d+J_s)=d'+J_r$, where $d'\in\D$ is chosen so
that $vd=d'v$.  
It is easy to check that $\alpha_v$ is a well-defined
$*$-isomorphism of $\D/J_s$ onto $\D/J_r$.

Let $d=d^* \in \D$.
Then $\alpha_v(d+J_s)=d'+J_r$ where $d'$ is chosen so that $d'=d'^* \in \D$.
Thus, we have the equality of sets,
\[ \begin{split}
\{v^*vd\}=v^*v(d+J_s) & =  v^*(d'+J_r)v =[(d'+J_r)v]^*v \\
	 & =[v(d+J_s)]^*v=(d+J_s)v^*v=\{dv^*v\}. 
\end{split} \]
Hence $v^*v$ commutes with the selfadjoint elements of $\D$ and so 
commutes with $\D$.
Since $v^*$ is also an intertwiner, we conclude similarly that $vv^*\in\D'$.  

If $\D$ is maximal abelian, then $v\D v^*=\D vv^*\subseteq \D$ and 
$v^*\D v= v^*v\D\subseteq \D$, and $v$ is a normalizer as desired.
\end{proof}

For $v \in \N(\C)$, let $S(v):=
\{\phi\in\hat{\D}: \phi(v^*v)> 0\}$; note this is an open set in $\hat{\D}$.
As observed by Kumjian (see~\cite[Proposition~6]{MR88a:46060}), 
there is a homeomorphism $\beta_v:S(v)\rightarrow S(v^*)$ given by
$$\beta_v(\phi)(d)=\frac{\phi(v^*dv)}{\phi(v^*v)}.$$
It is easy to show that $\beta_v^{-1}=\beta_{v^*}$.

\begin{proposition}\label{normalTOinter}
For $v\in \N(\C)$, if $\beta_{v^*}$ extends to a homeomorphism of
$\overline{S(v^*)}$ onto $\overline{S(v)}$, then $v$ is an intertwiner.  
Moreover, if $\I:=\{w\in\C: w\D=\D w\}$ is the set of intertwiners,
then $\N(\C)$  is contained in the norm-closure of $\I$, 
and when $\D$ is a MASA in $\C$,
$\N(\C)=\overline{\I}$. 
\end{proposition}

\begin{proof}
It is clear that the set of normalizers is norm-closed.

Regard $\C$ as sitting inside its double dual $\ddual{\C}$.  
Let $v=u|v|=|v^*|u$ be the polar decomposition for $v$.  
Since $u$ is the strong-$*$ limit of $u_n:=v(1/n+|v|)^{-1}$, we find 
that $u$ also normalizes $\ddual{\D}$.

Therefore, given any $d\in\D$, $vdv^*=u|v|du^*|v^*|=udu^*vv^*$.
Hence for $\phi\in S(v^*)$, we have
\begin{equation}\label{ualpha}
\beta_{v^*}(\phi)(d)= \phi(udu^*).
\end{equation}

Suppose now that $\beta_{v^*}$ extends to a homeomorphism of
$\overline{S(v^*)}$ onto $\overline{S(v)}$.
By Tietze's Extension Theorem, we may then find
an element $d_1\in\D$ such that for every $\phi\in S(v^*)$,
$\phi(d_1)=\beta_{v^*}(\phi)(d).$
Thus, for every $\phi\in\hat{\D}$,
$\phi(vdv^*)=\phi(udu^*)\phi(vv^*)=\phi(d_1vv^*),$ so that
$$(udu^*-d_1)vv^*(ud^*u^*-d_1^*)=0.$$  
This shows that $udu^*v=d_1v$ and so
$$vd=u|v|d=uu^*ud|v|=udu^*u|v|=udu^*v=d_1v. $$ 
Hence $v\D\subseteq \D v$.  
Since the adjoint of a normalizer is again a normalizer and
$\beta_v=(\beta_{v^*})^{-1}$, we may
repeat this argument to obtain $v^*\D \subseteq \D v^*$.  
Taking adjoints yields $\D v\subseteq v\D$.  
Hence $v$ is an intertwiner.

Given a general normalizer $v$, let $\eps>0$ and let 
$K=\{\phi\in\hat{\D}: \phi(vv^*)\geq \eps^2\}$. 
Then $K$ is a compact subset of $S(v^*)$.  
Choose $d_0\in\D$ so that $0\leq d_0\leq I$ and $\hat{d_0}$ is 
compactly supported in $S(v^*)$ and $\hat{d_0}|_K=1$.  
Since $\beta_{v^*d_0}=\beta_{v^*}|_{S(v^*d_0)}$ and $\beta_{v^*}$
is a homeomorphism, it extends to a homeomorphism of 
$\overline{S(v^*d_0)}$ onto $\overline{S(d_0v)}$.
Thus $d_0v$ is an intertwiner.
Further, $\norm{d_0v-v}=\norm{(d_0-1)vv^*(d_0-1)}^{1/2} <\eps$
so $d_0v$ approximates $v$ to within $\eps$.
Thus, $\N(\C)\subseteq \overline{\I}$.  When $\D$ is a MASA in $\C$,
Proposition~\ref{interTOnormal} shows every intertwiner of $\D$ is a
normalizer, so $\N(\C)=\overline{\I}$.
\end{proof}

\begin{remark}{Remarks}\label{nostar} 
(1) Taken together, these two propositions show that for a MASA $\D$
in a \cstar-algebra $\C$, a \textit{partial isometry} $v$ is a normalizer if and 
only if it is an intertwiner.  
Related results for partial isometries are known~\cite[Lemma~3.2]{MR1054818}.

\noindent
(2) The $*$-isomorphism $\alpha_v:\D/J_s\rightarrow \D/J_r$ appearing 
in the proof of Proposition~\ref{interTOnormal} induces a homeomorphism from the
zero set $Z_r:=\{\rho\in\hat{\D}: \rho|_{J_r}=0\}$ onto the zero set
$Z_s:=\{\rho\in\hat{\D}:\rho|_{J_s}=0\}$.  It is not hard to show
that when $v$ is an intertwiner and $\D$ is a MASA, $Z_s=\overline{S(v)}$,
$Z_r=\overline{S(v^*)}$ and that $\alpha_v$ is the extension of
$\beta_{v^*}=\beta_v^{-1}$ to $Z_r$.  Thus, it is possible to describe 
$\beta_v$ without explicit reference to the $*$-structure.
\end{remark}

We turn to constructing intertwiners in a module using
the following property.

\begin{defn}\label{D:extension}
Let $\C$ be a unital \cstaralg. 
A $\cstar$-subalgebra $\D\subseteq\C$ is said to have
the \textit{extension property} if 
 every pure state of $\D$ has a unique extension to a state on $\C$
and no pure state of $\C$ annihilates $\D$.
\end{defn}

If $\D \subseteq \C$ is abelian and $\D$ has the extension property
relative to $\C$, then the Stone-Weierstrass Theorem implies that
$\D$ is a MASA~\cite[p.~385]{MR0123922}.
We shall make essential use of the following result characterizing the
 extension property for abelian algebras.
\begin{theorem}[{\cite[Corollary~2.7]{ArchboldBunceGregsonExStC*AlII}}]
\label{exptHull}
Let $\C$ be a unital \cstaralg\ and let $\D$ be an abelian
\cstar-subalgebra of $\C$ which contains the unit of $\C$.  
Then $\D$ has the extension property if and only if 
$$\overline{\rm{co}}\{gxg^{-1}: g\in\D \text{ and $g$ is
  unitary}\}\cap \D\neq \emptyset.$$  Furthermore, when this occurs,
  $\D$ is a MASA and there  
exists a conditional expectation $E:\C\rightarrow \D$
such that 
\begin{equation*}
\overline{\rm{co}}\{gxg^{-1}: g\in\D \text{ and $g$ is
  unitary}\}\cap \D = \{E(x)\};
\end{equation*}

\end{theorem}

\begin{remark*}{Remark}  
Archibold, Bunce and Gregson in~\cite{ArchboldBunceGregsonExStC*AlII} also
show that the condition $$\C=\D+\overline{\rm{span}}\{cd-dc: c\in \C,\, d\in\D\}$$
also characterizes the extension property for an abelian subalgebra
$\D\subseteq \C$.  This characterization was important in Kumjian's
work on \cstardiag s.
\end{remark*}

\begin{defn}
For $v\in \N(\C)$, define $E_v:\C\rightarrow \N(\C)$ by 
\[ E_v(x)=vE(v^*x). \]
\end{defn}

Our first application of Theorem~\ref{exptHull} is essentially contained in
the proof of \cite[Proposition~4.4]{MR94i:46075}, but Theorem~\ref{exptHull}
provides a different (and simpler) proof.  The crucial implication of
Proposition~\ref{modmap} is that bimodules contain many normalizers.

\begin{proposition}\label{modmap}  Suppose $\D$ is a abelian
  $C^*$-subalgebra of the unital \cstaralg\ $\C$ which has the
  extension property, and $\M\subseteq \C$ is a $\D$-bimodule.  If $v
  \in\C$ is a normalizer (resp.\ intertwiner) and $x\in \M$, then
  $E_v(x)\in \M$ and is a normalizer (resp.\ intertwiner).
\end{proposition}

\begin{proof}
Fix $x\in \M$ and let $G$ be the unitary group of $\D$.
Since $v\in N_\D(\C)$,  for every $g\in G$,  $vgv^*\in \D$,
so that $(vgv^*) x g^{-1}\in\M$.
Thus the norm-closed convex hull,
$$H:=\overline{\text{co}}\{(vgv^*)xg^{-1}:g\in G\}\subseteq \M.$$
By Theorem~\ref{exptHull}, $E(v^*x)$
belongs to $K:=\overline{\text{co}}\{gv^*xg^{-1}: g\in G\}$.
Since $vK\subseteq H$, we conclude that $vE(v^*x)=E_v(x)\in \M$.
\end{proof}

Although we have no application for it, the following result, which is
a corollary of Proposition~\ref{modmap},
provides another means of constructing normalizers in a bimodule.

\begin{proposition}\label{less}
If $v\in\N(\C)$, $m\in\M$ and $|\phi(v)|\leq |\phi(m)|$ for all
$\phi\in\Eigone(\C)$, then $v\in\M$.
\end{proposition}

\begin{proof}
For $\rho\in\hat{\D}$ with $\rho(v^*v)\neq 0$, $\rho(v^*v)$ equals
\[ \rho(v^*v)^{1/2} [v,\rho](v) \leq |\rho(v^*v)^{1/2} [v,\rho](m)|
 =|\rho(v^*m)| = |\rho(E(v^*m))|, \] and so, for $\rho\in\hat{\D}$,
$\rho(v^*v)\leq |\rho(E(v^*m))|$.  Thus, there exists $d\in\D$ with
$v^*v=E(v^*m)d$.  Hence for $n\in\bbN$, $(v^*v)^{1/n}$ belongs to the
closed ideal of $\D$ generated by $E(v^*m)$.  By
Proposition~\ref{modmap}, $vE(v^*m)\in\M$ and we conclude
$v(v^*v)^{1/n}\in\M$.  But
$v=\lim\limits_{n\rightarrow\infty}v(v^*v)^{1/n}$, whence $v\in\M$.
\end{proof}

Our second application of Theorem~\ref{exptHull} is to provide an alternate
proof of, and slightly strengthen, a result of Kumjian.  

\begin{proposition}[{\cite[Lemma~9, p.~972]{MR88a:46060}}]\label{Kumjian}
Let $\D$ be an abelian $C^*$-subalgebra of the unital \cstaralg\ $\C$
with the extension property.  
For $v \in \N(\C)$, $v^*E(v)$ and $vE(v^*)$ both belong to $\D$.
If $\rho \in S(v)$, then the following are equivalent:
\begin{enumerate}
\item\label{KU1} $\rho(v) \ne 0$,
\item\label{KU2} $\beta_v(\rho)=\rho$,
\item\label{KU3} $\rho(E(v^*)) \ne 0$.
\end{enumerate}
\end{proposition}

\begin{proof}
Taking $x=I$ and $\M=\D$ in Proposition~\ref{modmap} we see that
$v^*E(v)$ and $vE(v^*)$ both belong to $\D$.

Suppose $\rho(v)\neq 0$.  An easy calculation shows that when $d\in\D$
and $\rho(d)\neq 0$, then $\beta_{vd}=\beta_v$.  Let $d=v^*E(v)$.  By
hypothesis, we find $\rho(d)=|\rho(v)|^2\neq 0$, and another
calculation shows that $\beta_{vd}=\beta_v$, so \eqref{KU1} implies
\eqref{KU2}.

Assume \eqref{KU2} holds.  Letting $G$ again be the unitary group of
$\D$, we have, for all $g\in G$, 
$$\rho(v^*gvg^{-1})=
\rho(v^*gv)\rho(g^{-1})=\beta_v(\rho)(g)\rho(g^{-1})\rho(v^*v)=\rho(v^*v).$$
Thus, $\rho(v^*\overline{\text{co}}\{gvg^{-1} \: g\in
G\})\subseteq\{\rho(v^*v)\}$, and therefore by
Theorem~\ref{exptHull} we obtain $\rho(v^*)\rho(v)=
\rho(v^*E(v))\neq 0$.  Thus $\rho(E(v^*))\neq 0$ and by our earlier
remarks, $vE(v^*)\in\D$. 

Finally, it is evident that \eqref{KU3} implies \eqref{KU1}.
\end{proof}

\section{Regular $C^*$-inclusions and \cstardiag s}
\label{S:inclusions}

We now turn to the key context for our subsequent work, that of
\cstar-diagonals and regular \cstar-inclusions.  After recalling
Kumjian's twist, we show that the elements of the twist are
eigenfunctionals on the \cstar-algebra and conversely
(Proposition~\ref{EigTwist} and Theorem~\ref{Eigequal}).  This allows
us to show that every regular \cstar-inclusion has a quotient which is
a \cstardiag\ with the same coordinate system,
Theorem~\ref{Dinclusion}, as well as strengthening various results
from Section~\ref{eigeninter} in this context.  The crucial result for
subsequent sections is Theorem~\ref{partA}, which shows that for a
bimodule $\M$, the intrinsically defined eigenfunctionals on $\M$ and
the restriction of the twist are the same.

\begin{defn}\label{cstardiag}
The pair $(\C,\D)$ will be called a \textit{regular \cstar -inclusion}
if $\D$ is a maximal abelian \cstar -subalgebra of the unital
\cstaralg\ $\C$ such that
\begin{enumerate}
\item $\D$ has the extension property in $\C$;
\item $\C$ is regular (as a $\D$-bimodule).
\end{enumerate}
Always, $E$ denotes the (unique) conditional expectation of $\C$ onto $\D$.
We call $(\C,\D)$ a \textit{\cstardiag} if, in addition, 
\begin{enumerate}\setcounter{enumi}{2}
\item  $E$ is faithful.
\end{enumerate} 
\end{defn}

\begin{remark}{Remarks} 
\begin{enumerate}
\item 
Renault~\cite{MR87a:46111} gives, without proof, an example of an
algebra satisfying only the first two conditions but not the third,
i.e., a regular \cstar-inclusion that is not a \cstardiag.
\item By Propositions~\ref{interTOnormal} and~\ref{normalTOinter},
regularity of a bimodule $\M$ is equivalent to norm-density of 
the $\D$-intertwiners.
\item As observed by Renault~\cite{MR87a:46111}, this definition 
of \cstardiag\ is equivalent to Kumjian's original definition,
namely,
\begin{enumerate}
\item[(a)] there is a faithful conditional expectation
$E:\C\rightarrow \D$
\item[(b)] the closed span of the free normalizers in $\C$ is $\ker E$.
\end{enumerate}  
A normalizer $v$ of $\D$ is \textit{free} if $v^2=0$.
Kumjian also required that $\C$ is separable and $\hat{\D}$ is second 
countable, but \cite{MR87a:46111} shows this is not necessary.
It is often easier to verify Kumjian's axioms when
working with particular examples.
\end{enumerate}

\end{remark}

Unless explicitly stated otherwise, for the remainder of this section,
we work in the following setting.

\begin{remark}{Context}
Let $(\C,\D)$ be a regular \cstar-inclusion
and $\M \subseteq \C$ be a norm-closed $\D$-bimodule, where the module
action is multiplication. 
That is, for $d,e\in \D$ and $m\in \M$, $d\cdot m\cdot e :=dme \in \M$.
\end{remark}

For $\rho\in\hat{\D}$, the unique extension on $\C$ is $\rho\circ E$,
which we will again denote by $\rho$.
Thus, we regard $\rho$ as either a multiplicative linear functional on $\D$
or a pure state on $\C$ which satisfies $\rho\circ E=\rho$.
In particular, for $d\in\D$ and $x\in\C$, $\rho(dx)=\rho(xd)=\rho(d)\rho(x)$.

We now summarize some results and definitions from~\cite{MR88a:46060}
and then relate eigenfunctionals to the elements of Kumjian's twist.
Note that these results from~\cite{MR88a:46060} hold for regular
$C^*$-inclusions.

\begin{defn}[Kumjian]\label{twist}
For $v \in \N(\C)$ and $\rho \in \hat{\D}$ with $\rho(v^*v) > 0$,
define a linear functional on $\C$, $[v,\rho]$, by
\[ [v,\rho](x) = \frac {\rho(v^*x)}{\rho(v^*v)^{1/2}}. \]
Kumjian denotes by $\Gamma$ the collection of all such linear
functionals.  We shall see in Proposition~\ref{EigTwist} below shows
that $\Gamma=\Eigone(\C)$.
\end{defn}

We follow Kumjian~\cite[p.~982]{MR88a:46060} in pointing out that
Proposition~\ref{Kumjian} implies

\begin{cor}[Kumjian]\label{KumjianCor}
The following are equivalent:
\begin{enumerate}
\item $[v,\rho]=[w,\rho]$,
\item $\rho(v^*w) > 0$,
\item there are $d,e\in \D$ with $\rho(d),\rho(e) > 0$ so that $vd=we$.
\end{enumerate}
\end{cor}

\begin{remark}{Remark}  Using Theorem~\ref{exptHull} and the
  techniques of the proof of Proposition~\ref{modmap}, one can show
  that when $\rho(v^*w)>0$, we may take $d=w^*wE(v^*w)$ and
  $e=w^*vE(v^*w)$ in the third part of Corollary~\ref{KumjianCor}.
\end{remark}

Kumjian shows that $\Gamma$, with a suitable operation and the
relative weak$^*$-topology, is a groupoid and admits a natural
$\bbT$-action, given by $\lambda [v,\rho]=[\overline{\lambda} v,
\rho]$.  The range map for
$\Gamma$ is $[v,\rho]\mapsto \beta_v(\rho)$ and the source map for
$\Gamma$ is $[v,\rho]\mapsto \rho.$   
The map $[v,\rho]\in\Gamma \mapsto (\beta_v(\rho), \rho)$ sends $\Gamma$ to a
Hausdorff equivalence relation (principal groupoid) on $\hat{\D}$,
denoted $\Gamma\backslash \bbT$, and $\Gamma$ is a 
locally trivial principal $\bbT$-bundle over $\Gamma\backslash \bbT$.

While the range and source maps of
this groupoid have, a priori, no connection to the range and source
maps for eigenfunctionals, they turn out to be the same.

We recall the multiplication on $\Gamma$ and use it to define
multiplication of eigenfunctionals for regular \cstar-inclusions.  The
partially defined multiplication on $\Gamma$ has $[v_1,\rho_1],[v_2,\rho_2]$
composable if $\rho_1=\beta_{v_2}(\rho_2)$, in which case,
\[ [v_1,\phi_1]\cdot [v_2,\rho_2] = [v_1v_2, \rho_2]. \]

\begin{defn}[Kumjian]\label{abstracttwist}
A twist is a proper $\bbT$-groupoid $\Gamma$ so that $\Gamma/\bbT$ is 
an r-discrete principal groupoid.
\end{defn}

Kumjian constructs, for each twist $\Gamma$, a \cstar-diagonal,
$(A(\Gamma),B(\Gamma))$.
The main result of~\cite{MR88a:46060} is that for $(\C,\D)$ a \cstar-diagonal, 
there is a unique (up to isomorphism) twist $\Gamma$ and an isomorphism 
$\Phi : A(\Gamma) \to \C$ such that $\Phi(B(\Gamma))=\D$.
Thus, every twist arises as a $\Gamma$ as above, for some \cstar-diagonal $(\C,\D)$.

 From our point of view, justified by the following proposition,
the twist $\Gamma$ associated to $(\C,\D)$ is $\Eigone(\C)$ equipped 
with this groupoid operation, $\bbT$-action, and topology.

\begin{prop}\label{EigTwist}
For all $\rho \in \hat{\D}$ and $v \in \N(\C)$ with $\rho(v^*v) > 0$,
$[v,\rho] \in \Eigone(\C)$.
Moreover, the range and source maps agree, that is, viewed as an
element of $\Eigone(\C)$, we have $\rho=s([v,\rho])$ and
$\beta_v(\rho)=r([v,\rho])$.

Conversely, if $\phi \in \Eig(\C)$ there exists a normalizer $v\in\C$ 
 such that $\phi(v)\neq 0$, and $v$ may be taken to be an intertwiner
 if desired.  For any normalizer (or intertwiner) $v$ with
 $\phi(v)\neq 0$, we have
 $s(\phi)(v^*v)\neq 0$ and 
$$\phi(v)\,\, [v,s(\phi)]=[v,s(\phi)](v)\,\, \phi.$$  In particular, if
$\phi\in\Eigone(\C)$ then
$\lambda:=\frac{\phi(v)}{[v,s(\phi)](v)}\in\bbT$ and $\phi=\lambda[v,s(\phi)].$
\end{prop}

\begin{proof}
Easily, we have $[v,\rho](xd)=[v,\rho](x)\rho(d)$ and,
as $v^*v(v^*dx)=(v^*dv)v^*x$, we can apply $\rho$ to this equation
and divide by $\rho(v^*v)^{3/2}$ to obtain
\[ [v,\rho](dx)= \frac{\rho(v^*dv)}{\rho(v^*v)}
   [v,\rho](x)=\beta_v(d)[v,\rho](x),
\]
so $[v,\rho]$ is a $\D$-eigenfunctional with
range $\beta_v(\rho)$ and source $\rho$.
Letting $w=\rho(v^*v)^{-1/2} v$, $[v,\rho](w)=1$,
so $[v,\rho] \in \Eigone(\C)$.

For the converse, let $\phi \in \Eigone(\C)$ and set $\rho=s(\phi)$.
Since the span of $N_\D(\C)$ is norm dense in $\C$, there exists $v
\in N_\D(\C)$ with $\phi(v)\neq 0$, which by
Proposition~\ref{normalTOinter} may be taken to be an intertwiner if
desired.  Fix such a $v$. 

For any $n\in \bbN$,
$\phi(v(v^*v)^{1/n})=\phi(v)\rho(v^*v)^{1/n}$.  Since $v(v^*v)^{1/n}$
converges to $v$, $\rho(v^*v) > 0$.  Also, for any $d\in \D$ we have, 
\begin{align*}
r(\phi)(d)\phi(v)=\phi(dv)& =\frac{\phi(dv(v^*v))}{\rho(v^*v)} 
=\frac{\phi(v(v^*dv))}{\rho(v^*v)}\\
&\qquad =\frac{\phi(v)\rho(v^*dv)}{\rho(v^*v)}= \phi(v)\beta_v(\rho)(d),
\end{align*}
so $r(\phi)=\beta_v(\rho)$. 

For any unitary element $g\in\D$ we have,
\begin{align*}
\phi(vgv^*xg^{-1}) &= r(\phi)(vgv^*)\phi(x)\rho(g^{-1})\\
&=\frac{\rho(v^*vgv^*v)}{\rho(v^*v)}
\phi(x)\rho(g)^{-1} =\rho(v^*v)\phi(x).
\end{align*}
Theorem~\ref{exptHull} implies that $vE(v^*x)$ 
belongs to the closed convex hull of $\{vgv^*xg^{-1}:
g\in \U(\D)\}$, and thus $\phi(vE(v^*x))=\rho(v^*v)\phi(x)$. 
Hence for any $x\in \C$, 
\begin{align*}
\phi(v)\,\, [v,\rho](x) &= \frac{\phi(v) \rho(v^*x)}{\rho(v^*v)^{1/2}}
	= \frac{\phi(v E(v^*x))}{\rho(v^*v)^{1/2}}\\
& \qquad = \rho(v^*v)^{1/2}\phi(x) = [v,\rho](v)\, \, \phi(x).
\end{align*}
This equality, together with the
 fact that $\norm{[v,\rho]}=1,$ shows that $\lambda\in\bbT$ when $\phi$
 has unit norm.
\end{proof}

Proposition~\ref{EigTwist} leads to a description of $\Eig(\M)$ in
terms of $[v,\rho]$'s, for a norm-closed $\D$-bimodule $\M$
(Theorem~\ref{partA}).  Before giving this description, we use
Proposition~\ref{EigTwist} as a tool in the following result, which
gives the precise relationship between the concepts of regular
\cstar-inclusion and \cstardiag.

\begin{thm}\label{Dinclusion}
Let $(\C,\D)$ be a regular \cstar-inclusion, and let $\fN:=\{x\in \C:
E(x^*x)=0\}$ be the left kernel of $E$.
Then $\fN$ is a closed (two-sided) ideal of $\C$ and 
 $$\fN = \{x\in\C:\phi(x)= 0 \hbox{ for all } \phi\in\Eig(\C)\}.$$ 

Let $\pi: \C\rightarrow \C/\fN$ be the quotient map. Then
$(\pi(\C),\pi(\D))$ is a \cstardiag, $\pi|_{\D}$ is a $*$-isomorphism
of $\D$ onto $\pi(\D)$, and the restriction of the adjoint map
$\dual{\pi}|_{\Eig(\pi(C))}$ is an isometric isomorphism of
$\Eig(\pi(\C))$ onto $\Eig(\C)$.
\end{thm}

In particular, the first part of the theorem shows that, in
a \cstar-diagonal, $\Eigone(\C)$ separates points.

\begin{proof}
Since $\fN = \{x \in \C : E(x^*x)=0\}$ is clearly a closed left ideal,
to show that it is also a right ideal, it suffices to show that, for
$x \in \fN$ and a normalizer $v \in \C$, $xv \in \fN$.  To do this, we
shall prove that for every $\rho \in \hat{\D}$, $\rho(v^*x^*xv)=0$.
When $\rho(v^*v)=0$, this holds since $\rho(v^*x^*xv) 
\le \norm{x}^2\rho(v^*v)$.  When $\rho(v^*v)\neq 0$, let $\psi \in \hat{\D}$ be
given by
\[ \psi(z) = \frac{ \rho(v^*zv)}{\rho(v^*v)} \]
and observe that, as $x \in \fN$, $\psi(x^*x)=0$ and hence
$\rho(v^*x^*xv)=0$ in this case as well.  Thus $\fN$ is a closed
two-sided ideal.

We now show that $\fN= \{x\in\C:\phi(x)= 0 \hbox{ for all } \phi\in\Eig(\C)\}.$
Suppose that $x \in \fN$.
Given $\phi \in \Eig(\C)$, by Proposition~\ref{EigTwist}
we may assume $\phi=[v,\rho]$ where $\rho \in \hat{\D}$ and $v$ is a 
normalizer with $\rho(v^*v) > 0$.
By the Cauchy-Schwarz inequality
we have $|\rho(v^*x) |^2 \le \rho(v^*v) \rho(x^*x)$.
But $\rho(x^*x)=\rho(E(x^*x))=0$, so $[v,\rho](x)=0$.

Conversely, suppose $x\neq 0$ and $\phi(x)=0$ for every $\phi\in
\Eigone(\C)$.  We shall show that $\rho(x^*x)=0$ for every
$\rho\in\hat{\D}$.  So fix $\rho\in\hat{\D}$.  Notice that for every
normalizer $v\in\C$, we have $\rho(v^*x)=0$: this follows from the
Cauchy-Schwartz inequality when $\rho(v^*v)=0$ and from the hypothesis
and Proposition~\ref{EigTwist} when $\rho(v^*v)\neq 0$.  Let $\eps>0$.
Since the span of the normalizers is dense in $\C$, we may find
normalizers $v_1,\ldots,v_n$ so that $\norm{x - \sum_{i=1}^n v_i} <
\eps/\norm{x}$.  Thus,
\[ |\rho(x^*x)|=\left|\rho\left(x^*x-\sum_{j=1}^n v_j^*x  \right)\right|
	\leq \norm{x^*-\sum_{j=1}^n v_j^*}\,\norm{x}<\eps,
\] and we conclude that $\rho(x^*x)=0$.  Since this holds for every
$\rho\in\hat{D}$, we have $E(x^*x)=0$.

Clearly, $\fN \cap \D = 0$, so $\pi|_{\D}$
is an isomorphism of $\D$ onto $\pi(\D)$.
To see that $(\pi(\C),\pi(\D))$ is a \cstar-diagonal, observe first
that $E$ gives rise to a faithful conditional expectation on $\C/\fN$,
by the definition of $\fN$.    Given a pure state
$\rho$ of $\pi(\D)$, let $\tau_1$ and $\tau_2$ be pure states of
$\pi(\C)$ which extend $\rho$.  Now $\tau_i\circ \pi$ are extensions
of the pure state $\rho\circ\pi$ of $\D$ and hence coincide because
$\D$ has the extension property in $\C$.  Therefore, as $\pi$ is onto,
$\tau_1=\tau_2$,
and $\pi(\D)$ has the extension property in $\pi(\C)$.

To show that $\pi(\C)$ is regular (relative to $\pi(\D)$), let
$x+\fN\in\pi(\C)$, and $\eps >0$.  Since $\C$ is regular, we may
find normalizers $v_i\in\C$ such that $y:=\sum_{i=1}^{n} v_i$
satisfies $\norm{x-y}<\eps$.  Then $\pi(v_i)$ is a normalizer, and
since $\pi$ is contractive, $\norm{(x+\fN)-(y+\fN)}<\eps$. Hence
$\pi(\C)$ is regular.  It is also clear that
$\dual{\pi}|_{\Eig(\pi(\C))}$ is a homeomorphism of $\Eig(\pi(\C))$
onto $\Eig(\C)$ which preserves the partially defined product
structure.  Finally, since the adjoint of a quotient map is always
isometric, the proof is complete.
\end{proof}

We turn now to additional consequences of
Proposition~\ref{EigTwist}.
Recall that for any $f\in\dual{\C}$, $f^*$ is the bounded linear
functional given by $f^*(x)=\overline{f(x^*)}$.  It is easy to see
that if $\phi$ is an eigenfunctional on $\C$, then so is $\phi^*$ and
also that $s(\phi^*)=r(\phi)$ and $r(\phi^*)=s(\phi)$.  Thus, $f
\mapsto f^*$ provides an involution on $\Eigone(\C)$.  The inverse of
$[v,\rho]$ is $[v^*,\beta_v(\rho)]$.  Thus, we can summarize
Proposition~\ref{EigTwist} and the discussion preceding it.

\begin{thm}\label{Eigequal}
For a regular \cstar-inclusion $(\C,\D)$, $\Eigone(\C)=\Gamma$ and
the range and source maps, the involution, and topology are all the same.
\end{thm}

\begin{cor}\label{multiples}
If $\phi,\psi \in \Eig(\M)$ satisfy $r(\phi)=r(\psi)$ and
$s(\phi)=s(\psi)$,  then there exists $\lambda\in\bbC$ such that $ \lambda\neq
0$ and  $\phi=\lambda\psi$.
\end{cor}

\begin{proof}
Without loss of generality, assume $\norm{\phi}=\norm{\psi}=1$.  Let
$\rho=s(\phi)=s(\psi)$.  Theorem~\ref{extension} shows $\phi$ and
$\psi$ extend to norm-one eigenfunctionals on $\C$, which we denote by
the same symbols.  By Proposition~\ref{EigTwist}, we have
$\phi=[v,\rho]$ and $\psi=[w,\rho]$ for some $v,w \in \N(\C)$. It
follows from the hypothesis that $\beta_v(\rho)=\beta_w(\rho)$.  Thus
$\beta_{v^*w}(\rho)=\beta_{v*}(\beta_{w}(\rho)) =\rho$, so we can find
$\lambda \in \bbT$ with $\rho(\lambda v^*w) > 0$.  Hence by
Corollary~\ref{KumjianCor}, $[\overline{\lambda} v,\rho]=[w,\rho]$ and
therefore $\lambda \phi = \psi$.
\end{proof}

\begin{remark}{Remark}\label{rangesourceform}
Using the expression for $\phi=[v,s(\phi)]$ from Proposition~\ref{EigTwist}
and a short calculation with $r(\phi)(x) = \rho(v^*xv)/\rho(v^*v)$, 
where $\rho=s(\phi)$, we obtain for all $x \in \C$,
\[ \phi(x) = \frac{ s(\phi)(v^*x) }{[s(\phi)(v^*v)]^{1/2}}
           = \frac{ r(\phi)(xv^*) }{[r(\phi)(vv^*)]^{1/2}}.
\]
\end{remark}

Also, we can sharpen the inequality of Proposition~\ref{seminormprop}~(2)
to an equality.

\begin{cor}\label{eigennorm}
If $\phi \in \Eig(\M)$, then 
$|\phi(m)| = \|\phi\| B_{r(\phi),s(\phi)}(m)$ 
for all $m\in\M$. 
\end{cor}

\begin{proof}
By Proposition~\ref{seminormprop}, we have 
$|\phi(m)|\leq \norm{\phi}B_{r(\phi),s(\phi)}(m)$
for every $m\in\M$.

To obtain the reverse inequality, fix
$m\in\M$ such that $B_{r(\phi),s(\phi)}(m)\neq 0$. 
Given $\eps\in (0,1)$, we may find elements $a,b\in \D$
such that $\sigma(a)=\rho(b)=1$ and $B_{r(\phi),s(\phi)}(m) > (1-\eps)
\norm{amb}$.  Let $m_0:=amb$ and define a linear functional $f$ on
$\bbC m_0$ by 
$$f(tm_0)=t\norm{\phi}B_{r(\phi),s(\phi)}(m_0)=t\norm{\phi}B_{r(\phi),s(\phi)}
(m).$$  By the Hahn-Banach Theorem and Proposition~\ref{seminormprop},
 $f$ extends to an eigenfunctional $F$ on
$\M$ such that for every $x\in\M$, 
$$|F(x)|\leq \norm{\phi}B_{r(\phi),s(\phi)}(x).$$  Thus $\norm{F}\leq
\norm{\phi}$, $s(F)=s(\phi)$
and $r(F)=r(\phi)$. By Corollary~\ref{multiples}, there exists a
nonzero scalar $\lambda$ with $|\lambda|\leq 1$ and $F=\lambda\phi$.
We obtain,
\begin{align*}
|\phi(m)|&\geq |F(m)|=|F(m_0)|=|f(m_0)|=\norm{\phi} \frac{
B_{r(\phi),s(\phi)}(m)}{\norm{amb}}\norm{amb} \\
& > (1-\eps) \norm{\phi} \norm{amb}\\ 
& \geq (1-\eps)\norm{\phi}B_{r(\phi),s(\phi)}(m).
\end{align*}   
Letting $\eps\rightarrow 0$, we obtain the result. 
\end{proof}

We now extend Proposition~\ref{EigTwist} to bimodules and show
that the set of eigenfunctionals on a bimodule can be written in the
form $[v,\rho]$ with $v\in\M$.  

\begin{thm}\label{partA}
If $\phi \in \Eigone(\M)$, then there is a intertwiner $v\in\M$ 
and $\rho \in \hat{\D}$ so that $\phi = [v,\rho] |_\M$.

Conversely, if $v \in \N(\M)$, and $\rho \in \hat{\D}$ satisfies
$\rho(v^*v) > 0$.  then $[v,\rho] |_\M \in \Eigone(\M)$.
\end{thm}

\begin{proof}
We prove the second assertion first.  
Given $v\in\N(\M)$ and $\rho\in\hat{\D}$ with $\rho(v^*v)\neq0$,
Proposition~\ref{EigTwist} shows that $\phi:=[v,\rho]|_{\M}$ is an
eigenfunctional on $\M$.  Clearly, $\norm{\phi}\leq
\norm{[v,\rho]}=1$.  To show that $\norm{\phi}=1$, fix $m\in \M$
with $\phi(m)\neq 0$ and set $w=E_v(m)/\rho(v^*v)$.
Proposition~\ref{modmap} shows that $w\in \N(\M)$.  Calculation then
yields
\begin{equation}\label{rhoeqs}
\rho(w^*w)=|[v,\rho](m)|^2\neq 0, \qquad
[v,\rho](w)=[v,\rho](m).
\end{equation}
Let $\eps>0$. Pick a norm-one positive element $d\in \D$ with
$\rho(d)=1$ so that $\norm{d^2w^*w}\leq (1+\eps)\rho(w^*w)$.
If we let $s= wd/\rho(w^*w)^{1/2}$, then $\norm{s} \le (1+\eps)^{1/2}$
and the equations~\eqref{rhoeqs} give
$$\left|\phi\left(\frac{wd}{\rho(w^*w)^{1/2}}\right)\right|
=\left|[v,\rho]\left(\frac{wd}{\rho(w^*w)^{1/2}}\right)\right| =
\left| \frac{[v,\rho](w)}{\rho(w^*w)^{1/2}}\right|=1.$$

To prove the first assertion, let $\phi \in \Eigone(\M)$
and fix $m\in \M$ with $\phi(m)\neq 0.$
By Theorem~\ref{extension}, $\phi$
extends to a norm-one eigenfunctional, also called $\phi$, on $\C$.
By Proposition~\ref{EigTwist}, there is an intertwiner $v$ and
$\rho\in\hat{\D}$ with $\rho(v^*v) > 0$ so that $\phi=[v,\rho]$.
By Proposition~\ref{modmap},
$$w := \frac{vE(v^*m)}{\phi(m)\rho(v^*v)^{1/2}}\in \M.$$
and $w$ is an intertwiner. Thus,
$$\phi(w)=\phi(v)\frac{\rho(v^*m)}{\phi(m)\rho(v^*v)^{1/2}}
=\phi(v)\frac{[v,\rho](m)}{\phi(m)}=\phi(v),$$
so that $[w,\rho]=[v,\rho]=\phi$.
\end{proof}

\begin{remark}{Remark} \label{productstructure}
Let $\A\subseteq \C$ be a norm-closed algebra which is also a $\D$-bimodule.
We can use the multiplication defined for elements of $\Gamma$ 
to give a (also partially defined) multiplication on $\Eigone(\A)$ 
and hence on $\Eig(\A)$. 
Indeed, call $\phi_1,\phi_2 \in \Eigone(\A)$  composable if 
$s(\phi_1)=r(\phi_2)$.  
By Theorem~\ref{partA}, we can write $\phi_i=[v_i,\rho_i]$, where 
$v_1$ and $v_2$ are in $\N(\A)$.
Define $\phi_1\phi_2$ to be the product of the $[v_i,\rho_i]$
restricted to $\A$, namely $[v_1v_2,\rho_2] |_{\A}$.
\end{remark}

We can also improve several of the results of the previous section.
First, we immediately have a unique extension in Theorem~\ref{extension}.

\begin{cor}\label{extensionunique}
Suppose $\M_1, \M_2$ are norm-closed $\D$-bimodules with
$\M_1 \subseteq \M_2$.
There is a unique isometric map
$\iota:\Eig(\M_1)\rightarrow \Eig(\M_2)$ so that, for every
$\phi\in\Eig(\M_1)$, $ \iota(\phi) |_{\M_1} = \phi $.
The image $\iota(\Eig(\M_1))$ is
an open subset of $\Eig(\M_2)$.  

If in addition, $\M_2$ is regular,
then $\iota$ is $\sigma(\dual{\M_1},\M_1)$--$\sigma(\dual{\M_2},\M_2)$
continuous on bounded subsets of $\Eig(\M_1)$. 
\end{cor}

\begin{proof}
Fix $\phi\in \Eig(\M_1)$. Theorem~\ref{partA} shows $\phi$
extends uniquely to an eigenfunctional $\iota(\phi)$ on $\M_2$ of the
same norm.  
Writing $\phi=[v,\rho]$ for some normalizer $v\in\M_1$ and $\rho\in\hat{\D}$,  
then $\{\psi\in\Eig(\M_2): \psi(v)\neq 0\}$ is a $\sigma(\dual{\M_2},
\M_2)$-open set containing $\iota(\phi)$, so $\iota(\Eig(\M_1))$ is
an open set in $\Eig(\M_2)$.

It remains to prove that $\iota$ is continuous on bounded subsets of
$\Eig(\M_1)$ when $\M_2$ is regular.
Suppose $\phi_\lambda$ is a bounded net in $\Eig(\M_1)$ converging
$\sigma(\dual{\M_1},\M_1)$ to $\phi\in\Eig(\M_1)$.  Let $\rho=s(\phi)$
and $\rho_\lambda=s(\phi_\lambda)$. Then $\rho_\lambda$ converges
in the $\sigma(\dual{\C}, \C)$-topology to $\rho$.  

By Theorem~\ref{partA}, there exists a normalizer $v\in\M_1$ such that 
 $\phi=\norm{\phi}[v,\rho]$.  For large enough $\lambda$,
$\phi_\lambda(v)\neq 0$, so there exist scalars $t_\lambda\in\bbC$
with $|t_\lambda|=\norm{\phi_\lambda}$ and $\phi_\lambda = t_\lambda
[v,\rho_\lambda]$. Since $\phi_\lambda(v)\rightarrow \phi(v)$ we have
$t_\lambda\rightarrow \norm{\phi}$.  

For any normalizer $w\in\M_2$, we have
$$\iota(\phi_\lambda)(w)= t_\lambda [v,\rho_\lambda](w)=t_\lambda
\frac{\rho_\lambda(v^*w)}{\rho_\lambda(v^*v)^{1/2}}\rightarrow
\norm{\phi}[v,\rho](w) =\iota(\phi)(w).$$ As $\M_2$ is the span of the
normalizers it contains and $\phi_\lambda$ is a bounded net, we conclude
that for any $x\in\M_2$, $\iota(\phi_\lambda)(x)\rightarrow
\iota(\phi)(x)$.
\end{proof}

\begin{remark}{Remark}
Given a bimodule $\M\subseteq \C$, the set
$\Gamma_\M:=\{\phi\in\Gamma: \phi|_\M\neq 0\}$ is an open subset of
$\Gamma$ which plays a crucial role in the study of $\M$ (see, for
example, \cite{MR94i:46075,MR95a:46080}).
Theorem~\ref{partA}, together  with
Corollary~\ref{extensionunique}, shows that the restriction map
$[v,\rho]\in \Gamma_\M\mapsto [v,\rho]|_{\M}$ is a homeomorphism of 
 $\Gamma_\M$ onto $\Eigone(\M)$.  Thus,
$\Gamma_\M$ can be defined directly in terms of the bimodule structure
of $\M$, without explicit reference to $\C$.
\end{remark}

Since the norm is only lower semi-continuous for weak-$*$ convergence, it
is not possible to show that $\Eigone(\M)$ is locally compact for
general modules.  
However, for regular \cstar-inclusions, we can show this.

\begin{proposition}\label{loccmpt}
With the relative weak$^*$-topology,
$\Eigone(\M) \cup\{0\}$ is compact.
Thus, $\Eigone(\M)$ is a locally compact Hausdorff space.
\end{proposition}

\begin{proof}
Suppose that $\phi_\lambda$ is a net in $\Eigone(\M)\cup \{0\}$
which converges to $\phi\in (\dual{\M})_1.$
If $\phi=0$, there is nothing to do, so we assume that $\phi\neq 0$
and show that $\norm{\phi}=1$.

Fix a normalizer $v\in \M$ with $\phi(v) > 0$.
 From Theorem~\ref{partA}, $\phi=\norm{\phi} [v,\rho]$.
Choose a positive element $d\in\D$ so that $0\leq dv^*vd\leq I$ and
$\widehat{dv^*vd}=1$ in a neighborhood of $\rho$.
For large enough $\lambda$, $\phi_\lambda(vd)\neq 0$, so there exists
$t_\lambda\in\bbT$ such that $\phi_\lambda=t_\lambda [vd,s(\phi_\lambda)]$.
Thus, $1=|\phi_\lambda(vd)|$.
As $\phi_\lambda$ converges to $\phi$, we obtain $|\phi(vd)|=1=\norm{vd}$,
so $\norm{\phi}=1$.
\end{proof}

As usual, we may regard an element $m\in\M$ as a function on
$\Eigone(\M)$ via $\hat{m}(\phi)=\phi(m),$ and $\Eigone(\M)$ can be
regarded as a set of coordinates for $\M$.  Thus we make the following
definition.   

\begin{defn}\label{coordsys}
For a norm-closed $\D$-bimodule $\M$, we call the set $\Eigone_\D(\M)$,
equipped with the relative weak-$*$ topology, the $\bbT$-action, and
the range and source mappings, a \textit{coordinate system} for $\M$.  

When $\A$ is both a norm-closed algebra and a $\D$-bimodule, the
coordinate system $\Eigone(\A)$ 
also has the additional structure of a continuous
partially defined product as described in
Remark~\ref{productstructure}.  In this case we will sometimes refer
to the coordinate system as a \textit{semitwist}.
\end{defn}

\begin{defn}\label{defnsim}
If $\M \subseteq \C$ is a $\D$-bimodule, let $R(\M):=\{|\phi|:
\phi\in\Eigone(\M)\}$.  Then $R(\M)$ may be identified with the
quotient $\Eigone(\M)\backslash\bbT$ of $\Eigone(\M)$ by the natural
action of $\bbT$. Obviously, $\phi\mapsto |\phi|$ is the quotient map,
and the topology on $R(\M)$ is the quotient topology.
Corollary~\ref{extensionunique} shows that we may regard $R(\M)$ as a
subset of $R(\C)$, and, as $\C$ is regular,
Corollary~\ref{extensionunique} also implies that if $v\in\M$ is an
intertwiner, then $G_v:=\{\phi\in R(\M): |\phi(v)|>0\}$ is an open set
and $\{G_v: v\in\M \text{ is an intertwiner}\}$ is a base for the
topology of $R(\M)$.  Thus $R(\M)$ is a locally compact Hausdorff
space.

We shall sometimes find it useful to view $R(\M)$ as a topological
relation on $\hat{\D}$.  The map $|\phi|\mapsto (r(\phi), s(\phi))$ is
a bijection between $R(\M)$ and
$\{(r(\phi),s(\phi))\in\hat{\D}\times\hat{\D}: \phi\in\Eigone(\M)\}$,
and we will sometimes identify these two sets under this bijection.
With this identification, $(\sigma,\rho)\in R(\M)$ if and only if
there is an intertwiner $v\in\M$ with $\rho(v^*v)\neq 0$ and $\sigma=
\beta_v(\rho)$; moreover, the set $G_v$ is the graph of $\beta_v$ and
the collection of such sets  gives a base for the topology.
We call $R(\M)$ the \textit{spectral relation} of
$\M$.
Also, $R(\M)$ is reflexive if $\D \subseteq \M$, is
symmetric if $\M=\M^*$ and is transitive if $\M$ is a subalgebra.

A \textit{topological equivalence relation} is a principal topological
 groupoid.  If $v, w$ normalize $\D$, then $G_{vw}=\{ |\phi \psi|:
 s(\phi)=r(\phi), \phi\in G_v, \psi\in G_w\}$ and $G_{v^*} = \{
 |\phi^*| : \phi \in G_v\}$.  It follows that the topology on $R(\C)$
 is compatible with the groupoid operations, so $R(\C)$ is a
 topological equivalence relation.  We will sometimes write
 $\sigma \sim_{\C} \rho$, or simply $\sigma\sim \rho$, when
 $(\sigma,\rho)\in R(\C)$.
\end{defn}

We now show that the regularity of $(\C,\D)$ and the faithfulness of $E$ 
imply that the span of eigenfunctionals is weak-$*$ dense in $\dual{\M}$, 
for any norm closed $\D$-bimodule $\M$.   

\begin{proposition}\label{Mdensity}  
Suppose $(\C,\D)$ is a \cstardiag,  and let
  $\M\subseteq \C$ be a norm-closed $\D$-bimodule.  Then 
$\spn\Eigone(\M)$ is $\sigma(\dual{\M},\M)$-dense in $\dual{\M}$.
\end{proposition}
\begin{proof}
Let $W$ be the $\sigma(\dual{\M},\M)$-closure of $\spn\Eigone(\M).$ If
  $W\neq \dual{\M}$, then there exists a nonzero
  $\sigma(\dual{\M},\M)$-continuous linear functional $\psi$ on
  $\dual{\M}$ which annihilates $W$.  Since $\psi$ is
  $\sigma(\dual{\M},\M)$-continuous, there exists $m\in\M$ such that
  $\psi(f)=f(m)$ for all $f\in \dual{\M}$.  But then for every
  $\phi\in\Eigone(\M)$, we have $\psi(\phi)=\phi(m)=0$.
  Since $\Eigone(\C)$ separates points (Theorem~\ref{Dinclusion}),
  there exists an eigenfunctional
  $\overline{\phi}\in \Eigone(\C)$ so that $\overline{\phi}(m)\neq 0$.
  The restriction $\phi:=\overline{\phi}|_\M$ is a 
  eigenfunctional on $\M$, so
  $0=\psi(\phi)=\phi(m)=\overline{\phi}(m)\neq 0$, a
  contradiction. Therefore, $W=\dual{\M}$.
\end{proof}

We conclude this section with two applications of the results in this
section.  For our first application,  we give a description of the
\cstar-envelope of an algebra satisfying $\D\subseteq \A\subseteq \C$.

\begin{theorem}
\label{envelope}
Let $(\C,\D)$ be a \cstardiag\ and suppose $\A$ is a norm closed
algebra satisfying $\D\subseteq\A \subseteq \C$. 
If $\B$ is the \cstar-subalgebra of $\C$ generated by $\A$, then
$\B$ is the \cstar-envelope of $\A$.

If in addition, $\B=\C$, 
 then $R(\C)$ is the topological equivalence 
relation  generated by $R(\A)$.
\end{theorem}

\begin{proof}
Let $\B_e$ be the \cstar-envelope of $\A$ and let $j:\A \rightarrow
\B_e$ be the canonical embedding (see, for example,
\cite[Section~4.3]{BlecherLeMerdyOpAlThMo}).  Then there exists a
unique $*$-epimorphism $\pi:\B\rightarrow \B_e$ such that
$\pi(a)=j(a)$ for every $a\in\A$.  In particular, $\pi$ is faithful on
$\D$.  Assume, to get a contradiction, that $\pi$ is not injective.
Then $\ker\pi$ is a $\D$-bimodule and let $x$ be a non-zero element of
$\ker\pi$.  By
Theorem~\ref{Dinclusion}, eigenfunctionals separate points, so there
exists an element $\phi\in\Eigone(\C)$ with $\phi(x)\ne 0$.  Writing
$\phi=[v,\rho]$, Proposition~\ref{modmap} shows that $u:=E_v(x)$ is a
nonzero normalizer belonging to $\ker\pi$. 
  But then $u^*u$ is a nonzero element
of $\D\cap \ker\pi$, a contradiction.  Thus, $\pi$ is faithful on
$\B$, and hence $\B$ is the \cstar-envelope of $\A$.

Suppose now that $\B=\C$.  By Corollary~\ref{extensionunique}, there
 is an inclusion $\Eigone(\A)\subseteq \Eigone(\C)$.  So
 $R(\A)\subseteq R(\C)$, and hence $R(\C)$ contains the equivalence
 relation generated by $R(\A)$.
 
For the other direction, assume that $(\sigma,\rho)\in R(\C)$.  Then
there is $\phi \in \Eigone(\C)$ with source $\rho$ and range $\sigma$.
Let $\W$ be the set of all finite products of intertwiners belonging
to $\A$ or to $\A^*$.
Then $\W\subseteq \N_\D(\C)$ and the set of finite sums from $\W$
is a $*$-algebra which, because $\A$ generates $\C$, is dense in $\C$.
Hence there is some $w \in \W$ so that $\phi(w) \ne 0$.
By Proposition~\ref{EigTwist}, $\phi=[w,\rho]$.

Suppose that $w$ factors as $v_{2n}^*v_{2n-1}\ldots v_2^* v_1$,
where each $v_i$ is an intertwiner in $\A$.
Let $\rho_1=\rho$ and for $i=2,\ldots,2n$, let $\rho_i$ be image of
$\rho$ under conjugation by the rightmost $i-1$ factors in the factorization
of $w$.
It follows that $\phi$ is the product
\[ [v_{2n},\rho_{2n}]^* [v_{2n-1},\rho_{2n-1}] \cdots
   [v_2,\rho_2]^* [v_1,\rho_1]
\]
and each $[v_i,\rho_i]$ is in $\Eigone(\A)$.
Thus, the equivalence relation generated by $R(\A)$ contains
$(\sigma,\rho)$.
Similar arguments apply for the other possible factorizations of $w$,
so the equivalence relation generated by $R(\A)$ contains $R(\C)$.

It remains to show that the usual topology on $R(\C)$ equals that generated by
$R(\A)$, i.e., the smallest topology containing the topology of $R(\A)$ which 
makes $R(\C)$ into a topological equivalence relation.  

As we noted in Definition~\ref{defnsim}, $R(\C)$ is already a
topological equivalence relation, and so its topology contains the
topology generated by $R(\A)$.  Since the norm-closed span of $\W$ is
$\C$, it follows that $\{G_w : w\in\W\}$ is a base for the topology of
$R(\C)$, where, as before, $G_w = \{ \phi \in R(\C) : |\phi(w)| > 0
\}$.  For a topological equivalence relation, the inverse map is a
homeomorphism.  Further, given two precompact open G-sets (i.e., a
subset of $R(\C)$ on which the two natural projection maps into
$\hat{\D}$ are injective), one can show that their product is again a
precompact open G-set (for example, adapt the proof of Proposition
I.2.8 in~\cite{MR82h:46075}).  Since, for $v$ an intertwiner in $\A$,
each $G_v$ is a precompact open G-set in $R(\A)$, and each $w \in W$
is a finite product of such $v$'s and their inverses, it follows that
each $G_w$, $w \in W$, is open in the topology generated by $R(\A)$.
Thus, the topology generated by $R(\A)$ contains $R(\C)$.
\end{proof}

Our second application is an application of
Theorem~\ref{Dinclusion}. We show
show that inductive limits of \cstardiag s are again \cstardiag s,
when the connecting maps satisfy a certain condition, which we now
define.  The difficulty in showing that these inductive limits are
again \cstardiag s is in showing that the expectation is faithful, and
this is where Theorem~\ref{Dinclusion} provides a key tool.

\begin{defn}
Given regular \cstar-inclusions $(\C_i,\D_i)$, $i=1,2$, and a $*$-homomorphism
$\pi : \C_1 \to \C_2$, we say $\pi$ is \textit{regular} if 
$\pi(\N(\C_1)) \subseteq \N(\C_2)$.
\end{defn}

Of course, if $\pi$ is regular, then $\pi(\D_1) \subseteq \D_2$.
Indeed, for $D \in \D_1$ with $D \ge 0$, $D^{1/2} \in \N(\C_2)$
and so $\pi(D)=\pi(D^{1/2})1\pi(D^{1/2}) \in \D_2$.

\begin{thm}\label{indlimit}
Let $(\C_\lambda,\D_\lambda)$, $\lambda \in \Lambda$, be a directed net
of regular \cstar-inclusions with regular $*$-monomorphisms
$\pi_{\lambda,\mu}: \C_\mu \to C_\lambda$.  Then
$\bigl(\indlimit(\C_\lambda,\pi_{\lambda,\mu}),
\indlimit(\D_\lambda,\pi_{\lambda,\mu})\bigr)$ is a regular
\cstar-inclusion.  Moreover, if each $(\C_\lambda,\D_\lambda)$ is a
\cstar-diagonal, then so is 
$(\indlimit(\C_\lambda,\pi_{\lambda,\mu}),
\indlimit(\D_\lambda,\pi_{\lambda,\mu})\bigr)$.
\end{thm}

\begin{proof}
The first part of the proof is routine.  We regard the $\C_\lambda$ as
a $*$-subalgebras of $\C := \indlimit(\C_\lambda,\pi_{\lambda,\mu})$
and identify $\pi_{\lambda,\mu}$ with the inclusion map from
$\C_\lambda$ to $\C_\mu$.  Then $\D =
\indlimit(\D_\lambda,\pi_{\lambda,\mu})$ is a subalgebra of $\C$.

Given a normalizer $v \in \C_\lambda$, by the regularity of the
inclusion maps, $v$ normalizes $\D_\mu$ for all $\mu \ge \lambda$.
Thus, $v$ normalizes $\D$ and so $\N_{\D_\lambda}(\C_\lambda)
\subseteq \N_{\D}(\C)$.  Since $\C$ is the closed union of the
$\C_\lambda$, and each $\C_\lambda$ is the span of
$\N_{\D_\lambda}(\C_\lambda)$, $\C$ is regular in $\D$.

Given $\rho\in \hat{\D}$, suppose $\phi$ and $\psi$ are 
extensions of $\rho$ to states of $\C$.  Then, for each $\lambda \in
\Lambda$, $D_\lambda \subseteq \D$ and so $\phi|_{\C_\lambda}$ and
$\psi|_{\C_\lambda}$ are  extensions of the pure state
$\rho|_{\D_\lambda}\in\hat{\D}_\lambda$ and so agree on $\C_\lambda$.
Since $\C$ is the closed union of the $\C_\lambda$, $\phi=\psi$.
Thus, $(\C,\D)$ is a regular \cstar-inclusion.

Let $E : \C \to \D$ be the expectation.  By Theorem~\ref{Dinclusion},
$\fN:=\{ x \in \C : E(x^*x)=0\}$ is an ideal of $\C$, and, if $q :
\C \to \C/\fN$ is the quotient map, then $(q(\C),q(\D))$ is a
\cstar-diagonal.  If $x \in \C_\lambda$ and $E(x^*x)=0$, then
$\rho(x^*x)=0$ for all $\rho \in \hat{\D}$.  Since every $\sigma \in
\hat{\D}_\lambda$ has at least one extension to an element of
$\hat{D}$, we have $E_\lambda(x^*x)=0$, where $E_\lambda$ is the
expectation for $(\C_\lambda,\D_\lambda)$.
Thus, 
if $(\C_\lambda,\D_\lambda)$ is a \cstar-diagonal, then we have
$x=0$, that is, $\fN\cap \C_\lambda=(0)$, 
so $q$ is faithful on $\C_\lambda$.  

Therefore, when each $(\C_\lambda,\D_\lambda)$ is a \cstardiag,
$q(\C)$ contains isomorphic copies of each $\C_\lambda$, and when
$\lambda\leq \mu$, $q(\C_\lambda)\subseteq q(\C_\mu)$.   By the
minimality of the inductive limit, $q$ is an isomorphism of $\C$
onto $q(\C)$, i.e.\ $\fN=0$.    Thus, $(\C,\D)$ is a \cstar-diagonal.
\end{proof}

\section{Compatible Representations of \cstardiag s}\label{repns}

Our goal in this section is to produce a faithful representation $\pi$
of a  \cstardiag\ $(\C,\D)$.  
Because we require the faithfulness of the expectation, we work with 
\cstardiag s instead of regular \cstar -inclusions.  

\begin{remark*}{Standing Assumptions for Section~\ref{repns}}
We assume $(\C,\D)$ is a \cstar-diagonal.
For $(\C,\D)$ a \cstardiag, we write 
$\A\subseteq (\C,\D)$ if $\A\subseteq \C$ is a norm-closed subalgebra 
with $\D\subseteq \A$.

For $\rho \in \hat{\D}$, we use $(\H_\rho,\pi_\rho)$ for the GNS
representation of $\C$ associated to the (unique) extension of $\rho$.
\end{remark*}

Eigenfunctionals can be viewed as normal linear functionals on $\ddual{\C}$
and we start by using the polar decomposition for such functionals
to obtain a `minimal' partial isometry for each eigenfunctional.
Although these results are implicit in the development of dual groupoids
(see~\cite[p.~435]{MR87a:46111}), we give a (mostly) self-contained
treatment.

Fix a norm-one eigenfunctional $\phi$ on $\C$.
By the polar decomposition for linear functionals
(see~\cite[Theorem~III.4.2, Definition~III.4.3]{TakesakiThOpAlI}),
there is a partial isometry $u^*\in\ddual\C$ and a positive linear functional
$|\phi| \in\dual{\C}$ so that $\phi=u^*\cdot |\phi|=|\phi^*|\cdot u^*.$
Applying the characterization given in
\cite[Proposition~III.4.6]{TakesakiThOpAlI}, we find that
\begin{equation}\label{absval}
r(\phi)=|\phi| \quad \text{and}\quad s(\phi)=|\phi^*|.
\end{equation}
Moreover, $uu^*$ and $u^*u$ are the smallest projections in
$\ddual{\C}$ which satisfy,
\begin{align*}
u^*u\cdot s(\phi)&=s(\phi)\cdot u^*u=s(\phi)\\
\intertext{and}
 uu^*\cdot r(\phi)&= r(\phi)\cdot uu^*=r(\phi).
\end{align*}

\begin{defn}\label{asparisom}
For $\phi \in \Eigone(\C)$, we call the partial isometry $u$ above 
the \textit{partial isometry associated to $\phi$} and denote it by $v_\phi$.
If $\phi \in \hat{\D}$, then $u$ is a projection and we denote it by
$p_\phi$.
\end{defn}

\begin{remark}{Remark}\label{minpifacts}
The above equations show that $v_\phi^* v_\phi = p_{s(\phi)}$
and $v_\phi v_\phi^* = p_{r(\phi)}$.  Moreover, given
$\phi\in\Eigone(\C),$ $v_\phi$, Proposition~\ref{combineeigen} 
below implies that may be characterized as the unique
minimal partial isometry $w\in\ddual{\C}$ such that $\phi(w)>0$.
\end{remark}

Our first goal is to show that the initial and final projections of this
partial isometry are minimal projections in $\ddual{\C}$ and
compressing by them gives $\phi$, in the following sense.

\begin{proposition}\label{combineeigen}
For $\rho\in\hat{\D}$, $p_\rho=p_{\rho \circ E}$ 
is a minimal projection in $\ddual{\C}$.
For all $\phi\in\Eigone(\C)$ and $x\in \ddual{\C}$,
$$ p_{r(\phi)} x p_{s(\phi)} =\phi(x)v_\phi.$$
\end{proposition}

\begin{proof}
First, we show that $p_\rho$ is a minimal projection in $\ddual{\D}$.
We know that $p_\rho$ is the smallest projection in $\ddual{\D}$ such that
$p_\rho\cdot \rho=\rho\cdot p_\rho =\rho$.
Suppose, to get a contradiction, that
$p_1, p_2$ are nonzero projections in $\ddual{\D}$ with
$0 \le p_1,p_2 \le p$ and $p_\rho=p_1+p_2$.

If $\rho(p_1)=0$, then for $d\geq 0$,
$\rho(p_1d)=\rho(p_1dp_1)\leq \rho(p_1)\norm{d}=0$
and so, for all $d\in \D$,  $\rho(p_1d)=0$.
But then $p_2\cdot\rho=p\cdot\rho=\rho=\rho\cdot p_2$,
which yields $p=p_2$, contrary to hypothesis.
Hence $\rho(p_1)\neq 0$ and, similarly, $\rho(p_2)\neq 0$.

This implies that $\rho$ can be written as a nontrivial convex
combination of states on $\D$, for
$$\rho= \rho(p_1) \frac{p_1\cdot \rho}{\rho(p_1)} + \rho(p_2)
\frac{p_2\cdot \rho}{\rho(p_2)}.$$
But this is a contradiction, since elements of of $\hat{\D}$ are 
pure states.

Since $p_\rho$ is minimal in $\ddual{\D}$, for $d\in \D$, 
$p_\rho dp_\rho =\rho(d) p_\rho$.

Now suppose that $q\in\ddual{\C}$ is a projection with $0< q\leq p_\rho$.
Since $q\neq 0$, there exists a state $g\in\dual{\C}$
such that $g(q)>0$.  Define
$$f:= \frac{q\cdot g \cdot q}{g(q)}.$$
Then $f$ is a state on $\C$ with $f(q)=1$.
As $p_\rho q=q$, we have for $d\in \D$,
$$f(d)=\frac{g(q(p_\rho dp_\rho )q)}{g(q)}=\rho(d).$$
Since pure states on $\D$ extend uniquely to pure states on $\C$,
we conclude that $f=\rho\circ E$.

If $p_\rho$ is not minimal, write $p_\rho=q_1+q_2$ where $q_i\in\ddual{\C}$ are
projections with $0 < q_1,q_2 \le p_\rho$.
Apply the argument of the previous paragraph to find states $h_1$ and
$h_2$ on $\C$ such that $h_i(q_i)=1$ and $q_i\cdot h_i\cdot q_i=h_i$.
Since $q_1q_2=0$, $h_1(q_2)=h_2(q_1)=0$.
But the previous paragraph shows that $h_1=\rho\circ E=h_2$,
contradicting the extension property.
So $p_\rho$ is minimal in $\ddual{\C}$.

The uniqueness of polar decompositions implies that $p_{\rho}=p_{\rho
\circ E}$.

To prove the second statement, first note that
$p_{r(\phi)} \ddual{\C} p_{s(\phi)}$ has dimension one, since
$p_{r(\phi)},p_{s(\phi)}$ are minimal projections in $\ddual{\C}$ and
$p_{r(\phi)} v_\phi p_{s(\phi)} = v_\phi \ne 0$. 
Hence there is a linear functional $g$ on $\C$ such that
for every $x\in\ddual{\C}$, $p_{r(\phi)}xp_{s(\phi)}=g(x) v_\phi$.
Then $g |_{\C}$ is an eigenfunctional with the same source
and range as $\phi$.
Since $g(v_\phi)=\phi(v_\phi)$, $g=\phi$.
\end{proof}

Recall (\cite[Lemma~III.2.2]{TakesakiThOpAlI}) that any
$*$-representation $\pi$ of a \cstaralg\ $\C$ has a unique extension
to a $*$-representation $\tilde{\pi}: \ddual{\C} \to \pi(\C)''$,
continuous from the $\sigma(\ddual{\C},\dual{\C})$-topology
to the $\sigma$-weak topology on $\pi(\C)''$, i.e.,
the $\sigma(\B(\H),\B(\H)_*)$-topology.
Let $\I=\ker\tilde{\pi}\subseteq \ddual{\C}$.
By the continuity of $\tilde{\pi}$, $\I$ is $\sigma(\ddual{\C},\dual{\C})$
closed, and hence (see~\cite[Proposition~1.10.5]{MR0442701})
there exists a unique central projection $P\in\ddual{\C}$ such that $\I=\ddual{\C}(I-P).$
Further, $\tilde{\pi}|_{\ddual{\C}P}$ is one-to-one (see
\cite[Definition~1.21.14]{MR0442701}) and is onto $\B(\H)$.
The projection $P$ is called the \textit{support projection} for $\pi$.

Recall that $\sigma,\rho \in \hat{\D}$ have $(\rho,\sigma) \in R(\C)$
if and only if there is $\phi \in \Eigone(\C)$ with $r(\phi)=\rho$
and $s(\phi)=\sigma$.
For brevity, we write $\sigma \sim \rho$ in this case.

\begin{proposition} \label{specialrep}
If $\rho\in\hat{\D}$, then $\pi_\rho(\D)''$ is an atomic MASA in
$\B(\H_\rho)$ and the support projection for $\pi_\rho$ is
$\sum_{\sigma \sim \rho} p_\sigma$.  Moreover, the map from
$\pi_\rho(\C)$ to $\pi_\rho(\D)$ given by $\pi_\rho(c)\mapsto
\pi_\rho(E(c))$ is well-defined and extends to a faithful normal
expectation $\tilde{E}:\pi_\rho(\C)'' \to \pi_\rho(\D)''$.
\end{proposition}

\begin{proof}
Let $\pi=\pi_\rho$ and $\H=\H_\rho$.
Since the extension of $\rho$ to $\C$ is pure, the representation
$\pi$ is irreducible, so $\pi(\C)''=\bh$.

Recall from~\cite[Corollary~1]{KadisonIrOpAl} that if $\M:=\{x\in \C:
\rho(x^*x)=0\}$, then $\C/\M$ is complete relative to the norm induced
by the inner product $\innerprod{x+\M,y+\M}=\rho(y^*x)$, and thus
$\H=\C/\M$.  Our first task is to obtain a convenient orthonormal
basis for $\H$ and, towards this end, we require the following observation.

If $v, w \in \N(\C)$ with $\rho(v^*v)=1$ and $\rho(v^*w)\neq 0$, then
\begin{equation}\label{Mpro}
w-[v,\rho](w)v\in\M.
\end{equation}
To see this,  let $w_1=
\rho(w^*v) w$.  Then $\rho(v^*w_1)>0$, so that by
Corollary~\ref{KumjianCor}, $[v,\rho](w_1)=[w_1,\rho](w_1)$, and thus
$\rho(v^*w_1)= \rho(w_1^*w_1)^{1/2}$.  It follows that $|\rho(v^*w)|=
\rho(w^*w)^{1/2}.$ A calculation then shows that $w-[v,\rho](w)v=
w-\rho(v^*w)v\in \M$, as required. 

Choose a set $\Z\subseteq \C$ of normalizers such that for each $z\in
\Z$, $\rho(z^*z)=1$ and the map $z\mapsto r([z,\rho])$ is a bijection
of $\Z$ onto $\O := \{ \sigma \in \hat{\D}: \sigma \sim \rho\}$.  

We claim that if $X=\sum_{j=1}^n w_j\in \C$ with each $w_j\in\N_\D(\C)$, then
\begin{equation}\label{MproB} 
X+\M = \sum_{z\in\Z} [z,\rho](X)(z+\M).
\end{equation}
To see that the summation is well-defined, first observe that if, for some 
$z \in \Z$, $\rho(z^*w_j)\neq 0$, then $r([w_j,\rho])=r([z,\rho]) \in \O$.
Thus, for each $j$, $\{z\in\Z: \rho(z^*w_j)\neq 0\}$ is a singleton,
and so $\{z\in\Z: \rho(z^*X)\neq 0\}$ is finite.
To prove the equality, interchange the order of summation and use~\eqref{Mpro}
as follows:
\begin{align*}
\sum_{z\in\Z} [z,\rho](X)(z+\M)
     &= \sum_{z\in\Z}\sum_{j=1}^n [z,\rho](w_j)(z+\M)\\
     &= \sum_{j=1}^n\sum_{z\in\Z} [z,\rho](w_j)(z+\M)\\
     &= \sum_{j=1}^n\sum_{\substack{z\in\Z\\ \rho(z^*w_j)\neq 0}} [z,\rho](w_j)(z+\M)\\
     &= \sum_{j=1}^n w_j + \M & \text{by \eqref{Mpro}} \\
     &= X+\M.
\end{align*} 

Next, we show that $\{z+\M: z\in\Z\}$ is an orthonormal basis for $\H$.
By Corollary~\ref{KumjianCor}, it is an orthonormal set.
Given any $Y\in\C$ and $\eps>0$, we may find a finite sum of
normalizers $X$ so that $\norm{X-Y}_\C <\eps$.
Then $\norm{X-Y+\M}_{\H_\rho}^2=\rho((X-Y)^*(X-Y))\leq \norm{X-Y}_\C^2 <\eps^2$.
As $X+\M$ is in $\spn\Z$ and $\eps$ is arbitrary, 
$Y+\M\in \overline{\spn\Z}$, so that $\{z+\M: z\in\Z\}$ is an orthonormal basis.

For $d\in\D$ and $z\in \Z$, using~\eqref{Mpro}
\begin{equation}\label{aMASA}
\pi(d)(z+\M)=dz+\M=[z,\rho](dz) z+\M = \rho(z^*dz) (z+\M).
\end{equation}
Thus, $\pi(d)$ is diagonal with respect to the basis $\{z+\M\}_{z\in\Z}$.
Fixing $z\in\Z$, let 
$\Lambda=\{d\in\D: d\geq 0, \text{ and } r([z,\rho])(d)=1\}$.  
Then $\Lambda$ is a directed set under the ordering $d_1\preceq d_2$
if and only if $d_1\geq d_2$. 
It is easy to see that the net $\{\pi(d)\}_{d\in\Lambda}$ decreases 
to the rank-one projection onto $z+\M$.  
Thus $\pi(\D)''$ is an atomic MASA.

Fix $z\in\Z$ and let $P_z$ be the orthogonal projection of $\H$ onto $z+\M$.  
Then a calculation shows that for $x\in\C$,
\begin{equation}\label{compressbypv}
P_z\pi(x)P_z= \rho(z^*xz) P_z.
\end{equation}
Since $r([z,\rho])$ is the vector state corresponding to $z+\M$ and 
$\tilde{\pi}$ is normal, \eqref{compressbypv} holds when 
$\pi$ is replaced by $\tilde{\pi}$ and $x\in\ddual{\C}$.  
As $p_{r([z,\rho])}$ is a minimal projection in $\ddual{\C}$, 
$\tilde{\pi}(p_{r([z,\rho])}) =P_z$.
If $Q=\sum_{\sigma\sim\rho} p_\sigma$, then $\tilde{\pi}(Q)=I$.
Letting $P$ be the support projection of $\pi$, this implies that $P \le Q$.
If $\sigma \sim \rho$, then $p_\sigma$ is a a minimal projection satisfying
$\tilde{\pi}(p_\sigma)\neq 0$.
Thus, $p_\sigma\leq P$ and so $Q\leq P$.
This shows that $Q$ is the support projection for $\pi$ and that
$Q\in\ddual{\D}$.

Finally, \eqref{aMASA} implies that for $x\in\C$, 
\begin{align*}
\pi(E(x))(z+\M)&=
\rho(z^*E(x)z)(z+M)=r([z,\rho])(E(x))(z+\M)\\
&=r([z,\rho])(x)(z+\M)=\rho(z^*xv)(z+\M).
\end{align*}  Hence, for $x\in\C$, we obtain
$$\pi(E(x))=  \sum_{z\in\Z} P_z \pi(x) P_z.$$
If $\tilde{E}:\bh\rightarrow \bh$ is defined by
$\tilde{E}(T)=\sum_{z\in\Z}P_z T P_z$, this shows that $\tilde{E}$ is a
faithful normal conditional expectation of $\bh=\pi(\C)''$ onto
$\pi(\D)''$ satisfying $\tilde{E}(\pi(x))=\pi(E(x))$ for $x\in\C$.
\end{proof}

We next record a simple consequence of the proof of
Proposition~\ref{specialrep}.   
\begin{corollary}\label{eigvecftn}
Suppose $\rho\in\hat{\D}$ and $\phi\in \Eigone(\C)$ satisfies
$\rho\sim s(\phi)$.  Then there exist unit orthogonal unit
vectors $\xi,\eta\in\H_\rho$ such that for every $x\in\C$,
$\phi(x)=\innerprod{\pi_\rho(x)\xi,\eta}$. 
\end{corollary}

\begin{proof}
With the same notation as in the proof of
Proposition~\ref{specialrep}, there exists $z,w\in \Z$ such that
$r([z,\rho])=s(\phi)$ and $r([w,\rho])=r(\phi)$.  Let $\xi_0=z+\M$ and
$\eta=w+\M$.  For $x\in\C$, let
$\psi(x)=\innerprod{\pi(x)\xi_0,\eta}$.  Observe that
$\psi(wz^*)=\rho(w^*wz^*z)=1$, so $\psi$ is nonzero.  Also, for
$x\in\C$ and $d\in\D$, we have
$$\psi(xd)=\rho(w^*xdz)=\frac{\rho(w^*xdzz^*z)}{\rho(z^*z)}=\rho(w^*xz)
\frac{\rho(z^*dz)}{\rho(z^*z)}=\psi(x)s(\phi)(d).$$
Similarly, $\psi(dx)=r(\phi)(d)\psi(x)$.  Thus $\psi$ is an
eigenfunctional with the same range and source as $\phi$.
Hence, there exists $t\in\bbT$ so that $\phi=t\psi$.  
Take $\xi=t\xi_0$.
\end{proof}

\begin{defn}
Given a \cstar-diagonal $(\C,\D)$, a
representation $\pi$ of $\C$ is $\D$-\textit{compatible}, (or simply
\textit{compatible}) if $\pi(\D)''$ is a
MASA in $\pi(\C)''$ and there exists a faithful conditional
expectation $\tilde{E}:\pi(\C)''\rightarrow \pi(\D)''$ such that for
every $x\in \C$, $\tilde{E}(\pi(x))= \pi(E(x))$.
\end{defn}

Proposition~\ref{specialrep} shows that the GNS representation of $\C$
associated to an element of $\hat{\D}$ is compatible.

\begin{remark}{Remark}
Suppose that $X$ is an r-discrete locally compact principal groupoid with unit
space $X^0$ and let $\sigma$ be a 2-cocycle.  Then
Drinen~\cite{DrinenViAFAlGrAl} 
shows that $(C_r(X,\sigma), C_0(X^0))$
is a \cstardiag.  We expect that when $\lambda^u$ is a Haar system on
$X$, and $\mu$ is a measure on $X^0$, the induced representation
$\rm{Ind}(\mu)$ of $C_r(X,\sigma)$ (see for example,
\cite[page~44]{MR90m:46098}) is a compatible representation of
$(C_r(X,\sigma), C_0(X^0))$, and it would not be surprising if every
compatible representation for $(C_r(X,\sigma), C_0(X^0))$ arises in
this way.  We do not pursue this issue here, however.
\end{remark}

\begin{lemma}\label{unitequiv}
Let $\rho,\sigma\in\hat{\D}$.  Then 
$\rho\sim \sigma$ 
if and only
if the  GNS representations $\pi_\sigma$ and $\pi_\rho$
are unitarily equivalent.
\end{lemma}
\begin{proof}
If $U\in\B(\H_\rho,\H_\sigma)$ is a unitary such that
$U\pi_\rho U^*=\pi_\sigma$, the fact that $\C/\M_\rho=\H_\rho$ allows
us to find $X\in\C$ such that $U^*(I+\H_\sigma)= X+\M_\rho$.  Then for
all $x\in\C$, we have $\sigma(x)=
\innerprod{\pi_\sigma(x)(I+\M_\sigma),(I+\M_\sigma)} =\rho(X^*xX)$.
As $\sigma$ and $\rho$ are normal states on $\ddual{\C}$, this
equality is valid for $x\in\ddual{\C}$ as well; in particular, $1=
\rho(X^*p_\sigma X)= \sigma(p_\sigma)$.  For $x\in\C$, define
$\phi(x)=\rho(X^*p_\sigma x)$.  Then $\phi\neq 0$ because $\phi(X)=1$
and $\phi$ is an eigenfunctional with $s(\phi)=\rho$ and
$r(\phi)=\sigma$.  Hence $\rho\sim\sigma$.  Conversely, if
$\rho\sim\sigma$, then find a normalizer $v$ with $\rho(v^*v)=1$ and
$\sigma(x)= \rho(v^*xv)$ for every $x\in\C$.  Then
$\sigma(x)=\innerprod{\pi_\rho(x)v+\M_\rho,v+\M_\rho}$.  By Kadison's
Transitivity Theorem, there exists a unitary $V\in\C$ such that
$V^*(1+N_\rho)=v+\M_\rho.$ Then $\pi_\sigma = V\pi_\rho V^*.$
\end{proof}

\begin{theorem}\label{faithfulcompatable}
If $\X\subseteq \hat{\D}$ contains exactly one 
element from each equivalence class in $R(\C)$,
then $\pi=\bigoplus_{\rho\in\X} \pi_\rho$ on
$\H=\bigoplus_{\rho\in\X}\H_\rho$
is a faithful compatible representation of $\C$ on $\B(\H)$ and
$\pi(\D)''$ is an multiplicity-free atomic MASA in $\bh$.

Moreover, if $\A\subseteq (\C,\D)$ and $P$ is the support projection for $\rho$, 
then $\ker\tilde{\pi}|_{\ddual{\A}}= P^\perp\ddual{\A}$ and 
$\tilde{\pi}(\ddual{\A})$ is a CSL algebra.
\end{theorem}

\begin{proof}
For $\rho\in\X$, let $P_\rho\in\ddual{\C}$ be the support projection
of $\pi_\rho$.  Then for $\sigma\in\hat{\D}$,  
$$P_\rho
p_\sigma=\begin{cases} p_\sigma&\text{if }\sigma\sim \rho\\ 0&
\text{otherwise.}
\end{cases}$$
Thus $\tilde{\pi}(p_\sigma)$ is a minimal projection for every
$\sigma\in\hat{\D}$.
Since $$I_\H=\sum_{\rho\in\X}\tilde{\pi}(P_\rho)=\sum_{\rho\in
  \X}\sum_{\sigma\sim \rho} \tilde{\pi}(p_\sigma),$$ 
which is a sum of minimal projections,  $\pi(\D)''$ is an multiplicity-free
atomic MASA in $\bh$.
If, for each $\rho \in \X$, $E'_\rho:\B(\H_\rho)\rightarrow\pi_\rho(\D)''$ is the expectation 
on $\B(\H_\rho)$ induced by $E$, 
 the map  $E'=\bigoplus_{\rho\in\X}
E'_\rho$ is faithful and satisfies $E'\circ \pi=
\pi\circ E$ so $\pi$ is compatible.  For any $x\in \C$ such that 
$\pi(x^*x)=0$, 
$$0=E'(\pi(x^*x))=\pi(E(x^*x)),$$ 
so $E(x^*x)=0$ since $\pi$ is faithful on $\D$.  Thus, as $E$ is
faithful, we conclude that $x^*x=0$, hence $\pi$ is faithful on $\C$.

Suppose $\A\subseteq (\C,\D)$.
As $P=\sum_{\rho\in\X}P_\rho$, where $P_\rho$ is the support projection
for $\pi_\rho$, Proposition~\ref{specialrep} implies that
$P\in\ddual{\D}\subseteq \ddual{\A}$.  
Thus $\ker\tilde{\pi}|_{\ddual{\A}}= P^\perp\ddual{\A}$.

We claim that $\tilde{\pi}(\ddual{\A})$ is weak-$*$ closed in $\bh$.
By the Krein-Smullian Theorem, it suffices to prove the
(norm-closed) unit ball of $\tilde{\pi}(\ddual{\A})$ is weak-$*$ closed.
Suppose $(y_\lambda)$ is a net in $\tilde{\pi}(\ddual{\A})$
converging weak-$*$ to $y\in \bh$ with
with $\norm{y_\lambda}\leq 1$ for all $\lambda$.
For each $\lambda$, choose $x_\lambda\in\ddual{\A}$ with
$\tilde{\pi}(x_\lambda)=y_\lambda$.
Since $\tilde{\pi}$ is isometric on $P\ddual{\C}$, the net $(Px_\lambda)$
in $\ddual{\A}$ satisfies $\norm{Px_\lambda}\leq 1$ for every $\lambda$.
Hence a subnet $Px_{\lambda_\mu}$ converges weak-$*$ to some $x\in \ddual{\A}$.
Then $\tilde{\pi}(x)=y$, as required.

Since $\tilde{\pi}(\ddual{\A})$ contains the atomic MASA $\pi(\D)''$,
$\L:=\lat(\tilde{\pi}(\ddual{\A}))$ is an atomic CSL.
We claim that $\tilde{\pi}(\ddual{\A})=\alg\L$.  As $\alg\L$ is the
largest algebra whose lattice of invariant subspaces is $\L$, 
$\tilde{\pi}(\ddual{\A})\subseteq \alg\L$.

To obtain the reverse inclusion, first observe that for minimal
projections $p,q \in\pi(\D)''$, $q\alg\L p=q\tilde{\pi}(\ddual{\A})
p$.  Indeed, the subspaces $\overline{\alg\L p\H}$ and
$\overline{\tilde{\pi}(\ddual{\A})p\H}$ are each the smallest element
of $\L$ containing the range of $p$, so the two subspaces coincide.
Thus $q\overline{\alg\L p\H}=q\overline{\tilde{\pi}(\ddual{\A})p\H}$,
which yields the observation.

Let $\bbA$ be the set of minimal projections of $\pi(\D)''$.  Given
$Y\in\alg\L$, we may write $Y$ as the weak-$*$ convergent sum,
$Y=\sum_{q,p\in \bbA}qYp$.  As each $qYp\in\tilde{\pi}(\ddual{\A})$,
and $\tilde{\pi}(\ddual{\A})$ is weak-$*$ closed, $Y\in
\tilde{\pi}(\ddual{\A})$, as desired.
\end{proof}

Muhly, Qiu and Solel \cite[Theorem~4.7]{MR94i:46075} prove that
if $\C$ is nuclear and $\A \subset (\C,\D)$ with $\A$ triangular, then
the expectation $E|_{\A}$ is a homomorphism.  
The connection with
CSL algebras provided by Theorem~\ref{faithfulcompatable} allows us to
remove the hypothesis of nuclearity.

\begin{theorem}\label{ExpectHomomorph} If $\A$ is a triangular
subalgebra of the \cstardiag\ $(\C,\D)$, then $E|_{\A}$ is a
homomorphism of $\A$ onto $\D$. 
\end{theorem}

\begin{proof}
Let $\pi$ be the faithful compatible representation of $\C$ provided
by Theorem~\ref{faithfulcompatable}, and again write $\alg\L$ for
$\tilde{\pi}(\ddual{\A})$.  Theorem~\ref{faithfulcompatable} also
shows that $\tilde{\pi}(\ddual{\D})$ is a MASA in $\bh$.  We claim
that $\L$ is multiplicity free.  To show this, we prove that $\alg\L
\cap(\alg\L)^*=\tilde{\pi}(\ddual{\D})$.  Clearly,
$\tilde{\pi}(\ddual{\D})\subseteq \alg\L \cap(\alg\L)^*$.  To show the
reverse implication, suppose $q_1,q_2\in \tilde{\pi}(\ddual{\D})$ are
distinct nonzero minimal projections and $q_2(\alg\L) q_1\neq (0).$   It
suffices to show that 
\begin{equation}\label{reverseQ}
q_1(\alg\L) q_2 =(0).
\end{equation}
We may 
find $\rho,\sigma\in \hat{\D}$ so that $q_1=\tilde{\pi}(p_\rho)$ and
$q_2=\tilde{\pi} (p_\sigma)$, where $p_\rho$ and $p_\sigma$ are as in
Definition~\ref{asparisom}. Since $\pi(\A)$ is weak-$*$ dense in
$\alg\L$, we see that $p_\sigma\A p_\rho\neq (0)$.  Therefore,
$B_{\sigma,\rho}$ is nonzero on $\A$, so that there exists an
eigenfunctional $\phi\in \Eigone(\A)$ with $s(\phi)=\rho$ and
$r(\phi)=\sigma$.  

We claim that $\phi^*$ vanishes on $\A$.  Indeed,
as $q_1\neq q_2$, $\rho\neq\sigma$, so that we may find a normalizer
$v\in\A$ so that $(v^*v)(vv^*)=0$ and $\phi=[v,\rho]$.  Suppose to
obtain a contradiction, that $\phi^*=[v^*,\sigma]$ does not vanish on
$\A$, and let $y\in\A$ satisfy $\phi^*(y)\neq 0$.
Proposition~\ref{modmap} shows that $w:=v^*E(vy)$ is a nonzero element
of $\A$.  But since $v\in\A$ and $E(vy)\in\D$, we also have $w^*\in
\A$.  However, 
\begin{align*}  
    ww^*&=v^*E(vy)E(vy)^*v\leq \norm{E(vy)}^2v^*v,\\
    w^*w&=E(vy)^*vv^*E(vy)\leq \norm{E(vy)}^2vv^*.
\end{align*}
Thus $w$ is a
nonnormal element of $\A\cap \A^*$, violating triangularity of $\A$.

Therefore, $\phi^*$ vanishes on $\A$, so that $p_\rho\A
p_\sigma=(0)$.  Applying $\tilde{\pi}$, we obtain $q_1(\alg\L)
q_2=(0)$.  Therefore, $\alg\L$ is multiplicity free as desired.

 Each minimal projection $e$ in $\tilde{\pi}(\ddual{\D})$ is the
difference of elements of $\L$, and hence the compression $x\mapsto
exe$ is a homomorphism on $\alg\L$.  The extension of $\tilde{E}$ of
$E$ to all of \bh\ is faithful and normal, and is the sum of such
compressions.  Thus, $\tilde{E}$ is a homomorphism on $\alg\L$ and, by
restriction, on $\A$.
\end{proof}

We describe the maximal ideals of $\A$ and identify $\ker E|_{\A}$
in algebraic terms.

\begin{proposition}\label{maxidkerE}
Let $(\C,\D)$ be a \cstardiag\ and $\A\subseteq (\C,\D)$ be triangular.
The map $J\mapsto J\cap \D$ is a bijection between the proper
maximal ideals of $\A$ and the proper maximal ideals of $\D$.
Further, 
$$\ker E|_\A=\bigcap \{J\subseteq \A: J \text{ is a maximal ideal of } \A\}.
$$
\end{proposition}

\begin{proof}
Let $J$ be a proper maximal ideal of $\A$.  For any $x\in\C$, $E(x)$
belongs to the norm-closed convex hull of $\{gxg^{-1}: g \in\D,
g\text{ is unitary}\}$
(Theorem~\ref{exptHull}).
Since $J$ is a closed ideal, we see that
$E(J)\subseteq J\cap \D \subseteq E(J)$, so $E(J)=J\cap \D$.  
Since $J$ is proper, $J\cap \D\neq \D$.  Hence there exists
$\rho\in\hat{\D}$ such that $\ker\rho\supseteq J\cap \D$.  The unique
extension of $\rho$ to $\A$ is $\rho\circ E|_\A$, so we have
$J\subseteq \ker( \rho\circ E|_\A)$.  
Since $J$ is maximal, $J=\ker\rho\circ E|_\A$.  
Therefore, $J\cap \D=\ker\rho$, which is a proper maximal ideal of $\D$.

To show that the map $J\mapsto J\cap \D$ is a bijection, we need only
consider proper maximal ideals.
The previous paragraph shows that if $J_1$ is another proper maximal 
ideal of $\A$ and $J_1\cap \D=J\cap D$, then $J=J_1$.  
Also, if $K\subseteq \D$ is a proper maximal ideal, then $K=\ker\rho$ 
for some $\rho\in\hat{\D}$, and then $J:=\ker\rho\circ E_\A$ is a proper 
maximal ideal of $\A$ with $J\cap \D=K$.  
Thus the map is onto.  

For $x \in \A$, we have $x\in \ker E|_\A$ if and only if 
$x\in\ker\rho\circ E|_\A$ for every $\rho\in \hat{\D}$.
The map above shows this is equivalent to
$x\in\bigcap \{J\subseteq \A: J \text{ is a maximal ideal of } \A\}$.
\end{proof}

\section{Invariance under Diagonal-Preserving Isomorphisms}\label{S:Iso}

As an application of the results obtained so far, we will show that
coordinate systems are preserved under isomorphisms of algebras which
preserve the diagonal.  These results extend those for isometric
isomorphisms and we compare our results with them.

\begin{defn}\label{D:DiagPres}
For $i=1,2$, suppose $(\C_i,\D_i)$ are regular \cstar-inclusions and
that $\A_i$ are subalgebras with $\A_i\subseteq (\C_i,\D_i)$, i.e.,
$\A_i \subseteq \C_i$ is a norm-closed subalgebra with $\D_i \subset \A_i$.
We say that a (bounded) isomorphism $\theta:\A_1\rightarrow\A_2$ is
\textit{diagonal preserving} if $\theta(\D_1)=\D_2$.
\end{defn}
 
\begin{remark}{Remarks} \begin{enumerate}
\item If $\theta$ is diagonal preserving, then $\theta|_{\D_1}$ is a 
$*$-isomorphism of $\D_1$ onto $\D_2$. 
\item While we are interested in coordinates for modules, when
  studying invariance properties, we note that it suffices to consider
  subalgebras of regular \cstar-inclusions.  This is because of the
  well-known ``$2\times 2$ matrix trick.''  If $(\C,\D)$ is a regular
  \cstar-inclusion, so is $(M_2(\C),\D\oplus\D)$. For a
  $\D$-bimodule $\M\subseteq \C$, let $\fT(\M)$ be the subalgebra of
  $M_2(\C)$,
\[ \fT(\M):= \left\{ \begin{bmatrix} d_1 & m \\ 0 & d_2 \end{bmatrix} 
	: d_1,d_2 \in \D, m \in \M \right\}\] contained in
$(M_2(\C),\D \oplus\D)$.  An isomorphism of 
$\D_i$-bimodules $\M_i$ ($i=1,2$), that is, a bounded map
$\theta:\M_1\rightarrow\M_2$ together with an isomorphism
$\alpha:\D_1\rightarrow \D_2$ satisfying
$\theta(dme)=\alpha(d)\theta(m)\alpha(e)$ ($d,e\in\D_1,m\in\M$) can be
equivalently described as a diagonal preserving isomorphism of
$\fT(\M_1)$ onto $\fT(\M_2)$.
\end{enumerate}
\end{remark}

We have noted that $\theta^\#$ is a bicontinuous isomorphism from $\Eig(\A_2)$ 
onto $\Eig(\A_1)$ (Proposition~\ref{adjointmap}).
The next result shows that normalizing $\theta^\#$ pointwise gives an 
isomorphism of the norm-one eigenfunctionals.

\begin{thm}\label{bimodiso}
For $i=1,2$, let $(\C_i,\D_i)$ be regular \cstar-inclusions, let
$\A_i \subseteq (\C_i,\D_i)$, and suppose 
$\theta : \A_1 \to \A_2$ is a bounded diagonal-preserving isomorphism.
There exists a bicontinuous isomorphism of coordinate systems 
$\gamma:\Eigone(\A_1)\rightarrow \Eigone(\A_2)$ given by
\[
\gamma(\phi)=\frac{\phi \circ \theta^{-1}}{\|\phi \circ \theta^{-1}\|}.
\]
Moreover, if $\phi\in\Eigone(\A_1)$ is written as $[v,\rho]$,
then $(\rho\circ\theta^{-1})(\theta(v)^*\theta(v))$ is nonzero, and 
$$ \gamma(\phi)=[\theta(v),\rho\circ\theta^{-1}].$$
\end{thm}

When necessary for clarity, we use $\gamma_\theta$ to denote the
dependence of $\gamma$ on $\theta$.

\begin{remark}{Remark} 
As a special case, if $\theta: \sA_1 \to \sA_2$ is a contractive
isomorphism, then it is diagonal preserving, as its restriction to
$\D_1$ is then a $*$-homomorphism.  Thus, Theorem~\ref{bimodiso}
extends previous work for isometric isomorphisms of TAF algebras and
of subalgebras of (nuclear) groupoid \cstar-algebras;
see~\cite[Theorem~3]{MR91e:46078} and~\cite[Theorem~2.1]{MR95a:46080}.
\end{remark}
  
\begin{proof}[Proof of Theorem~\ref{bimodiso}]
Given $\phi\in\Eigone(\A_1)$, write $\phi=[v,\rho]$ where $v\in\A_1$
is an intertwiner and $\rho(v^*v)\neq 0$.   That 
$r(\gamma(\phi))=r(\phi)\circ\theta^{-1}$ and
$s(\gamma(\phi))=s(\phi)\circ\theta^{-1}$ follows from the definition
of $\gamma$.

Clearly $\theta(v)$ is a $\D_2$-intertwiner, so
$\theta(v)^*\theta(v)\in\D_2$.  For $\rho\in\hat{\D}_1$ with
$\rho(v^*v)\neq 0$,
\begin{align*}
(\rho\circ\theta^{-1})(\theta(v)^*\theta(v))&=
 \inf\{\norm{\theta(d)^*\theta(v)^*\theta(v)\theta(d)}:d\in\D_1, \rho(d)=1\}\\
&=\inf\{\norm{\theta(vd)}^2: d\in\D_1, \rho(d)=1\}\\  
&\geq\norm{\theta^{-1}}^{-2}\inf\{\norm{vd}^2:d\in\D_1, \rho(d)=1\}\\
&=\norm{\theta^{-1}}^{-2}\rho(v^*v) \neq 0.
\end{align*}  
Simple calculations show that $[\theta(v),\rho\circ\theta^{-1}]$ and
$\phi\circ\theta^{-1}/\norm{\phi\circ\theta^{-1}}$ are elements of 
$\Eigone_{\D_2}(\A_2)$ with the same range.
Since both are positive on $\theta(v)$, 
$$\gamma(\phi)=[\theta(v),\rho\circ\theta^{-1}].$$ 

This formula implies that
$\gamma(\phi_1\phi_2)=\gamma(\phi_1)\gamma(\phi_2)$ whenever
$\phi_1\phi_2$ is defined, so that 
$\gamma$ is an algebraic isomorphism of coordinate systems.

To show continuity, let $(\phi_\lambda)$ be a net in $\Eigone(\A_1)$
converging weak-$*$ to $\phi\in\Eigone(\A_1)$.  Write $\phi=[v,\rho]$
for some intertwiner $v\in\A_1$ and $\rho\in\hat{\D}_1$.  Since
$\rho(v^*v)\neq 0$, there exists $d\in\D_1$ and a neighborhood of
$\rho$, $G\subseteq \hat{\D}_1$, with $d \ge 0$ and
$\sigma(d^*v^*vd)=1$ for all $\sigma\in G$.  Since $\rho(d)\ne 0$,
$\phi=[vd,\rho]$.  Replacing $v$ with $vd$, we may assume that
$\sigma(v^*v)=1$ for every $\sigma$ in $G\subseteq\hat{\D}_1$, a
neighborhood of $\rho$.

Since $\rho_\lambda:=s(\phi_\lambda)$ converges weak-$*$ to $s(\phi)=\rho$, by
deleting the first part of the net, we may assume that
$\rho_\lambda\in G$ and  $\phi_\lambda(v)\neq 0$ for all
$\lambda$.  By
Proposition~\ref{EigTwist}, there exist  scalars $t_\lambda\in\bbT$
such that $\phi_\lambda=[t_\lambda v,\rho_\lambda].$
Since
$$\overline{t}_\lambda= 
[t_\lambda v,\rho_\lambda](v)=\phi_\lambda(v)\rightarrow\phi(v)=
1,$$ and $[\theta(v),\rho_\lambda\circ\theta^{-1}]$ converges weak-$*$ to
$[\theta(v),\rho\circ\theta^{-1}]$, 
we conclude 
$$\gamma(\phi_\lambda)=[t_\lambda\theta(v),\rho_\lambda\circ\theta]
=\overline{t}_\lambda[\theta(v),\rho_\lambda\circ\theta^{-1}]
\rightarrow [\theta(v),\rho\circ\theta^{-1}]=\gamma(\phi).$$
Thus, $\gamma$ is continuous.
Similarly, $\gamma^{-1}$ is continuous.
\end{proof}

We give several applications of Theorem~\ref{bimodiso} to isomorphisms of  
subalgebras.
The following corollary is immediate.

\begin{cor}\label{bimodisocor}
For $i=1,2,3$, let $(\C_i,\D_i)$ be regular \cstar-inclusions, and let
$\A_i \subseteq (\C_i,\D_i)$ be norm-closed algebras.  For $j=1,2$,
let $\theta_{j}:\D_j\rightarrow \D_{j+1}$ be bounded diagonal
preserving isomorphisms.  Then
$\gamma_{\theta_2\circ\theta_1}=\gamma_{\theta_2}\circ\gamma_{\theta_1}.$
\end{cor}

\begin{defn}\label{cocycledefn}
Suppose $(\C,\D)$ is a regular \cstar-inclusion and $\A \subseteq
(\C,\D)$.  By a \textsl{(1-)cocycle} on $\Eigone(\A)$, we mean a map
$c : \Eigone(\A) \to \bbC_*$, the group of nonzero complex numbers
under multiplication, satisfying, for all composable elements
$\phi,\psi \in \Eigone(\A)$,
\[ c(\phi \psi)=c(\phi) \cdot c(\psi). \]
\end{defn}

\begin{corollary}\label{EigoneUnique}
Suppose $(\C,\D)$ is a  regular \cstar-inclusion, 
$\A\subseteq (\C,\D)$ is a  norm-closed algebra, and
$\theta:\A\rightarrow \A$ is a bounded automorphism fixing $\D$
elementwise.  Assume further that $\gamma_\theta$ is the identity map
on $\Eigone(\A)$.  
Then, for all $\phi\in\Eigone(\A)$,
\begin{equation}\label{AMC}
\phi(\theta(x)) =\phi(x)\norm{\phi\circ\theta}.
\end{equation} 
If $c:\Eigone(\A)\rightarrow \bbR$ is defined by
$c(\phi)=\norm{\phi\circ\theta}$, then $c$ is a 
positive cocycle on $\Eigone(\A)$.  

If, in addition, $\theta$ is isometric, then $\theta=\rm{id}_\A$. 
\end{corollary}

\begin{remark}{Remark}
In essence, \eqref{AMC} shows that $\theta$ is given by
multiplication by a cocycle on $\Eigone(\A)$.
Thus, we will call $\theta$ a \textit{cocycle automorphism}.
\end{remark}

\begin{proof}
Formula \eqref{AMC} follows immediately from Theorem~\ref{bimodiso} 
applied to $\theta^{-1}$. 
If $\theta$ is isometric, then \eqref{AMC} shows that $\phi(x)=\phi(\theta(x))$ 
for all $\phi \in \Eigone(\A)$.
By Theorem~\ref{Dinclusion}, $\Eigone(\A)$ separates points
and so $x=\theta(x)$ for all $x \in \A$.

It remains to show that $c$ is a cocycle.
Observe that if $v\in\C$ is a $\D$-intertwiner, then $v^*\theta(v)\in\D$.  
Indeed, for a self-adjoint $d\in\D$, using Remark~\ref{nostar}(2), we may find
$d'\in\D$ with $d'$ also self-adjoint and so that $vd=d'v$.  
Then $v^*\theta(v)d=v^*d'\theta(v)=dv^*\theta(d)$, so that $v^*\theta(v)$ 
commutes with the self-adjoint elements of $\D$ and hence belongs to $\D$.

Finally, suppose for $i=1,2$, that $\phi_i=[v_i,\rho_i] \in \Eigone(\A)$ 
such that the product $\phi_1\phi_2$ is defined.  
Then $\phi_1\phi_2=[v_1v_2,\rho_2]$, and 
$\rho_1=r([v_2,\rho_2])$, i.e., $\rho_1(x)=\rho_2(v_2^*xv_2)/\rho_2(v_2^*v_2)$. 
Using \eqref{AMC} and these facts, we have
\begin{align*}
c(\phi_1\phi_2) 
&=\frac{[v_1v_2,\rho_2](\theta(v_1v_2))}{[v_1v_2,\rho_2](v_1v_2)}\\
&=\frac{\rho_2(v_2^*v_1^*\theta(v_1) 
      \theta(v_2))}{\rho_2(v_2^*v_1^*v_1v_2)}\\
&=\frac{\rho_2(v_2^*v_2)\rho_2(v_2^*v_1^*\theta(v_1) 
      \theta(v_2))}{\rho_2(v_2^*v_2)\rho_2(v_2^*v_1^*v_1v_2)}\\
&=\frac{\rho_1(v_1^*\theta(v_1))\rho_2(v_2^*\theta(v_2))}
{\rho_2(v_2^*v_1^*v_1v_2)}\\
&=\frac{\rho_1(v_1^*\theta(v_1))\rho_2(v_2^*\theta(v_2))}
{\rho_1(v_1^*v_1)\rho_2(v_2^*v_2)}
 =c(\phi_1)c(\phi_2)
\end{align*}
\end{proof}

For our next application, we need two technical lemmas.

\begin{lemma}\label{starmap}
For $i=1,2$, let  $(\C_i,\D_i)$ be regular \cstar-inclusions and let 
$\B_i \subseteq (\C_i,\D_i)$ be selfadjoint.
If $\gamma : \Eigone(\B_1) \to \Eigone(\B_2)$ is an isomorphism of 
coordinate systems, then, for all $\phi \in \Eigone(\B_1)$,
$\gamma(\phi^*)=\gamma(\phi)^*$.
\end{lemma}

\begin{proof}
Since each $\B_i$ is selfadjoint, $\tau^* \in \Eigone(\B_i)$ if 
$\tau \in \Eigone(\B_i)$, and so $\gamma(\phi^*)$ and $\gamma(\phi)^*$ are
defined.
As $\gamma(\phi^*)$ and $\gamma(\phi)^*$ have the same range and domain,
by Corollary~\ref{multiples}, there is $\lambda \in \bbT$ with
$\gamma(\phi^*)=\lambda \gamma(\phi)^*$.
As $\phi \cdot \phi^* = r(\phi)$, we have
\[
\gamma(r(\phi))=\gamma(\phi \cdot \phi^*) = \gamma(\phi) \cdot \gamma(\phi^*)
	= \gamma(\phi) \cdot (\lambda \gamma(\phi)^*) = \lambda r(\gamma(\phi)).
\]
As $\gamma(r(\phi))=r(\gamma(\phi)$, the result follows.
\end{proof}

\begin{lemma}\label{regularguy}
For $i=1,2$, let  $(\C_i,\D_i)$ be regular \cstar-inclusions, let 
$\A_i \subseteq (\C_i,\D_i)$ be regular, and let $\pi : \C_1 \to C_2$
be a $*$-isomorphism with $\pi(\D_1)=\D_2$.
If $\gamma_\pi$ maps $\Eigone(\A_1) $ into $\Eigone(\A_2)$, then 
$\pi(\A_1)\subseteq \A_2$.
\end{lemma}

\begin{proof}
Since each $\A_i$ is regular, it is enough to show that for an
intertwiner $v\in\A_1$, we have $\pi(v)\in\A_2$.  Clearly, $\pi(v)$ is
an intertwiner in $\C_2$, and it is easy to see that the closed
$\D_2$-bimodule generated by $\pi(v)$ is isometrically isomorphic to
$\overline{|\pi(v)|\D_2}$ via the map $\pi(v)d\mapsto |\pi(v)|d$.  Given
$\rho\in\hat{\D}_2$ with $\rho(\pi(v^*v))\neq 0$, we have
$\gamma([v,\rho\circ \pi^{-1}])= [\pi(v),\rho]\in \Eigone(\A_2)$.
Applying Proposition~\ref{modmap}, we conclude that there exists
$d\in\D_2$ with $d \ge 0$, $\rho(d)=1$, and $\pi(v)d\in\A_2$. Given
$\eps >0$, let $f_{\eps} \in\D$ be chosen so that $0\leq f_{\eps} \leq
I$,\, 
$\hat{f}_{\eps}$ is
compactly supported in $G:=\{ \rho\in\hat{\D}_2: \rho(|\pi(v)|)>0\}$
and $\norm{\pi(v)- \pi(v)f_{\eps}}<\eps$.  
Compactness of $F_{\eps}:=\overline{\rm{supp}\,f_{\eps}}$ ensures that
there exist $d_1,\dots, d_n\in\D_2$ such that $\pi(v)d_i\in\A_2$ and
$\rho(\sum_{i=1}^n d_i) > 0$ for every $\rho\in F_{\eps}$.  Hence there
exists an element $g\in\D_2$ such that for every $\rho\in F_{\eps}$,
$\rho(g\sum_{i=1}^n d_i)=1$.  Then 
$\pi(v)f_{\eps}= \sum_{i=1}^n \pi(v)d_igf_{\eps}\in
\A_2$.  Letting $\eps\rightarrow 0$, we conclude that $\pi(v)\in \A_2$
as well.
\end{proof}

Recall that a subalgebra  $\A\subseteq \C$ is said to
be \textit{Dirichlet} if $\A+\A^*$ is norm dense in $\C$.  This
implies that $\Eigone(\C)=\Eigone(\A) \cup \Eigone(\A)^*$.  To see
this, let $\phi \in \Eigone(\C)$.  By density, $\phi$ does not vanish
on one of $\A$ or $\A^*$ and hence is in either $\Eigone(\A)$ or
$\Eigone(\A)^*$.  

\begin{thm}\label{trialgcor}
For $i=1,2$, let  $(\C_i,\D_i)$ be regular \cstar-inclusions and let 
$\A_i$ be Dirichlet subalgebras with 
$\D_i\subseteq \A_i\subseteq \C_i$.
Consider the following statements.
\begin{enumerate}
\item \label{rone} $\A_1$ and $\A_2$ are isometrically isomorphic.
\item \label{rtwo} There exists a bounded isomorphism
  $\theta:\A_1\rightarrow\A_2$ such that $\theta(\D_1)=\D_2$.
\item \label{rthree} $\Eigone(\A_1)$ and $\Eigone(\A_2)$ are
  isomorphic coordinate systems.
\item \label{rfour} $\Eigone(\C_1)$ and $\Eigone(\C_2)$ are
  isomorphic twists and the isomorphism maps $\Eigone(\A_1)$ 
  onto $\Eigone(\A_2)$.
\end{enumerate}
Then statement~$(n)$ implies statement $(n+1)$, $n=1,2,3$.

If, in addition, $(\C_i,\D_i)$ are \cstardiag s and $\A_i$ are regular, then 
the four statements are equivalent.
\end{thm}

\begin{proof}
That \eqref{rone} implies \eqref{rtwo} is obvious.  To show that
\eqref{rtwo} implies \eqref{rthree}, let $\alpha=\theta|_{\D_1}$, a
$*$-isomorphism, and apply Theorem~\ref{bimodiso}.

To show \eqref{rthree} implies \eqref{rfour}, suppose $\gamma :
\Eigone(\A_1) \to \Eigone(\A_2)$ is a continuous isomorphism.  We
extend $\gamma$ to a map, call it $\delta$, on all of $\Eigone(\C_1)$
by mapping $\phi \in \Eigone(\A_1)^*=\Eigone(\A_1^*)$ to
$\gamma(\phi^*)^*$.  By Lemma~\ref{starmap}, $\delta$ is well-defined
on $\Eigone(\A_1 \cap \A_1^*)$, and hence well-defined on all of
$\Eigone(\C_1)$.

Since the adjoint map is continuous, $\delta$ is a continuous
homomorphism restricted to $\Eigone(\A_1)$ or to $\Eigone(\A_1^*)$.
By Corollary~\ref{extensionunique}, these are open sets in
$\Eigone(\C_1)$, and so $\delta$ is continuous on their union,
$\Eigone(\C_1)$.  Applying the same argument to $\delta^{-1},$ we see
that $\delta^{-1}$ is continuous as well. 

The restriction of $\delta$ to $\Eigone(\A_1^*)$ or to $\Eigone(\A_1)$ 
is a homomorphism.
To show that $\delta$ is a homomorphism, fix $\phi, \psi\in \Eigone(\A_1)$.
We claim that if $\phi^*\psi$ is defined, then
$\delta(\phi^*\psi)=\delta(\phi^*)\delta(\psi)$ and
if $\phi\psi^*$ is defined, then
$\delta(\phi\psi^*)=\delta(\phi)\delta(\psi)^*$.
We only show the first equality; the proof of the second is similar.

As $\delta$ maps $\hat{\D}_1$ to $\hat{\D}_2$, we have,
for $\phi \in \Eigone(\A_1)$, 
\begin{equation}\label{deltarange}
\delta(\phi\phi^*)=\gamma(r(\phi))=r(\gamma(\phi))=
\gamma(\phi)\gamma(\phi^*).
\end{equation}

Let $\eta=\phi^*\psi$.  If $\eta\in\Eigone(\A_1)$, we have
$\phi\eta=\phi\phi^*\psi=\psi$, and hence
$\gamma(\phi)\gamma(\eta)=\gamma(\psi)$.  Multiplying each side of
this equality by $\gamma(\phi^*)=\delta(\phi)^*$ and
using~\eqref{deltarange} yields
$\delta(\phi^*\psi)=\delta(\phi)^*\delta(\psi)$.  If
$\eta\in\Eigone(\A_1^*)$, then apply the previous argument to
$\eta^*=\psi^*\phi$ and take adjoints to obtain the equation.  Thus
$\delta$ is a homomorphism and  
part~\eqref{rfour} holds.

Finally, if $(\C_i,\D_i)$ are \cstar-diagonals and \eqref{rfour}
holds, then Kumjian's theorem implies that there is an isomorphism
$\Phi : \C_1 \to \C_2$ with $\Phi(\D_1)=\D_2$.  Since $\Phi$ is
induced by an isomorphism of twists mapping $\Eigone(\A_1)$ onto
$\Eigone(\A_2)$, Lemma~\ref{regularguy} shows that $\Phi(\A_1) = \A_2$.
Thus, we obtain \eqref{rone} with $\Phi|_{\A_1}$ the isometric
isomorphism.
\end{proof}

\begin{remark}{Remark}
Theorem~\ref{trialgcor} is related to a result of Muhly, Qiu and
Solel, \cite[Theorem~2.1]{MR95a:46080}.  Their result, expressed in
terms of eigenfunctionals, takes the following form.  When the $\C_i$
are nuclear and $(\C_i,\D_i)$ are \cstardiag s, the following are
equivalent, for subalgebras $\A_i\subseteq (\C_i,\D_i)$ so that $\A_i$
generates $\C_i$ as \cstar-algebras (with \textsl{no} Dirichlet
hypothesis):
\begin{enumerate}
\item[(a)] there is an isometric isomorphism $\theta:\A_1\rightarrow \A_2$
  (necessarily, $\theta(\D_1)=\D_2$);
\item[(b)] there is a coordinate system isomorphism,
  $\gamma:\Eigone(\A_1)\rightarrow \Eigone(\A_2)$ which extends to a
  coordinate system isomorphism $\gamma': \Eigone(\C_1)\rightarrow
  \Eigone(\C_2)$;
\item[(c)] there is a $*$-isomorphism $\tau:\C_1\rightarrow \C_2$
  such that $\tau|_{\A_1}$ is an isomorphism of $\A_1$ onto $\A_2$.
\end{enumerate}
Theorem~\ref{trialgcor} extends the Muhly-Qiu-Solel result to 
not-necessarily-isometric diagonal preserving  
isomorphisms, assuming the Dirichlet 
condition instead of the hypothesis in (b) that $\gamma$ extends to an
isomorphism of $\Eigone(\C_1)$ onto $\Eigone(\C_2)$.
Example~\ref{cost} shows that in the absence of the Dirichlet
condition, isomorphisms of coordinate systems need
not extend to isomorphisms of the enveloping twists.  Thus,  the
hypothesis that $\gamma$ extends in (b) is essential. 

Also, since Theorem~\ref{trialgcor} did not use the Spectral Theorem
for Bimodules~\cite[Theorem~4.1]{MR94i:46075}, we do not need 
the $\C_i$ to be nuclear.
In Theorem~\ref{easyMQS}, we prove the full Muhly-Qiu-Solel 
result, without requiring nuclearity.
\end{remark}

\begin{remark}{Example} \label{cost} 
Without the Dirichlet hypothesis, an isomorphism of coordinate systems
need not induce an isometric isomorphism of the algebras.
Let $(\C,\D)=(M_4(\bbC),\D_4)$, where $\D_4$ is the algebra of 
diagonal matrices.

Let $\A_1=\A_2$ be the algebra 
$$\A:=\left\{ \begin{pmatrix} t_{11}&0&t_{13}&t_{14}\\
0&t_{22}&t_{23}&t_{24}\\ 0&0&t_{33}&0\\ 0&0&0& t_{44}\end{pmatrix}:
t_{ij}\in\bbC \right\}. $$
The automorphism 
$$\begin{pmatrix} t_{11}&0&t_{13}&t_{14}\\
0&t_{22}&t_{23}&t_{24}\\ 0&0&t_{33}&0\\ 0&0&0& t_{44}\end{pmatrix}
 \mapsto \begin{pmatrix} t_{11}&0&t_{13}&t_{14}\\
0&t_{22}&t_{23}&-t_{24}\\ 0&0&t_{33}&0\\ 0&0&0& t_{44}\end{pmatrix}
$$
is not isometric, and induces an automorphism of $\Eigone(\A)$ which
does not extend to an automorphism of $\Eigone(\C)$.
\end{remark}

The Dirichlet condition can be removed if one assumes a continuous
section from $R(\C_i)$ into $\Eigone(\C_i)$.  
Since TAF algebras always admit such sections,
the following result generalizes~\cite[Theorem~7.5]{MR94g:46001}.

\begin{theorem}\label{bigSteve}  Suppose for $i=1,2$, $(\C_i,\D_i)$
  are \cstardiag s and $\A_i\subseteq (\C_i,\D_i)$ are norm closed
  subalgebras such that $\A_i$ generates $\C_i$ as a
  \cstaralg.   Consider the following statements:
\begin{enumerate}
\item \label{bone} $\A_1$ and $\A_2$ are isometrically isomorphic;
\item \label{btwo} there exists a bounded isomorphism
  $\theta:\A_1\rightarrow\A_2$ such that $\theta(\D_1)=\D_2$;
\item \label{bthree} $R(\A_1)$ and $R(\A_2)$ are isomorphic topological
  binary relations,
\item \label{bfour} $R(\C_1)$ and $R(\C_2)$ are isomorphic topological 
  equivalence\\ relations, and the isomorphism maps $R(\A_1)$ onto $R(\A_2)$.
\end{enumerate}
Then for $i=1,2,3$, statement $(i)$ implies statement $(i+1)$.  

If, in addition, $\A_i$ are regular and there exist continuous sections 
$h_i: R(\C_i)\rightarrow \Eigone(\C_i)$, then the statements are equivalent.
\end{theorem}

\begin{proof}
That \eqref{bone}$\Rightarrow$\eqref{btwo} is obvious and
\eqref{btwo}$\Rightarrow$\eqref{bthree} follows as in the proof of 
Theorem~\ref{trialgcor}.  

Suppose \eqref{bthree} holds.  By Proposition~\ref{envelope}, the
isomorphism of $R(\A_1)$ onto $R(\A_2)$ extends uniquely to an
isomorphism of $R(\C_1)$ onto $R(\C_2)$, so \eqref{bfour} holds.

To complete the proof, we show that, when the $\A_i$ are regular and 
there exist continuous sections $h_i:R(\C_i)\rightarrow \Eigone(\C_i)$,
then \eqref{bfour} implies \eqref{bone}.
The existence of the sections and \eqref{bfour} gives
a coordinate system isomorphism $\gamma$ of $\Eigone(\C_1)$ and 
$\Eigone(\C_2)$ such that $\gamma|_{\Eigone(\A_1)}$ is an isomorphism of 
$\Eigone(\A_1)$ onto $\Eigone(\A_2)$.  
By Kumjian's theorem, there is a (regular) $*$-isomorphism $\pi$ of 
$(\C_1,\D_1)$ onto $(\C_2,\D_2)$.

Finally, Lemma~\ref{regularguy} implies $\pi(\A_1)=\A_2$. 
\end{proof}

\begin{remark}{Remarks} A modification of Example~\ref{cost} shows that
  in general, there may exist a section for $R(\A)$ which cannot be
  extended to a section of $R(\C)$.  Thus, one cannot replace the
  hypothesis of a section for $R(\C_i)$ with a hypothesis of a section
  for $R(\A_i)$ in Theorem~\ref{bigSteve}.

To drop the Dirichlet condition in Theorem~\ref{trialgcor} without
assuming the existence of a continuous section, we would need to
replace $\gamma$ with a new isomorphism from $\Eigone(\A_1)$ to
$\Eigone(\A_2)$ that could extend to the twist of $\C_1$, which would
be larger than $\Eigone(\A_1) \cup \Eigone(\A_1)^*$.  In
particular, while isometric bimodule maps on finite relations are
always $*$-extendible (a key ingredient in~\cite{MR95a:46080}), this
is not true for general bimodule maps~\cite[Theorem~1.2,
Proposition~1.4]{MR92h:47058}.  Indeed, based on~\cite{MR92h:47058},
there will homological obstructions to be considered.
\end{remark}

\section{Invariance under General Isomorphisms}
\label{isocoordsys}

In this section and the next, we come to the core of the paper, the study of 
bounded isomorphisms of 
triangular algebras which need not map the diagonal to the diagonal.
The principal result of this section is Theorem~\ref{algisom}, which shows that 
such isomorphisms induce \textit{algebraic} isomorphisms of the corresponding 
coordinate systems.
We would particularly like to know if this algebraic isomorphism is always 
continuous.  
We can prove that it is continuous in certain cases, Theorem~\ref{algpresgamma}, 
extending results of Donsig-Hudson-Katsoulis.

\begin{remark*}{Standing Assumption for Section~\ref{isocoordsys}}
\label{standing}
For $i=1,2$, let $(\C_i,\D_i)$ be \cstardiag s, and 
let $\A_i\subseteq (\C_i,\D_i)$ be (norm-closed) triangular subalgebras.  
Let $E_i : \C_i \to \D_i$ be the unique conditional expectations.
Suppose $\theta:\A_1\rightarrow \A_2$ is a bounded isomorphism.
\end{remark*}

Our first task is to show that $\theta$ induces
an \textit{algebraic} isomorphism 
$\gamma:\Eigone(\A_1)\rightarrow \Eigone(\A_2)$.  
We have been unable to show that $\gamma$ is continuous in general.  
Since we do not assume $\theta(\D_1)=\D_2$, it is not possible to use
Theorem~\ref{bimodiso}, so we proceed along different lines.

\begin{defn}\label{alphaDefn}
Define $\alpha:\D_1\rightarrow\D_2$ by $\alpha(d)=E_2(\theta(d))$.
\end{defn}

\begin{prop}\label{alphaMap}
The map $\alpha$ is a $*$-isomorphism of $\D_1$ onto $\D_2$.
\end{prop}

\begin{proof}
Theorem~\ref{ExpectHomomorph} shows $E_2 |_{\A_2}$ is a homomorphism; 
hence $\alpha$ is a homomorphism.  Since any algebraic
isomorphism of commutative \cstaralg s is a $*$-isomorphism, it
suffices to  show that $\alpha$ is bijective.
  
Since $E_1$ is idempotent, $\D_1\cap \ker E_1=\{0\}$, so
$\A_1=\D_1+\ker E_1|_{\A_1}$ is a direct sum decomposition.  By
Proposition~\ref{maxidkerE}, $\theta(\ker E_1|_{\A_1})=\ker
E_2|_{\A_2}$, so that we have two direct sum decompositions of $\A_2$:
$$\A_2=\D_2 +\ker E_2|_{\A_2}=\theta(\D_1)+\ker E_2|_{\A_2}.$$
Therefore, $\ker E_2|_{\theta(\D_1)}=\ker\alpha$ is trivial.  If
$d\in\D_2$ we may write $d= x+y$ where $x\in\theta(\D_1)$ and
$y\in\ker E_2$.  Then $E_2(x)=d$, so $\alpha$ is onto.
\end{proof}

Given Banach spaces $X$ and $Y$, and a bounded linear map
$R:X\rightarrow Y$, the double transpose map
$\ddual{R}:\ddual{X}\rightarrow\ddual{Y}$ is a norm-continuous extension
of $R$ which is also
$\sigma(\ddual{X},\dual{X})$--$\sigma(\ddual{Y},\dual{Y})$ continuous.   
In light of our standing assumptions,
$\ddual{\theta}:\ddual{\A_1}\rightarrow \ddual{\A_2}$ and
$\ddual{\alpha}:\ddual{\D_1}\rightarrow \ddual{\D_2}$ are also
isomorphisms. 
Similarly, the $\ddual{E_i}$ are homomorphisms of $\ddual{\A_i}$ onto 
$\ddual{\D_i}$. 
For notational ease, we will sometimes identify the double
transpose map with the original map.
Thus we often write $\theta$ or $\alpha$ instead of 
$\ddual{\theta}$ or $\ddual{\alpha}$.

\begin{remark}{Remark}
Notice that $\A_2$ may be regarded as a $\D_1$-bimodule in two ways:
for $d, e\in\D_1$ and $x\in \A_2$, define $d\cdot_\alpha x \cdot_\alpha
e:= \alpha(d)x\alpha(e)$ and $d\cdot_\theta x\cdot_\theta e:=
\theta(d)x\theta(e)$.  When these two modules are (boundedly)
isomorphic, methods similar to those used in the proof of
Theorem~\ref{bimodiso} show there exists an isomorphism of the
coordinate systems $\Eigone(\A_1)$ and $\Eigone(\A_2)$.
Unfortunately, we do not
know in general whether the $\alpha$ and $\theta$ module actions of
$\D_1$ on $\A_2$ are isomorphic.  However, there is enough structure
present to show that these modules are ``virtually'' isomorphic. 
\end{remark}

\begin{defn}\label{STdefn}
Let $G=\U(\D_1)$ be the group of unitary elements of $\D_1$, regarded
as a discrete abelian group and fix, once and for all, an invariant
mean $\Lambda$ on $G$.

Define elements $S, \, T\in\ddual{\A_2}$ by
requiring that for every $f\in\dual{\A_2}$,
\[ f(S)=\underset{g\in G}{\Lambda} f(\theta(g)\alpha(g^{-1}))
    \quad \text{and}\quad
   f(T)=\underset{g\in G}{\Lambda}  f(\alpha(g)\theta(g^{-1})).\]
\end{defn}

Notice that $S\in \chull \{\theta(g)\alpha(g^{-1}):g\in G\}$ and
$T\in\chull\{\alpha(g)\theta(g^{-1}):g\in G\}$, where $\chull Z$ is the
$\sigma(\ddual{\A_2},\dual{\A_2})$-closed convex hull of the set $Z$.
(We implicitly embed $\A_2$ into $\ddual{\A_2}$ using the canonical inclusion.)

We next collect some properties of $S$ and $T$.  Of particular interest
to us is the fact that they intertwine $\alpha(\D_1)$ and $\theta(\D_1)$.

\begin{proposition}\label{STFacts}
For $S$ and $T$ as above, we have
\begin{enumerate}
\item For every $d\in\D_1$,
$$ T\theta(d)=\alpha(d)T\quad \text{and}\quad S\alpha(d)=\theta(d)S.$$
\item $\ddual{E_2}(S)=I=\ddual{E_2}(T)$  and
$\ddual{E_1}(\theta^{-1}(S))=I=\ddual{E_1}(\theta^{-1}(T))$.
\item Given $\rho\in\hat{\D}_1$, let $p=p_\rho$
  (see Definition~\ref{asparisom}). Then
\begin{align*}
\alpha(p)T &= \alpha(p)T\theta(p)=T\theta(p)=\alpha(p)\theta(p)\quad \text{and}\\
\theta(p)S &=\theta(p)S\alpha(p)=S\alpha(p)=\theta(p)\alpha(p).
\end{align*}
\item For all $x\in\A_2$ and for all $\phi\in\Eigone(\A_1)$,
$\phi(\theta^{-1}(STxST))=\phi(\theta^{-1}(x)).$ 
\end{enumerate}
\end{proposition}

\begin{proof}
If $h\in G=\U(\D_1)$, then $T\theta(h)=\alpha(h)T$ follows from the
invariance of $\Lambda$.
Indeed, for every $f\in\dual{\A_2}$ we have,
\begin{align*}
f(\alpha(h)T)&=\Lambda_g (f\cdot\alpha(h)) (\alpha(g)\theta(g)^{-1})\\
&=\Lambda_g f(\alpha(hg)\theta(hg)^{-1}\theta(h))\\
&=\Lambda_g f(\alpha(g)\theta(g^{-1})\theta(h))\\
&=\Lambda_g (\theta(h)\cdot f) (\alpha(g)\theta(g^{-1}))\\
&=f(T\theta(h)).
\end{align*}
Since $\dual{\A_2}$ separates points of $\ddual{\A_2}$, we see that
$T\theta(h)=\alpha(h)T$ for every $h\in G$.  But the  span of
$G$ is norm dense in $\D_1$, which yields $ T\theta(d)=\alpha(d)T$ for
$d\in\D_1$.  The proof that $\theta(d)S=S\alpha(d)$ is similar.

To prove $\ddual{E_1}(\theta^{-1}(S))=I$, first observe that weak-$*$ 
continuity of $\ddual{E_1}$ (and $\ddual{{\theta^{-1}}}$) implies 
\[ \ddual{E_1}(\theta^{-1}(S))\in \ddual{E_1}\chull
\{g\theta^{-1}(\alpha(g^{-1}))\}=\chull
\{E_1(g\theta^{-1}(\alpha(g^{-1})))\}. \]
Modifying the proof of Proposition~\ref{alphaMap} yields 
$\alpha^{-1}=E_1\circ\theta^{-1}$ and $E_1(g\theta^{-1}(\alpha(g^{-1})))=I$
for each $g\in G$, and the equality follows.
The remaining equalities in part (2) have similar proofs.

For part (3), the first two equalities follow from statement (1),
as $p$ is a $\sigma(\ddual{\A_1},\dual{\A_1})$-limit of elements of
$\D_1$.  For the third equality, first observe that for $g\in G$,
Theorem~\ref{combineeigen} gives
$\alpha(p)\alpha(g)=\rho(g)\alpha(p)$; similarly,
$\theta(g)^{-1}\theta(p)=\rho(g^{-1})\theta(p).$ Hence,
$$\alpha(p)\alpha(g)\theta(g)^{-1}\theta(p)=\alpha(p)\theta(p).$$
Since $T\in\chull\{\alpha(g)\theta(g)^{-1} : g \in G \}$, we find that
\[ \alpha(p)T\theta(p) \in
   \chull\{\alpha(p)\alpha(g)\theta(g)^{-1}\theta(p) : g \in G \}.\]
This set is a singleton, so
$\alpha(p)T\theta(p)=\alpha(p)\theta(p)$.  The proofs of the
equalities involving $S$ are similar.  

For part (4), fix $\phi\in\Eigone(\A_1)$, and let $q$ and
$p$ be the minimal projections in $\ddual{\D_1}$ corresponding to
$r(\phi)$ and $s(\phi)$, respectively.  Part (1) implies that $p$ and
$q$ commute with $\theta^{-1}(ST)$ and by part (2), we have
$r(\phi)(\theta^{-1}(ST))= r(\phi)(E_1(\theta^{-1}(ST))) =
r(\phi)(I)=1$.  Hence by Proposition~\ref{combineeigen},
$q\theta^{-1}(ST)=q$.  Likewise $p\theta^{-1}(ST)=p$.  As, again by
Proposition~\ref{combineeigen}, $\phi(a) = \phi(q a p)$, we have
\begin{align*}
\phi(\theta^{-1}(STxST))
& = \phi(q\theta^{-1}(ST)\theta^{-1}(x)\theta^{-1}(ST))p)\\
& = \phi(q\theta^{-1}(x)p) =\phi(\theta^{-1}(x)),
\end{align*}
as desired.
\end{proof}

We now obtain a bijective mapping between the eigenfunctionals of $\A_1$ and 
those of $\A_2$.

\begin{proposition}\label{simpleeigenmap}
For $\phi \in \Eig(\A_1)$, let
$$f=T\cdot (\phi\circ\theta^{-1})\cdot S.$$
Then $f$ is an eigenfunctional for $\A_2$ with $r(f)=r(\phi)\circ
  \alpha^{-1}$ and $s(f)=s(\phi)\circ\alpha^{-1}$.

Moreover, $\phi\circ\theta^{-1}=S\cdot f\cdot T$.
\end{proposition}

\begin{proof}
For clarity, let $\psi = \phi\circ\theta^{-1}$.
For all $d, e\in\D_1$ and $x\in\A_2$,
\begin{align*} f(\alpha(d)x\alpha(e))
&= \psi(S\alpha(d)x\alpha(e)T)\\
&= \psi(\theta(d)SxT\theta(e))\\
&= r(\phi)(d)\psi(SxT)s(\phi)(e)\\
&= (r(\phi)\circ\alpha^{-1})(\alpha(d))\,f(x)\,(s(\phi)\circ\alpha^{-1})
(\alpha(e)),
\end{align*}
showing $f$ is an eigenfunctional with the claimed range and source.

The last equality follows from part (4) of Proposition~\ref{STFacts}.
\end{proof}

We  now show the existence of an algebraic isomorphism between
the coordinate systems which is the non-diagonal preserving analog of
Theorem~\ref{bimodiso}.

\begin{thm}\label{algisom}
The map $\gamma : \Eigone(\A_1) \to \Eigone(\A_2)$ given by
$$\gamma(\phi) = \frac{T \cdot (\phi \circ \theta^{-1}) \cdot S}
{\norm{T \cdot (\phi \circ \theta^{-1}) \cdot S}}$$ is an algebraic
isomorphism of coordinate systems such that for every
$\rho\in\hat{\D}_1$, $\gamma(\rho)=\rho\circ\alpha^{-1}$.
\end{thm}

\begin{proof}
The fact that $\gamma$ is a bijection between $\Eigone(\A_1)$ and
$\Eigone(\A_2)$ such that for every $\phi\in\Eigone(\A_1)$,
$r(\gamma(\phi))=r(\phi)\circ\alpha^{-1}$ and
$s(\gamma(\phi))=s(\phi)\circ\alpha^{-1}$ 
follows immediately from Proposition~\ref{simpleeigenmap}, so we need
only show $\gamma$ is multiplicative on composable elements.

Given $\phi\in\Eigone(\A_1)$ with minimal partial isometry  $v_\phi$ (see
Definition~\ref{asparisom}), we first identify the minimal
partial isometry for $\gamma(\phi)$.  In fact, we claim that
\begin{equation}\label{v}
v_{\gamma(\phi)}=\frac{T\theta(v_\phi)S}{\norm{T\theta(v_\phi)S}}.
\end{equation}
 To see this, first observe that  part (4) of
 Proposition~\ref{STFacts} also holds for $x\in\ddual{\A_2}.$  Thus, 
$$(T\cdot (\phi\circ\theta^{-1})\cdot S)
(T\theta(v_\phi)S)=\phi(\theta^{-1}(ST\theta(v_\phi)ST))=
\phi(v_\phi)=1.$$  Therefore,
$$\gamma(\phi)(T\theta(v_\phi)S)>0.$$  Moreover, if $q=p_{r(\phi)}$ and
$p=p_{s(\phi)}$, then by part (3) of Proposition~\ref{STFacts},
$$T\theta(v_\phi)S=T\theta(q)\theta(v_\phi)\theta(p)S=
\alpha(p)\theta(v_\phi)\alpha(q).$$ Hence
$\ds\frac{T\theta(v_\phi)S}{\norm{T\theta(v_\phi)S}}$ is a minimal
partial isometry in $\ddual{\A_2}$ on which $\gamma(\phi)$ takes a
positive value.  Thus, equation~\eqref{v} holds by
Remark~\ref{minpifacts}.

Now suppose $\phi_1,\phi_2\in\Eigone(\A_1)$ are such that
$\phi_1\phi_2$ is defined.  Notice that the minimal partial isometry
for the product $\phi_1\phi_2$ is the product of the minimal partial
isometries for $\phi_1$ and $\phi_2$, that is,
$v_{\phi_1\phi_2}=v_{\phi_1}v_{\phi_2}$.  To show that
$\gamma(\phi_1\phi_2)=\gamma(\phi_1)\gamma(\phi_2)$, it suffices to
show that
\begin{equation}\label{prodgam}
v_{\gamma(\phi_1\phi_2)}=v_{\gamma(\phi_1)}v_{\gamma(\phi_2)}.
\end{equation} 

To do this, we first show that  for all $\rho\in\hat{\D}_1$, we have
\begin{equation}\label{thetacom}
\theta(p_\rho)\alpha(p_\rho)\theta(p_\rho)=\theta(p_\rho).
\end{equation}  
Indeed, by Proposition~\ref{combineeigen}, 
\begin{align*}
p_\rho\theta^{-1}(\alpha(p_\rho))p_\rho
&=\rho(\theta^{-1}(\alpha(p_\rho))) p_\rho \\
&=\rho((E_1\circ\theta^{-1})(\alpha(p_\rho))p_\rho\\
&=\rho(\alpha^{-1}(\alpha(p_\rho)))p_\rho= 
\rho(p_\rho)p_\rho=p_\rho.
\end{align*}
Applying $\theta$ to the ends of this equality yields~\eqref{thetacom}.

Noting that $p_{s(\phi_1)}=p_{r(\phi_2)}$, we have
\begin{align*}
&(T\theta(v_{\phi_1})S)(T\theta(v_{\phi_2})S)\\
&= [\alpha(p_{r(\phi_1)})\theta(v_{\phi_1}) \alpha(p_{s(\phi_1)}) ]
[\alpha(p_{r(\phi_2)})\theta(v_{\phi_2})
  \alpha(p_{s(\phi_2)})]\\
&=
[\alpha(p_{r(\phi_1)})\theta(v_{\phi_1}) \theta(p_{s(\phi_1)}) 
\alpha(p_{s(\phi_1)})]
[\alpha(p_{r(\phi_2)})\theta(p_{r(\phi_2)})\theta(v_{\phi_2})
  \alpha(p_{s(\phi_2)})]\\
&=
\alpha(p_{r(\phi_1)})\theta(v_{\phi_1}) [\theta(p_{s(\phi_1)})
  \alpha(p_{s(\phi_1)})
\alpha(p_{r(\phi_2)})\theta(p_{r(\phi_2)})]\theta(v_{\phi_2}) 
 \alpha(p_{s(\phi_2)})\\
&=
\alpha(p_{r(\phi_1)})\theta(v_{\phi_1}) 
 \theta(v_{\phi_2})\alpha(p_{s(\phi_2)})\\
&=
\alpha(p_{r(\phi_1)})\theta(v_{\phi_1\phi_2}) \alpha(p_{s(\phi_2)})\\
&=
T\theta(v_{\phi_1\phi_2})S.
\end{align*}
This relation, together with equation~\eqref{v}, shows that
equation~\eqref{prodgam} holds, and the proof is complete.
\end{proof}

The continuity of the map $\gamma$ appearing in Theorem~\ref{algisom}
is a particularly vexing issue; in general, we do not know whether it
is continuous.  Theorem~\ref{algisom} does imply the restriction of
$\gamma$  to $\hat{\D}$ is continuous, a fact we will use in
Example~\ref{refined}.

In the following corollary, we show that in some circumstances,
$\gamma$ is ``nearly continuous'', in the sense that it is possible to
alter $\gamma$ by multiplying by an appropriate $\bbT$-valued cocycle
to obtain a continuous isomorphism of coordinate systems.  The key
hypothesis, that $\alpha$ from Definition~\ref{alphaDefn} extends to a
$*$-isomorphism of \cstar-envelopes, is in part motivated by
Theorem~\ref{easyMQS}, which shows that an \textit{isometric}
isomorphism $\theta$ between triangular algebras is $*$-extendible to
their \cstar-envelopes, and in particular, $\alpha=\theta|_{\D_1}$ is
$*$-extendible to the envelopes.

\begin{corollary}\label{envelopeiso}  
If there is a $*$-isomorphism $\pi:\C_1\rightarrow \C_2$ so that
$\pi|_{\D_1}=\alpha$, where $\alpha=E_2\circ\theta|_{\D_1}$, then the
map $\delta : \Eigone(\A_1) \to \Eigone(\A_2)$, defined by $\phi
\mapsto \phi\circ\pi^{-1}$, is an isomorphism of the coordinate
systems.  Moreover, there exists a cocycle $c:\Eigone(\A_1)\rightarrow
\bbT$ such that for every $\phi\in\Eigone(\A_1)$, 
$$\delta(\phi)=c(\phi)\gamma(\phi).$$
\end{corollary}

\begin{proof}
Clearly, $\phi\mapsto \phi\circ\pi^{-1}$ is a bicontinuous isomorphism
of $\Eigone(\C_1)$ onto $\Eigone(\C_2)$.  We must show that this map
sends $\Eigone(\A_1)$ into $\Eigone(\A_2)$.  
Fix $\phi \in \Eigone(\A_1)$ and let $\gamma$ be the map from 
Theorem~\ref{algisom}.  
By Proposition~\ref{simpleeigenmap},
\[ 
r(\phi\circ\pi^{-1})=r(\phi)\circ \pi^{-1}
	=r(\phi)\circ\alpha^{-1}=r(\gamma(\phi)),
\]
and similarly, $s(\phi\circ\pi^{-1}) =s(\gamma(\phi))$.  
By Corollary~\ref{multiples}, there is $c(\phi)\in \bbT$ so that 
$\phi\circ\pi^{-1}=c(\phi) \gamma(\phi)$, and so 
$\phi \circ \pi^{-1} \in \Eigone(\A_2)$.
Since both $\gamma$ and $\delta$ are multiplicative on composable
elements of $\Eigone(\A_1)$, so is $c$, whence $c$ is a cocycle.
\end{proof}

We now introduce a new class of algebras for
which $\gamma$ is continuous for bounded isomorphisms between
triangular algebras in the class.  This class includes those algebras
$\A\subseteq (\C,\D)$ where $\C$ admits a cover by monotone $G$-sets
with respect to $\A$ (see~\cite[p.~57]{MR90m:46098}).  In the context
of limit algebras, this class includes limit algebras generated by
their order-preserving normalizers
(see~\cite{MR96k:46099,MR2000d:47103}).

As the definition of the class does not require our standing
assumptions for the section, we relax them momentarily.

\begin{defn}\label{alp}
Let $(\C,\D)$ be a \cstar-diagonal and  $\A\subseteq (\C,\D)$ be a
subalgebra (not necessarily triangular).
We say a normalizer $v \in \A$ is \textit{algebra-preserving} if 
either $v \A v^* \subseteq \A$ or $v^* \A v \subseteq \A$.
\end{defn}

This is related to the notion of order-preserving normalizers, which
can be described as those $v \in \N_{\D}(\A)$ satisfying both 
$v \A v^* \subseteq \A$ and $v^* \A v \subseteq \A$.

As a trivial example of an algebra-preserving normalizer that is not
order preserving, let $\C$ be $M_4(\bbC)$, $\D$ the diagonal matrices,
and $\A$ the span of all upper-triangular matrix units except $e_{1,2}$.
Then $v = e_{1,3}+e_{2,4}$ normalizes $\D$ but is not order preserving,
since $v e_{3,4} v^* = e_{1,2} \notin \A$.
However, $v$ is algebra preserving, since $v^* \A v = v^* \D v \subset \A$.

\begin{lemma}\label{radicalzero}
Let $(\C,\D)$ be a \cstardiag\ and let $\A\subseteq (\C,\D)$ be
triangular.  If $\phi \in \Eigone(\A)$ and $v \in \A$ is an
algebra-preserving normalizer with $\phi(v) \ne 0$, then for  $x,y
\in \ddual{\A} \cap \ker \ddual{E} $,
\[ \phi(xv)=\phi(vy)=\phi(xvy)=0. \]
\end{lemma}

\begin{proof}
By Remark~\ref{rangesourceform} and the fact that $\ker\ddual{E}\cap
\ddual{\A}$ is a $\D$-bimodule,
\[ \phi(xv) = \frac{ r(\phi)(xvv^*) }{[r(\phi)(vv^*)]^{1/2}} = 0,
    \qquad
   \phi(vy) = \frac{ s(\phi)(v^*vy) }{[s(\phi)(v^*v)]^{1/2}} = 0,
\]
as $xvv^*, v^*vy \in \ker\ddual{E}$.

For the last equality, assume first that $v \A v^* \subseteq \A$.
Then $v\ddual{\A} v^*\subseteq \ddual{\A}$ and, as $\ker
\ddual{E}\cap\ddual{\A}$ is an
ideal in $\ddual{\A}$ (Theorem~\ref{ExpectHomomorph}),
\[ \phi(xvy) 
    = \frac{ r(\phi)(x(vyv^*)) }{[r(\phi)(vv^*)]^{1/2}}=
\frac{ r(\phi)(\ddual{E}(x(vyv^*))) }{[r(\phi)(vv^*)]^{1/2}} =0.
\]

If $v^* \A v \subseteq \A$, then, similarly,
$\phi(xvy) = s(\phi)((v^*xv)y)/[s(\phi)(v^*v)]^{1/2}$
shows that $\phi(xvy)=0$.
\end{proof}

We now reimpose the Standing Assumptions for Section~\ref{isocoordsys};
they remain in force through the remainder of the section.

\begin{theorem}\label{algpresgamma}
If $\phi \in \Eigone(\A_1)$ and $v \in \A_1$ is an algebra-preserving 
normalizer such that $\phi(v) \ne 0$,
then $\gamma^{-1}$ is continuous at $\gamma(\phi)$.

In particular, if $\A_1$ and $\A_2$ are the closed span of their
algebra-preserving normalizers, then $\gamma$ is a homeomorphism.
\end{theorem}

\begin{proof}
Letting $\rho=s(\phi)$, $\phi=\lambda [v,\rho]$ for some $\lambda \in \bbT$,
by Theorem~\ref{EigTwist}.
Without loss of generality, we may replace $v$ by $\lambda v$.

Fix $\psi \in \Eigone(\A_1)$.
By Proposition~\ref{STFacts}, $\theta^{-1}(S)=I+X$, $\theta^{-1}(T)=I+Y$ 
where $X,Y \in \ker \ddual{E_1}|_{\ddual{\A_1}}$.  
Thus,
\[ 
\psi( \theta^{-1}(S \theta(v) T)) = \psi( v + Xv + vY + XvY) = \psi( v ),
\]
by Lemma~\ref{radicalzero}.
Putting $n(\psi)=\norm{T\cdot (\psi\circ\theta^{-1}) \cdot S}^{-1}$,
we can conclude that, for all $\psi \in \Eigone(\A_1)$,
\begin{equation}\label{nbd}
\gamma(\psi)(\theta(v))=n(\psi)\psi(v).
\end{equation}

Let $(\phi_\lambda)$ be a net in $\Eigone(\A_1)$ such that
$\gamma(\phi_\lambda)\rightarrow\gamma(\phi)$
and let $\rho_\lambda=s(\phi_\lambda)$.
Then $s(\gamma(\phi_\lambda))\rightarrow s(\gamma(\phi))$,
and so, by Proposition~\ref{simpleeigenmap}, $\rho_\lambda\rightarrow \rho$.
Since $\gamma(\phi_\lambda)(\theta(v))\rightarrow
\gamma(\phi)(\theta(v))=n(\phi)\phi(v)>0$, by \eqref{nbd}, we may assume
that $\phi_\lambda(v)\neq 0$ for all $\lambda$.  
Thus, there exist $t_\lambda\in\bbT$ such that 
$\phi_\lambda = t_\lambda [v,\rho_\lambda]$.  
Using~\eqref{nbd} and the convergence of $\gamma(\phi_\lambda)$,
\begin{align*}
n(\phi_\lambda)t_\lambda [v,\rho_\lambda](v)
   &=\gamma(\phi_\lambda)(\theta(v))\\
   &\rightarrow \gamma(\phi)(\theta(v))
   =n(\phi)\phi(v)= n(\phi)[v,\rho](v).
\end{align*}
Since $\rho_\lambda\rightarrow\rho$, $[v,\rho_\lambda](v)\rightarrow
[v,\rho](v)\neq 0$, and hence $t_\lambda n(\phi_\lambda)\rightarrow n(\phi)$. 
Taking absolute values shows that $n(\phi_\lambda)\rightarrow n(\phi)$, 
and hence $t_\lambda \rightarrow 1$.  
Therefore,
$$\phi_\lambda=t_\lambda[v,\rho_\lambda]\rightarrow [v,\rho]=\phi,$$
as desired.
\end{proof}

\begin{cor}
Suppose, in addition, that each
$\A_i$ is Dirichlet and is the norm closure of the span of 
its algebra-preserving normalizers.
Then $\A_1$ and $\A_2$ are boundedly isomorphic if and only if
they are isometrically isomorphic.
\end{cor}

\begin{proof}
If $\A_1$ and $\A_2$ are boundedly isomorphic,
Theorem~\ref{algpresgamma} shows $\Eigone(\A_1)$ and $\Eigone(\A_2)$
are isomorphic coordinate systems. As the $\A_i$ are Dirichlet and
regular by hypothesis,  the result follows from an application of 
Theorem~\ref{trialgcor}. 
\end{proof}

This corollary extends~\cite[Theorem~2.3]{MR2001k:47103}, which proves
the corresponding result for strongly maximal TAF algebras generated
by their order-preserving normalizers.  The cited theorem does
somewhat more, as it shows that algebraic isomorphism implies
isometric isomorphism.

In light of Theorem~\ref{algisom}, it is natural to ask what
implications can be drawn from the existence of the algebraic
isomorphism of coordinate systems.  It has been known for more than a
decade that algebraic isomorphism of coordinate systems does not imply
isometric isomorphism~\cite[Remark on page 120]{MR94k:47068}.  It is
easily shown that the
algebras in their example fail to be boundedly
isomorphic and the algebraic isomorphism of
coordinates they exhibit is continuous on the diagonals.

We give a somewhat different example, where the algebras are anti-isomorphic,
have no minimal projections, and, most importantly, continuity on the diagonal
can be exploited to show that the algebras are not boundedly isomorphic.

\begin{remark}{Example}\label{refined}
Let $M_k$, $D_k$, $T_k$ be the algebra of $2^k \times 2^k$ matrices and
the subalgebras in $M_k$ of diagonal and upper-triangular matrices,
respectively.  Where necessary, we will equip an $M_k$ with a matrix
unit system $\{ e_{i,j}\}$.  Let $A_k$ (resp., $B_k$) be the
permutation unitary in $M_k$ that interchanges the first (resp., last)
two entries of a vector in $\bbC^{2^k}$.  We consider three embeddings
from $M_k$ to $M_{k+1}$.  First, we have $\pi_k$ that sends a matrix
$[a_{ij}]$ to the block matrix $[a_{ij}I_2]$ where $I_2$ is the $2
\times 2$ identity matrix.  Let $\alpha_k = \Ad A_{k+1} \circ \pi_k$
and
$\beta_{k} = \Ad B_{k+1} \circ \pi_k$.  Then $\M_a =
\indlim(M_k,\alpha_k)$ and $\M_b = \indlim(M_k,\beta_k)$ are both the
$2^\infty$ UHF \cstar-algebra.  Since $\alpha_k|_{D_k} =
\beta_k|_{D_k} = \pi_k|_{D_k}$, we have
$\indlim(D_k,\alpha_k|_{D_k})=\indlim(D_k,\alpha_k|_{D_k})$, which we
denote $D$.

The operator algebras $T_a = \indlim(T_k,\alpha_k|_{T_k})$ and
$T_b = \indlim(T_k,\beta_k|_{T_k})$ are anti-isomorphic.
Indeed, if $\phi_k : T_k \to T_k$ sends $e_{i,j}$ to 
$e_{2^k+1-j,2^k+1-i}$,
then $\beta_{k+1} \circ \phi_k = \phi_{k+1} \circ \alpha_k$, so the limit
of the $\phi_k$ defines an anti-isomorphism between $T_a$ and $T_b$.

Suppose that $T_a$ and $T_b$ were boundedly isomorphic.  By
Proposition~\ref{alphaMap}, we would have a $*$-isomorphism $\alpha :
D \to D$ and by Theorem~\ref{algisom}, there would be an algebraic
isomorphism $\gamma : \Eigone(T_a) \to \Eigone(T_b)$.  This induces
$\delta : R(T_a) \to R(T_b)$, an (algebraic) isomorphism of spectral
relations, namely $\delta(|\phi|)=|\gamma(\phi)|$.  Although $\delta$
need not be continuous, we know that on the diagonals of $R(T_a)$ and
$R(T_b)$, $\delta$ can be identified with $\hat{\alpha}: \hat{D} \to
\hat{D}$, the map induced by $\alpha$, and so is continuous.
Moreover, by the range and source condition in Theorem~\ref{algisom},
$(\rho,\sigma) \in R(T_a)$ if and only if
$(\hat{\alpha}(\rho),\hat{\alpha}(\sigma)) \in R(T_b)$.

Let $f$ (resp., $l$) be the element of $\hat{D}$ that equals $1$ on
the $e_{1,1}$ matrix unit (resp., $e_{2^k,2^k}$ matrix unit) in each
$D_k$.  Now $(f,l)$ is in both $R(T_a)$ and $R(T_b)$.  Let $\O_a = \{
\rho \in \hat{\D} : (f,\rho) \in R(T_a)$ and define $\O_b$ similarly.
The existence of $\delta$ implies that $\hat{\alpha}$ maps $\O_a$ onto
$\O_b$.  We claim that this is impossible for a continuous $\alpha$.
The essence of the following argument is that every element of
$\O_a\backslash\{f\}$ has a neighborhood $N$ where it is maximal in $N
\cap \O_a$, while no element of $\O_b\backslash\{l\}$ has such a
neighborhood.

The basic neighborhoods for $f$ are given by elements of $\hat{D}$ that 
are
nonzero on a $e_{1,1}$ matrix unit in some $D_k$.
If we consider some $\rho \in \O_a\backslash\{f\}$, then there is some $k$
and some $j \in \{2,\ldots,2^k\}$ so that $\rho = e_{j,1} \cdot f \cdot 
e_{1,j}$.
where $e_{1,j}$ is a matrix unit in $M_k$.
In particular, $N := \{ \sigma \in \hat{D} : \sigma(e_{j,j}) \ne 0 \}$ is 
a basic
neighborhood of $\rho$.
Every element of $N \cap \O_a$ is smaller than $\rho$ in the ordering 
induced by
$R(T_a)$, since $N$ is given by conjugating a basic neighborhood of $f$
by $e_{1,j}$ and conjugation by $e_{1,j}$ reverses the diagonal ordering 
for
pairs $(f,\psi)$, $\psi \in \O_a\backslash\{f\}$.
This last fact follows from considering the image of $e_{1,j}$ in $M_l$, 
$l > k$.

On the other hand, if $\sigma \in \O_b\backslash\{l\}$, then every
neighborhood of $\sigma$ contains elements $\tau \in \O_b$ with
$(\sigma,\tau) \in R(T_b)$.  This follows from repeating the argument
of the previous paragraph, observing that conjugation by $e_{1,j}$
preserves the $R(T_b)$-ordering for pairs $(f,\psi)$.

Pick some $\rho \in \O_a\backslash\{f,l\}$ and a neighborhood, $N$, of
$\rho$ with all elements of $N \cap \O_a$ less than $\rho$ in the
diagonal order.  Now $\hat{\alpha}$ maps $\rho$ to an element of
$\O_b\backslash\{f,l\}$, call it $\sigma$.  By the previous paragraph,
every neighborhood of $\sigma$, including $\hat{\alpha}(N)$, contains
points of $\O_b$ greater than $\sigma$ in the $R(T_b)$-ordering.  This
contradicts $\hat{\alpha}$ mapping $R(T_a)$ onto $R(T_b)$ and so $T_a$
and $T_b$ are not boundedly isomorphic.  In fact, by the automatic
continuity result of~\cite{MR2001k:47103}, there is not even an
algebraic isomorphism between $T_a$ and $T_b$.

We show there is an algebraic isomorphism from $R(T_a)$ to $R(T_b)$,
and hence, using the continuous section, between $\Eigone(T_a)$ and
$\Eigone(T_b)$.  It suffices to construct a map $h:\hat{D}\rightarrow
\hat{D}$ so that $(\sigma,\tau)\in R(T_a)$ if and only if
$(h(\sigma),h(\tau))\in R(T_b)$.

For $u\in\{a,b\},$ define
$$X_u:=\{(\sigma,\tau)\in R(T_u): \sigma=f\}\cup\{(\sigma,\tau)\in
R(T_u): \tau=l\}.$$
Before defining $h$, we observe that $R(T_a)\backslash X_a
= R(T_b)\backslash X_b$.  To see this, first let $\tilde{T}_k:=
(e_1^{}e_1^*)^\perp T_k (e_{2^k}^{}e_{2^k}^*)^\perp$.  If
$(\sigma,\tau)\in R(T_a)\backslash X_a$, there is some $p\in\bbN$
and some matrix unit $e\in \tilde{T}_p$ such that $\sigma=e \cdot \tau
\cdot e^*$.  For $k\geq p$, $\alpha_k$ and $\beta_k$ agree on the
image of $e$ in $\tilde{T}_k.$ Thus, the image of $e$ in $T_a$ and the
image of $e$ in $T_b$ induce the same partial homeomorphism of
$\hat{D}$.  In particular, $(\sigma,\tau)\in R(T_b)\backslash X_b$.
The reverse inclusion is similar.

Thus, each of $\O_a\backslash\{f,l\}$ and $\O_b\backslash\{f,l\}$ is
ordered the same way by both $R(T_a)$ and $R(T_b)$.  In fact, we can
show that each set is ordered like $\bbQ$---for example, given
$(\sigma,\tau) \in \O_a\backslash\{f,l\}$ with $\sigma \ne \tau$,
find, as above, an off-diagonal matrix unit $e$ and use the images of
$ee^*$ and $e^*e$ in a later matrix algebra to show there is $\eta \in
\hat{D}$ with $(\sigma,\eta),(\eta,\tau) \in R(T_a)$ and $\eta$
different from $\sigma$ and $\tau$.

Define $h$ to be the identity map everywhere except
$\O_a\backslash\{f,l\}$ and $\O_b\backslash\{f,l\}$.  Since these two
sets have the same order type, there is a bijection $g :
\O_a\backslash\{f,l\} \to \O_b\backslash\{f,l\}$ so that
$(\sigma,\tau) \in \O_a\backslash\{f,l\}$ if and only if
$(g(\sigma),g(\tau)) \in \O_b\backslash\{f,l\}$.  As
$\O_a\setminus\{f,l\}$ and $\O_b\setminus\{f,l\}$ are disjoint sets,
we define $h$ on $\O_a\backslash\{f,l\}$ to be $g$ and on
$\O_b\backslash\{f,l\}$ to be $g^{-1}$.  Using the observation above
and the definitions of $\O_a$ and $\O_b$, it is straightforward to
show that $h$ has the required property.
\end{remark}

\section{Structure of General Isomorphisms}
\label{S:isostruct}

We now turn to the structure of isomorphisms of triangular algebras.
First, we build on the representation results from
Section~\ref{repns}, obtaining Theorem~\ref{CSLextend}, which extends
such an isomorphism to an isomorphism of CSL algebras.  After several
results about isomorphisms of CSL algebras, we obtain an analogue of a
factorization result of Arveson-Josephson, Theorem~\ref{isostructure}.
A main result of the paper is Theorem~\ref{CB}, which shows that such
an isomorphism is completely bounded.  Finally, we extend a result of
Muhly-Qiu-Solel, Theorem~\ref{easyMQS}, showing that an isometric
isomorphism extends to $*$-isomorphism of the \cstar-diagonals.

Our standing assumptions are the same as those of the previous section.

We start with a technical lemma.

\begin{lemma}\label{techie}  Suppose that $C^*(\A_i)=\C_i$. 
Given $\rho_2\in\hat{\D}_2$, let $\rho_1=\rho_2\circ\alpha\in\hat{\D}_1$, 
and let $(\pi_i,\H_i)$ be the $($compatible$)$ GNS representations of 
$(\C_i,\D_i)$ on $\H_i$ corresponding to $\rho_i$. 

If, for $i=1,2$, $P_i$ is the support projection for $\pi_i$,  
then $\theta(P_1)=P_2$.
\end{lemma}

\begin{proof}
By Proposition~\ref{specialrep}, $P_i \in \ddual{\D_i}$ and
$P_i=\sum_{q\in\O_i}q$, where 
$$\O_i:= \{q\in\ddual{\D_i}: q=p_\sigma \text{ for some } 
    \sigma\in\hat{\D}_i \text{ such that } (\rho_i,\sigma) \in R(\C_i) \}.
$$ 

Since $\A_i$ generates $\C_i$, Theorem~\ref{envelope} shows that
the equivalence relation generated by $R(\A_i)$ is $R(\C_i)$. 
By Theorem~\ref{algisom}, we have $(\sigma,\rho_2)\in R(\A_2)$ 
if and only if $(\sigma\circ\alpha,\rho_1)\in R(\A_1)$. 
Hence $(\sigma,\rho_2) \in R(\C_2)$ if and only if
$(\sigma\circ\alpha,\rho_1) \in R(\C_1)$.
Theorem~\ref{algisom} also shows that for $\sigma_1\in\hat{\D}_1$,
$\alpha(p_{\sigma_1})=p_{\sigma_1\circ\alpha^{-1}}.$  Therefore,
$\alpha(\O_1)=\O_2$ and we obtain,  
\begin{equation}\label{Pfine}
P_2=\sum_{p\in\O_1}\alpha(p).
\end{equation}

Fix $p\in\O_1$.  Since $p$ is a minimal projection in $\ddual{\C_1}$,
it is a minimal idempotent in $\ddual{\A_1},$ so $\theta(p)$ is a
minimal idempotent in $\ddual{\A_2}$.  As $P_2$ is a central
projection in $\ddual{\C_2}$, $\theta(p)P_2$ is an idempotent in
$\ddual{\A_2}$ and so $\theta(p)P_2$ is either $\theta(p)$ or $0$.  As
done earlier, we again use $\tilde{\pi}_i$ for the unique extension of
$\pi_i$ to $\ddual{\C_i}$.  Since
\[ 
\tilde{\pi}_2(E_2(\theta(p)P_2)) = \tilde{\pi}_2(\alpha(p)) \neq 0,
\]
we must have $\theta(p)P_2=\theta(p)$.  

The $\sigma(\ddual{\A_1},\dual{\A_1})$-$\sigma(\ddual{\A_2},\dual{\A_2})$
continuity of $\theta$ yields
\begin{equation}\label{thetaorbitA}
\theta(P_1)=\theta\left(\sum_{p\in\O_1} p\right) =\sum_{p\in\O_1}\theta(p)
	=\sum_{p\in\O_1}P_2\theta(p)=
\theta(P_1)P_2.
\end{equation}
Similar considerations show that for every $p\in\O_1$,
$\theta^{-1}(\alpha(p))P_1 =\theta^{-1}(\alpha(p))$ and
\begin{equation}\label{thetaorbitB}
\theta^{-1}(P_2)=
\sum_{p\in\O_1}P_1\theta^{-1}(\alpha(p))=P_1\theta^{-1}(P_2).
\end{equation}

Applying $\theta$ to \eqref{thetaorbitB} and 
using~\eqref{thetaorbitA} yields
$\theta(P_1)= P_2$.  
\end{proof}

The support projection of a direct sum of inequivalent representations
of a \cstaralg\ 
is the sum of the support projections of the individual
representations.  Thus, Lemma~\ref{techie} and
Theorem~\ref{faithfulcompatable} combine to produce the following
result, which  connects our context to the theory of CSL algebras.

When the \cstar-envelope of $\A_i$ is $\C_i$, Theorem~\ref{envelope}
shows that $R(\C_1)$ and $R(\C_2)$ are isomorphic as topological
equivalence relations.  
Thus, the assumption on $\X_2$ below implies that $\X_2$ also has exactly 
one element from each $R(\C_2)$-equivalence class.

\begin{theorem}\label{CSLextend}  
Suppose that $C^*(\A_i)=\C_i$.  Let $\X_2\subseteq \hat{\D}_2$ contain
exactly one element from each $R(\C_2)$ equivalence class and let
$\X_1=\{\rho\circ\alpha: \rho\in\X_2\}$.  Let
$\pi_i=\bigoplus_{\rho\in\X_i} \pi_\rho$ be the faithful, compatible
representations of $(\C_i,\D_i)$ as constructed in
Theorem~\ref{faithfulcompatable}.

If $\theta':\pi_1(\A_1)\rightarrow \pi_2(\A_2)$ is the map given by
$\theta'(\pi_1(a))=\pi_2(\theta(a))$, then
$\theta'$ extends uniquely to a $($bounded\,$)$ isomorphism $\overline{\theta} : 
\tilde{\pi}_1(\ddual{\A_1}) \to \tilde{\pi}_2(\ddual{\A_2})$.
\end{theorem}

\begin{proof}
Let $P_i$ be the support projections of $\pi_i$.  
By Proposition~\ref{specialrep} and Lemma~\ref{techie},
$P_i\in\ddual{\D_i}$ and $\theta(P_1)=P_2$.  
By Theorem~\ref{faithfulcompatable},
$\ker\tilde{\pi}_i|_{\ddual{\A_i}}= P_i^\perp \ddual{\A_i}$
and so $\tilde{\pi}_i$ is faithful on $P_i\ddual{\A_i}$.
As $\tilde{\pi}_i|_{P_i\ddual{\A_i}}$ 
has image $\tilde{\pi}_i(\ddual{\A_i})$, the map
$\overline{\theta}:\tilde{\pi}_1(\ddual{\A_1})\rightarrow
\tilde{\pi}_2(\ddual{\A_2})$ given by
$\tilde{\pi}_1(a)\mapsto \tilde{\pi}_2(\theta(a))$ is well defined.

Uniqueness follows from the weak-$*$ density of $\A_i$ in $\ddual{\A_i}$. 
\end{proof}

For a representation, $\pi$, of $\C_2$, we suspect that $\tilde{\pi}(S)$ and 
$\tilde{\pi}(T)$ are inverses of 
each other whenever $\pi$ is a compatible representation.
The next two propositions offer some evidence for this.
Indeed, Theorem~\ref{Sinvertible} proves it when $\pi_2$ is the faithful
compatible 
atomic representations of Theorem~\ref{CSLextend}.

\begin{proposition}\label{leftinverse}
If  $\pi$ is a compatible representation of $\C_2$ on $\H$,
then $\tilde{\pi}(TS)=I$.
\end{proposition}

\begin{proof}
{}From Proposition~\ref{STFacts}, we have
$\tilde{E_2}(\tilde{\pi}(S))=\tilde{E_2}(\tilde{\pi}(T))=I$.
Proposition~\ref{STFacts} also shows that $\tilde{\pi}(TS)\in\pi(\alpha(\D_1))'$,
so $\tilde{\pi}(TS)\in\pi(\D_2)''$ since $\pi(\D_2)''$ is a MASA in
$\pi(\C)''$.
Since $\tilde{E}_2$ is a homomorphism on $\pi(\A_2)$ it is also a
homomorphism on $\tilde{\pi}(\ddual{\A_2})$.
Hence,
$$\tilde{\pi}(TS)=\tilde{E}_2(\tilde{\pi}(TS))=
\tilde{E}_2(\tilde{\pi}(T))\tilde{E}_2(\tilde{\pi}(S))=I.$$
\end{proof}

\begin{theorem}\label{Sinvertible}  
For $\pi_2$ as in Theorem~\ref{CSLextend}, 
$\tilde{\pi}_2(S)^{-1}=\tilde{\pi}_2(T)$.
\end{theorem}

\begin{proof}
We use the same notation as in Lemma~\ref{techie},
Theorem~\ref{CSLextend} and their proofs.
Fix  
$\rho_2 \in \hat{\D}_2$.  We claim that 
$\tilde{\pi}_{\rho_2}(S)$ is invertible and
$\tilde{\pi}_{\rho_2}(S)^{-1}=\tilde{\pi}_{\rho_2}(T)$.

Applying $\tilde{\pi}_{\rho_2}$ to $\sum_{p\in\O_1}\theta(p)=P_2$
(obtained from Lemma~\ref{techie} and the first equality
in~\eqref{thetaorbitA}) yields the important equality,
\begin{equation}\label{sumthetap}
\sum_{p\in\O_1}\tilde{\pi}_{\rho_2}(\theta(p))=I_{\H_{\rho_2}}.
\end{equation}

By Proposition~\ref{STFacts}, $S\alpha(p)=\theta(p)\alpha(p)$ for all
$p$, and using ~\eqref{sumthetap} gives
\[
\tilde{\pi}_{\rho_2}(S)=\sum_{p\in\O_1}
                \tilde{\pi}_{\rho_2}(S\alpha(p)) =\sum_{p\in\O_1}
                \tilde{\pi}_{\rho_2}(\theta(p)\alpha(p))
\]
A similar calculation with $T$ gives
\[ \tilde{\pi}_{\rho_2}(T)
=\sum_{p\in\O_1}\tilde{\pi}_{\rho_2}(\alpha(p)\theta(p)).
\]
Finally, we then have
$$\tilde{\pi}_{\rho_2}(ST)=\sum_{p\in\O_1}\tilde{\pi}_{\rho_2}(\theta(p))=I_{\H_{\rho_2}}.$$
Proposition~\ref{leftinverse} established
$\tilde{\pi}_{\rho_2}(TS)=I_{\H_{\rho_2}}$, and hence our claim holds.

As $\pi_2=\bigoplus_{\rho_2\in\X_2}\pi_{\rho_2}$, the result follows. 
\end{proof}

We need two structural results for CSL algebras.
The factorization result, Lemma~\ref{atomicfactor}, is
well known, and we only sketch its proof.

\begin{theorem}[{\cite[Theorem~2.1]{GilfeatherMooreIsCeCSLAl}}]
\hspace*{-8pt}\footnote{Gilfeather and Moore
attribute this result to Ringrose in the nest algebra case and to
Hopenwasser for CSL algebras.  However, Gilfeather and Moore show that
$\beta$ is a bounded automorphism.}
\label{CSLisoform}\hspace*{0pt}  Suppose $\L_1$ and $\L_2$ are CSLs on 
Hilbert spaces $\H_1$ and $\H_2$ and that $\pi:\alg\L_1\rightarrow \alg\L_2$ 
is an algebra isomorphism.  
Then, given a MASA $\fM\subseteq \B(\H_1)$ which is also contained in $\alg\L_1$, 
there exist an invertible 
operator $X\in\B(\H_1,\H_2)$ and an automorphism 
$\beta:\alg\L_1\rightarrow\alg\L_1$ such that, for every $T\in\alg\L_1$,
$$\pi(T)=X\beta(T)X^{-1}\quad \text{and} \quad
  \beta|_{\fM}=\text{Id}_{\fM}.$$
\end{theorem} 

\begin{lemma}\label{atomicfactor}
Suppose $\L_1$ and $\L_2$ are atomic CSLs on Hilbert spaces $\H_1$ and
$\H_2$ and that $X\alg\L_1 X^{-1} =\alg\L_2$ for an invertible
operator $X \in \B(\H_1,\H_2)$.  There exists a unitary operator
$U\in\B(\H_1,\H_2)$ and an invertible operator $A\in\alg\L_1$ such
that $A^{-1}\in\alg\L_1$ and $X=UA$.
\end{lemma}

\begin{proof}  
Regard $\L_i$ as a commuting family of projections in $\B(\H_i)$.  
Let $\bbA_i$ be the set of minimal projections in $\L_i''$.  
By hypothesis, $I=\sum_{a\in\bbA_i}a$. 
Define $\mu:\L_1\rightarrow \L_2$ by $\mu(P)= [XP]$, where $[XP]$ denotes
  the projection onto the range of $XP$.  
Then $\mu$ is a complete lattice isomorphism. 
As each minimal projection in $\bbA_1$ has the form $PQ^\perp$ for some 
$P,Q\in\L_1$, we see that $\mu$ induces a map, $\mu':\bbA_1\rightarrow\bbA_2$ 
given by $\mu'(PQ^\perp)=\mu(P)\mu(Q)^\perp$.  
This map is well defined and bijective.  
Also, for each $P\in\L_1$, 
\begin{equation} \label{conta} 
  P=\sum_{\{a\in\bbA_1: a\leq P\}}a, \qquad 
  \mu(P) =\sum_{\{a\in\bbA_1: a\leq P\}}\mu'(a).
\end{equation} 

Since $\L_i$ atomic and the isomorphism between $\alg\L_1$ and
$\alg\L_2$ is given by an invertible element $X$, $\dim
a\H_1=\dim \mu'(a)\H_2$ for every $a\in\bbA_1$.  For $a\in\bbA_1$,
let $u_a: \H_1 \to \H_2 $ be a partial isometry with $u_au_a^*=\mu'(a)$ and
$u_a^*u_a=a$.  
Put $$U=\sum_{a\in\bbA_1} u_a.$$  
Then $U$ is unitary and it follows from \eqref{conta}
that $A:=U^*X$ satisfies $A,A^{-1}\in \alg\L_1$.
\end{proof}

The next two results are structural results for bounded isomorphisms
of triangular algebras.  The first is an analog of a result of Arveson
and Josephson \cite[Theorem 4.10]{ArvesonJosephsonOpAlMePrAuII}
appropriate to our setting.  Briefly, Arveson and Josephson study a
variant of the crossed product algebra associated to a homeomorphism
of a locally compact Hausdorff space.  If the homeomorphism has
no periodic points, then results
in~\cite{ArvesonJosephsonOpAlMePrAuII} show easily that the resulting
algebra is a triangular subalgebra of a \cstar-diagonal (see
also~\cite[Section~4]{MR1376551}).  Arveson and Josephson show that a
bounded isomorphism of these algebras factors into three maps, the
first an isometric map arising from a homeomorphism of the underlying
spaces, the second an isometric map arising from a diagonal unitary,
and the third a weakly inner automorphism, i.e., one implemented by an
invertible in the ultraweak closure of a suitable representation.

The main difference in the form of the Arveson-Josephson factorization
and the factorization in Theorem~\ref{isostructure} below is that we
do not know if the approximately inner part of our factorization
carries $\A_1$ to itself, so we need to introduce an algebra $\A_3$.
We also remark that isomorphisms which fix the diagonal pointwise are
essentially cocycle automorphisms (see Definition~\ref{cocycledefn}).

\begin{theorem}\label{isostructure} Assume that $C^*(\A_i)=\C_i$.  Let
  $\pi:\C_1\rightarrow \bh$ be the faithful compatible representation
  of $\A_1$ constructed in Theorem~\ref{faithfulcompatable}, and let
  $\alg\L$ be the weak-$*$ closure $($in \bh$)$ of $\pi(\A_1)$.  
Then $\theta$ factors as 
\[ \theta= \tau\circ \Ad A \circ \beta\circ\pi|_{\A_1}, \]
where $\beta\in\Aut(\alg\L)$ with $\beta(x) = x$ for $x \in
\pi(\D_1)''$, $A\in \alg\L$ with $A$ invertible and $A^{-1}\in\alg\L$,
and, finally, if $\A_3:=(\Ad A \circ \beta)(\pi(\A_1)),$ then
$\A_3\subseteq\alg\L$ and $\tau:\A_3\rightarrow\A_2$  is
an isometric isomorphism.
\end{theorem}

\begin{proof}  
Apply Theorem~\ref{CSLextend} and Theorem~\ref{CSLisoform} to obtain
an invertible operator $X\in\B(\H_1,\H_2)$ which implements a
similarity between 
$\tilde{\pi}_1(\ddual{\A_1})$ and $\tilde{\pi}_2(\ddual{\A_2}).$  
Factor $X$ as $UA$ where $U$ is unitary and $A,A^{-1}\in\alg\L$, as
in Lemma~\ref{atomicfactor}.
Take $\tau=\Ad U|_{\A_3}$.  
The result follows from Theorem~\ref{CSLisoform}.
\end{proof}

We come now to a main result.

\begin{theorem}\label{CB} Suppose that $C^*(\A_i)=\C_i$.
If $\theta:\A_1\rightarrow \A_2$ is a bounded isomorphism,
then $\theta$ is completely bounded and 
$\norm{\theta}_{cb} =\norm{\theta}$.
\end{theorem}

\begin{proof}
Using Theorem~\ref{CSLextend} (and its notation), we have
$\overline{\theta}:\tilde{\pi}_1(\ddual{\A_1})\rightarrow
\tilde{\pi}_2(\ddual{\A_2})$.
By~Theorem~\ref{CSLisoform}, $\overline{\theta}$ factors as
$\overline{\theta}=\Ad X\circ\beta$, where $X:\H_{\pi_1}\rightarrow
\H_{\pi_2}$ is a bounded invertible operator and $\beta$ is an
automorphism of $\tilde{\pi}_1(\ddual{A_1})$ fixing
$\tilde{\pi}_1(\ddual{\D_1})$ pointwise.
By Lemma~\ref{atomicfactor}, $X=UA$ where $A$ and $A^{-1}$ both
belong to $\tilde{\pi}_1(\ddual{\A_1})$ and $U$ is a unitary operator.
Then $\Ad A\circ\beta$ is an automorphism of
$\tilde{\pi}_1(\ddual{\A_1})$ whose norm is $\norm{\theta}$.  

By~\cite[Corollary~2.5 and Theorem~2.6]{MR92h:47058},
$\norm{\Ad A\circ\beta}_{cb}=\norm{\Ad A\circ\beta}=\norm{\theta}$.  
Thus, $\norm{\overline{\theta}}_{cb} =\norm{\theta}$.   
Therefore, for $\theta'$ as in Theorem~\ref{CSLextend}, 
$\theta'$ is completely bounded. 
Since $\norm{\theta'}_{cb} \leq \norm{\overline{\theta}}_{cb}$, 
we have $\norm{\theta'}_{cb}=\norm{\theta}$.  
Noting that each $\pi_i$ is a complete isometry of $\A_i$ onto 
its respective range completes the proof.
\end{proof}

Finally, we use the universal property of \cstar-envelopes to generalize
a result of Muhly, Qiu, and Solel, \cite[Theorem~1.1]{MR95a:46080}.  
Their result includes a corresponding statement for anti-isomorphisms, which
can be deduced from the statement below by considering appropriate opposite algebras.
Our generalization does not require nuclearity of the $\C_i$ or the second
countability of the $\hat{\D}_i$, as we do not use the spectral
theorem for bimodules, \cite[Theorem~4.1]{MR94i:46075}.

This result also generalizes Corollary~\ref{trialgcor} for isometric $\theta$ 
from triangular subdiagonal algebras to general triangular subalgebras. 

\begin{theorem}\label{easyMQS} 
 For $i=1,2$, let $\A_i$ be a triangular
  subalgebra of the \cstardiag\ $(\C_i,\D_i)$ and
  assume that $\A_i$ generates $\C_i$ as a \cstar-algebra.  If
  $\theta:\A_1\rightarrow \A_2$ is an isometric isomorphism, then
there is a unique $*$-isomorphism $\pi:\C_1\rightarrow\C_2$ with 
$\pi|_{\A_1}=\theta$.
\end{theorem}

\begin{proof}
By Proposition~\ref{envelope}, we know that $\C_i$ is the
\cstar-envelope of $\A_i$.  
Since $\theta$ is completely isometric by Theorem~\ref{CB}, the universal 
property for \cstar-envelopes shows that there exist unique $*$-epimorphisms
$\pi_{12}:\C_2\rightarrow \C_1$ and $\pi_{21}:\C_1\rightarrow \C_2$ 
such that 
$$\pi_{12}\circ i_2=i_1\circ\theta^{-1}\quad \text{and}\quad \pi_{21}\circ
i_1= i_2\circ\theta.$$  
where, for $k=1,2$, $i_k$ is the inclusion mapping of $\A_k$ into $\C_k$.  
Thus, $\pi_{12}\circ \pi_{21}\circ i_1=\text{Id}|_{i_1(\A_1)}$, and hence
$\pi_{12}\circ \pi_{21}=\text{Id}|_{\C_1}$.  
Thus $\pi_{21}$ is injective, and is the required $*$-isomorphism 
of $\C_1$ onto $\C_2$.
\end{proof}

\begin{remark}{Remark}\label{intrinsicdef}
In Theorems~\ref{CSLextend},~\ref{isostructure},~\ref{CB} and~\ref{easyMQS}, 
we require that $\C_i$ is the \cstar-envelope of $\A_i$,
which is somewhat unsatisfying, as we would prefer conditions in terms
of $\A_i$ alone.
The hypothesis that $C^*(\A_i)=\C_i$ could be removed if we knew that every 
\cstaralg\ $\B\subseteq (\C,\D)$ is regular, for then $(C^*(\A_i),\D_i)$ would 
again be a \cstar-diagonal.   
\end{remark}

\section{Bounded Isomorphism to $*$-Extendible Isomorphism}\label{S:IndLim}

Roughly speaking, Theorem~\ref{easyMQS} states that an isometric
isomorphism of triangular algebras is $*$-extendible.  Clearly a
bounded, nonisometric isomorphism between triangular algebras cannot
be extended to a $*$-isomorphism of the \cstar-envelopes, but it still
may be the case that the \cstar-envelopes of the triangular algebras
are $*$-isomorphic.  

\begin{remark}{Question}\label{bdiso}
Suppose $\A_i\subseteq (\C_i,\D_i)$ are triangular algebras such that
$C^*(\A_i)=\C_i$.  If $\A_1$ and $\A_2$ are (boundedly) isomorphic,
are $\C_1$ and $\C_2$ $*$-isomorphic?
\end{remark}
 In view of Theorems~\ref{bigSteve} and~\ref{algisom}, 
one might expect an affirmative answer when there exists a continuous
section of $R(\C_i)$ into $\Eigone(\C_i)$, since $R(\A_i)$ generates
$R(\C_i)$ by Theorem~\ref{envelope}.  However, 
Theorem~\ref{algisom} only implies an \textit{algebraic} isomorphism of
$R(\A_1)$ onto $R(\A_2)$; to apply Theorem~\ref{bigSteve}, we need
to know that the isomorphism of $R(\A_1)$ onto $R(\A_2)$ is
continuous.   Establishing the continuity of the map $\gamma$ from
Theorem~\ref{algisom} would immediately do this.

In this section, we provide an affirmative answer to
Question~\ref{bdiso} for the class of triangular limit algebras.
Perhaps surprisingly, our proof of this result uses $K$-theory.  We do
not know whether the isomorphism obtained satisfies the hypotheses of
Corollary~\ref{envelopeiso}, so we cannot use that corollary to
establish the existence of a continuous mapping between the coordinate
systems or spectral relations of the triangular limit algebras.

We start with a theorem about Murray-von Neumann equivalence.
The proof uses the ideas developed in the previous section.  
Recall that for any Banach algebra $\B$, two idempotents $e,f\in\B$ are 
\textit{algebraically equivalent} if there exist $x,y\in\B$ such that 
$xy=e$ and $yx=f$.

\begin{thm}\label{bigtechie}
For $i=1,2$, suppose $(\C_i,\D_i)$ are \cstar-diagonals, $\A_i \subseteq
\C_i$ are triangular subalgebras, and $\theta : \A_1 \to \A_2$ is a bounded
isomorphism.  If $u \in \A_1$ is a partial isometry intertwiner, then
$\theta(uu^*)$ and $\theta(u^*u)$ are algebraically equivalent in
$\C_2$.
\end{thm}

To prove the theorem, we need the following well-known result.
We give a proof to be self-contained.

\begin{lemma}\label{rangeproj}
Let $\C\subseteq \bh$ be a concrete unital \cstar-algebra.  If $e$ is
an idempotent in $\C$, then the projection $P_e$ onto the range of $e$ is
given by $P_e=(ee^*+(1-e)^*(1-e))^{-1}ee^*$ and so is in $\C$.  Moreover,
if $z=I-eP_e^\perp$, then $z$ is an invertible element of $\C$ and 
 $zP_ez^{-1}=e$.

Suppose $e$ and $f$ are idempotents in $\C$ and there exists an
element $x\in\C$  so that
$xe=x=fx$ and, as a map from $e\H$ to $f\H$, $x$ is invertible. 
If  $x = v |x|$ is the polar decomposition for $x$, then  $v \in \C$.
\end{lemma}

\begin{proof}
As $ee^*$ and $(1-e)^*(1-e)$ commute and their product is
0, they are self-adjoint elements which are supported on orthogonal
subspaces of $\H$.  Since $ee^*$ is bounded below on $e^*\H$, it is
invertible on $e^*\H$.  Similarly, $(1-e)^*(1-e)$ is invertible on
$(1-e)\H$.  Since the kernel of $e^*$ is the range of $(1-e)$,
$ee^*+(1-e)^*(1-e)$ is invertible.  By decomposing $\H$ into the sum
of the ranges of $e$ and $1-e^*$ and looking at the block matrix
decomposition of $ee^*+(1-e)^*(1-e)$, we see that
$(ee^*+(1-e)^*(1-e))^{-1}ee^*$ is the projection onto the range of
$ee^*$, which is the range of $e$.  

Since $e$ and $P_e$ have the same range, $eP_e=P_e$ and $P_ee=e$.  A
calculation now shows that $z^{-1}= I+eP_e^\perp$ and $zP_e=P_e= ez.$

Turning to $x$, let $P \in \C$ be the range projection of $e$.
Then $|x|$ is invertible from $e\H$ to $e\H$, so $|x|+(1-P)$ is an invertible
operator on $\H$.
As $v^*v$ is the projection onto $e\H=|x|\H$, $v(1-P)=0$.
Thus, $v=x(|x|+(1-P))^{-1} \in \C$.
\end{proof}

\begin{proof}[Proof of Theorem~\ref{bigtechie}]
By Theorem~\ref{CSLextend}, there exists an isomorphism
$\overline{\theta}$ between the  atomic CSL
algebras $\tilde{\pi}_1(\ddual{\A_1})$ and
$\tilde{\pi}_2(\ddual{\A_2})$.  Invoking Theorem~\ref{CSLisoform}, we
can factor $\overline{\theta}$ as $\Ad X \circ \beta$ where $\beta$ is
an automorphism of $\tilde{\pi}_1(\ddual{\A_1})$ that fixes
$\pi(\D_1)''$ pointwise and $X$ is invertible.

For ease of notation, identify $\C_i$ with its image $\pi_i(\C_i)$; in
particular, we write $u$ instead of $\pi_1(u)$, etc..

Since $\beta$ fixes $\D_1''$ and $u$ is a partial isometry
intertwiner, for every $d\in\D_1$, we have $\beta(u)d= udu^*\beta(u)$;
and hence $u^*\beta(u)\in\D_1'=\D_1''$.  
Let $r=u^*\beta^{-1}(u) u^*.$  We claim that 
\begin{equation}\label{mvnform}
r\beta(u)=u^*u\quad \text{and} \quad \beta(u)r=uu^*.
\end{equation} 
Indeed, 
$$r\beta(u) = u^*\beta^{-1}(u)u^*\beta(u)= u^*\beta^{-1}(uu^*\beta(u))
= u^*\beta^{-1}(\beta(u))=u^*u.
$$  
The other equality is similar.  We have 
$Xu^*uX^{-1}=X\beta(u^*u)X^{-1}=\theta(u^*u)$, and similarly,
$Xuu^*X^{-1}=\theta(uu^*)$.  Thus  \eqref{mvnform} yields,
$$(XrX^{-1})\theta(u)=\theta(u^*u)\quad \text{and} \quad
\theta(u)(XrX^{-1})= \theta(uu^*),$$
so that $\theta(u)$ is invertible as an operator from
the range of $\theta(u^*u)$ onto the range of $\theta(uu^*)$.  

Invoking the second part of Proposition~\ref{rangeproj},
$\theta(u)=v|\theta(u)|$ with $v \in \C_2$.  Thus, the range
projections of $\theta(u^*u)$ and $\theta(uu^*)$ are algebraically 
equivalent, and hence $\theta(u^*u)$ and $\theta(uu^*)$ are
algebraically 
equivalent in $\C_2$.
\end{proof}

\begin{remark}{Remark}  Uniqueness of inverses shows that actually
  $XrX^{-1}\in\C_2$. 
\end{remark}

\begin{remark}{Definition}  \label{limitalgdefn}
For $n\in\bbN$, let $(\C_n,\D_n)$ be a \cstardiag, where $\C_n$ is a
unital 
finite dimensional \cstaralg, and suppose each $\alpha_n:\C_n\rightarrow
\C_{n+1}$ is  a regular $*$-monomorphism.  Theorem~\ref{indlimit}
shows that the inductive
limit, $(\indlimit \C_n,\indlimit\D_n)$ is a \cstardiag, which we call an
\textit{AF-\cstardiag}.  
The MASA $\indlimit \D$ is often called a \textit{canonical MASA}.
A norm-closed  subalgebra $\A\subseteq
(\indlimit\C_n,\indlimit \D_n)$ is called a \textit{limit algebra.}
\end{remark}

We reprise some of the results on limit algebras we require.
For a more detailed exposition, see~\cite{MR94g:46001} or the
introduction to~\cite{MR2002k:47148}.  

Let $(\C,\D)=(\indlimit \C_n,\indlimit\D_n)$ be an AF \cstardiag, and
 let $\A\subseteq (\C,\D)$ be a limit algebra.  For $v \in
 N_{\D}(\C)$, there is some $i$ so that we can write $v=dw$ where $d
 \in \D$ and $w \in N_{\D_i}(\C_i)$.  We should also point out that,
 if $v \in \bigcup_{k=1}^\infty N_{\D_k}(\C_k)$, the sets $S(v)$ and
 $S(v^*)$ are closed and open, so by Propositions~\ref{normalTOinter}
 and~\ref{interTOnormal}, $v$ is also an intertwiner.  Also,
 $\A_k := \C_k \cap \A$ is a finite-dimensional CSL
 algebra in $\C_k$, and  $\A$ is the closed union of the $\A_k$.

The \cstar-subalgebra $\B$ of $\C$ generated by $\A$ is again an
AF-algebra containing $\D$, and $(\B,\D)$ is again an AF-\cstardiag.
By Proposition~\ref{envelope}, $\B$ is the \cstar-envelope of $\A$.
Thus by replacing $\C$ with $\B$ if necessary, we may, and shall,
always assume that $\A$ generates $\C$ as a \cstaralg.

The spectrum, or fundamental relation, of $\A$, was first defined
in~\cite{MR92a:47053}, as pairs $(\sigma,\rho) \in \hat{\D} \times
\hat{\D}$ for which there is a partial isometry normalizer $v\in\A$
with $\sigma = v \cdot \rho \cdot v^*$.  In our notation, this is
$R(\A)$.  The spectrum can also be described by picking systems of
matrix units for each $\C_n$ so that matrix units in $\C_n$ are sums
of matrix units in $\C_{n+1}$ and then considering those elements of
$\dual{\A}$ that are either $0$ or $1$ on all matrix units.  These
elements of $\dual{\A}$ are eigenfunctionals and this description
provides  a continuous section from $\R(\A)$ to
$\Eigone(\A)$.

We require a technical result on normalizing idempotents in a
 triangular subalgebra  of a \textit{finite dimensional} \cstardiag.
The method is similar to that of Proposition~\ref{STFacts}, and is in
 fact what led to the constructions of $S$ and $T$.

\begin{lemma}\label{Boolean} Suppose that $(\C,\D)$ is a \cstardiag\ with $\C$
  finite dimensional, and $\A\subseteq (\C,\D)$ is triangular.  Let
 $\fB\subseteq\A$ be a (necessarily finite) Boolean algebra of
 commuting idempotents.  Then there exists an invertible element
 $A\in\A$ such that $A\fB A^{-1}\subseteq\D$ is a Boolean algebra of
 idempotents. 
\end{lemma}

\begin{proof} 
Let $G=\{I-2e: e\in \fB\}$; then $G$ is a finite group whose elements are all
  square roots of the identity.  Clearly $G$ is in bijective
  correspondence with $\fB$.
Define
$$S=\frac{1}{|G|}\sum_{g\in G} E(g)g^{-1}.$$
A calculation shows that for any $g\in G$, $E(g)S=Sg$, and we have
$E(S)=I$.  Thus, $S=I+Y$ where $Y\in \A$ is nilpotent, and we conclude
that $S$ is invertible.  Then for every $e\in \fB$, $SeS^{-1}= E(e)$,
and we are done.
\end{proof}

\begin{corollary}\label{idemsim}
Suppose $(\C,\D)$ is an AF-\cstardiag\ and $\A\subseteq (\C,\D)$ is a
triangular subalgebra.  If $e\in\A$ is an idempotent, then there exist
$A\in\A$ such that $AeA^{-1}=E(e).$
\end{corollary}
\begin{proof}
Write $(\C,\D)=\indlimit (\C_n,\D_n)$ where $(\C_n,\D_n)$ are finite
dimensional \cstardiag s, and let $\A_n=\C_n\cap \A$, so that
$\A=\indlim \A_n$.  By \cite[Proposition~4.5.1]{BlackadarKThOpAl2ed}, there
exists $n\in\bbN$, an idempotent $f\in\A_n$ and an invertible
element $X\in\A$ so that $XeX^{-1}=f$.   Lemma~\ref{Boolean} (applied
to $\{0,e, I-e, I\}$) shows
that there exists $S\in \A_n$ so that $SfS^{-1}=E(f)$.  Thus
$(SX)e(SX)^{-1}=E(f)$.  Since $E|_{\A}$ is a homomorphism, applying $E$
to the previous equality yields $E(e)=E(f)$, and the proof is complete.
\end{proof}

\begin{remark}{Remark}\label{ktheormk} 
Let $\iota:\D\rightarrow \A$ be the inclusion map
  of $\D$ into the triangular limit algebra $\A$.  As in
\cite{PittsKGrNeAl}, Corollary~\ref{idemsim} implies $\iota_*:
K_0(\D)\rightarrow K_0(\A)$ is an isomorphism of scaled dimension
groups and $\iota_*^{-1}=E_*.$
\end{remark}

We now show that an isomorphism of triangular limit algebras implies the
existence of a $*$-isomorphism of the \cstar-envelopes.  
\begin{theorem}\label{trienvelopeiso}
Suppose $\theta:\A_1\rightarrow \A_2$ is an algebra isomorphism of the
triangular limit algebras $\A_i$.  For $i=1,2$, let $\C_i$ be the
\cstar-envelope of $\A_i$ and let $h_i:\D_i\rightarrow \A_i$ and
$k_i:\A_i\rightarrow \C_i$ be the inclusion maps.  Then there exists
a $*$-isomorphism  $\tau: \C_1\rightarrow
\C_2$ such that the following diagram of scaled dimension groups
commutes.
\begin{equation}\label{comdia}
\begin{CD}
K_0(\D_1)@>h_{1*}>> K_0(\A_1) @>k_{1*}>> K_0(\C_1)\\
@VV\alpha_*V      @VV\theta_*V    @VV\tau_*V\\  
K_0(\D_2)@>h_{2*}>> K_0(\A_2) @>k_{2*}>> K_0(\C_2)
\end{CD}
\end{equation}
\end{theorem}

\begin{proof} 
Recall from~\cite{MR2001k:47103} that algebraic isomorphisms
of limit algebras are necessarily bounded.

That $\theta_*\circ h_{1*}= h_{2*}\circ \alpha_*$ follows from the
fact that  $\alpha=E_2\circ\theta|_{\D_1}$ and Remark~\ref{ktheormk}.  
For $j=1,2$, let $\iota_j=k_j\circ h_j$ be the inclusion mapping of
  $\D_j$ into $\C_j$.  To complete the proof, we shall show the
existence of $\tau_*:K_0(\C_1)\rightarrow K_0(\C_2)$ so that
$\tau_*\circ\iota_{1*}=\iota_{2*}\circ\alpha_*$.  

 Write $\C_j=\indlim_k(\C_{jk},\D_{jk})$ where $(\C_{jk},\D_{jk})$ are
  finite dimensional \cstar-inclusions.  Without loss of generality,
  we may assume that $\A_{jk}:=\A_j\cap \C_{jk}$ satisfies
  $C^*(\A_{jk})=\C_{jk}$.  Any projection $p\in\C_1$ is algebraically
  equivalent to a projection in $\C_{1k}$ for some $k$, so $p$ is
  algebraically equivalent to a projection $p'\in \D_{1k}$.  It
  follows that the induced mapping of scaled dimension groups,
  $(\iota_j)_*: K_0(\D_j)\rightarrow K_0(\C_j)$ is onto.

We claim that if $p$ and $q$ are projections in $\D_1$ which are
algebraically equivalent in $\C_1,$ then $\alpha(p)$ and $\alpha(q)$
are algebraically equivalent in $\C_2$, and we modify ideas of 
\cite[Lemma~2.2]{MR1241116} for this.  We may assume $p,q \in
\D_{1k}$ for some $k$, and are algebraically equivalent in $\C_{1k}$.
In fact, we shall show that they are equivalent via an element of
$\N_{\D_{1k}}(\C_{1k}).$ 

Since $p$ and $q$ are algebraically equivalent, they have the same
center-valued trace, and hence there exists a positive integer $r$ and
minimal projections $p_i$, $q_i$ belonging to $\D_{1k}$ so that $p=p_1
+ \cdots + p_r$, $q=q_1 + \cdots + q_r$.  By relabeling if necessary,
we may assume that for each $i$ with $1\leq i\leq r$, $p_i$ and $q_i$
are algebraically equivalent. Let $w_i\in\C_{1k}$ be a partial
isometry so that $q_iw_ip_i=w_i$, $w_i^*w_i=p_i$ and $w_iw_i^*=q_i$.
Since $p_i$ and $q_i$ are minimal projections in $\C_{1k}$, $w_i$ are
minimal partial isometries in $\C_{1k}$.  Moreover, $w=\sum_{i=1}^r
w_i\in\N_{\D_{1k}}(\C_{1k})$ satisfies $w^*w=p$ and $ww^*=q$.

 Since $C^*(\A_{1k})=\C_{1k}$, we can write
$w_i=v_{i_1}v_{i_2}\dots v_{i_l}$ as a finite product of partial
isometries, with each $v_{i_j}$ a normalizer of $\D_{1k}$ and
belonging to either $\A_{1k}$ or $\A_{1k}^*$.  Theorem~\ref{bigtechie}
shows that $\alpha(v_{i_j}^*v_{i_j})$ and $\alpha(v_{i_j}v_{i_j}^*)$
are algebraically equivalent in $\C_2.$ Thus, $\alpha(p_i)$ and
$\alpha(q_{i})$ are equivalent in $\C_2$ also, whence $\alpha(p)$
and $\alpha(q)$ are algebraically equivalent in $\C_2$ as desired.

Thus, if $p\in \C_1$ is a projection, we may define
$\tau_*([p])=(\iota_{2*}\circ\alpha_*)([p'])$, where $p'\in\D_1$ is any
projection in $\D_1$ with $\iota_{1*}([p'])=[p].$  The previous
paragraph shows $\tau_*$ is well-defined, and it determines an
isomorphism of the scaled dimension groups $K_0(\C_1)$ and
$K_0(\C_2)$ satisfying \eqref{comdia}.

An application of Elliott's Theorem now completes the proof.
\end{proof}

We would very much like to know whether it is possible to choose
$\tau$ in the conclusion of Theorem~\ref{trienvelopeiso} so that
$\tau|_{\D_1}=\alpha$.  When this is the case,
Corollary~\ref{envelopeiso} implies the existence of a continuous
isomorphism of coordinate systems, and hence spectra.  The next
example shows that more than the $K$-theoretic data provided by the
conclusion of Theorem~\ref{trienvelopeiso} is required to prove the
existence of such a $*$-isomorphism.

\begin{remark}{Example}\label{ex:noextend}
Suppose, for $j=1,2$, that $(\C_j,\D_j)$, are AF \cstardiag s 
and $i_j : \D_j \to \C_j$ are the natural inclusions.
Given an isomorphism $\alpha : \D_1 \to \D_2$ with
an isomorphism of scaled dimension groups, $h : K_0(\C_1) \to K_0(\C_2)$
with $i_{2*} \circ \alpha_* = h \circ i_{1*}$, there need not
exist a $*$-isomorphism $\tau : \C_1 \to \C_2$ with $\tau |_{\D_1} = \alpha$.

To see this, we use two well-known direct systems for triangular AF algebras,
the refinement system and refinement with twist system.
Define $\rho_k : M_{2^k} \to M_{2^{k+1}}$ by sending a matrix $A=[a_{ij}]$
to $\rho_k(A)=[a_{ij} \otimes I_2]$, i.e., replacing each entry of $A$ with
the corresponding multiple of a $2 \times 2$ identity matrix.
Define $\phi_k :  M_{2^k} \to M_{2^{k+1}}$ to be $\Ad U_k \circ \rho_k$,
where $U_k$ is the $2^{k+1} \times 2^{k+1}$ permutation unitary which is the
direct sum of a $2^{k+1}-2$ identity matrix and 
$\left[\begin{smallmatrix} 0 & 1 \\ 1 & 0 \end{smallmatrix}\right]$.

Let $\C_1 = \indlim(M_{2^k}, \rho_k)$ and $\C_2 = \indlim(M_{2^k}, \phi_k)$.
Let $\D_1$ and $\D_2$ be the direct limits of the diagonal matrices in each
direct system.
Since $\rho_k$ and $\phi_k$ agree on $D_{2^k}$, the direct limit of the
identity maps $\id : D_{2^k} \to D_{2^k}$ defines an isomorphism,
$\alpha$, from $\D_1$ to $\D_2$.
Now, $\C_1$ and $\C_2$ are isomorphic, as they are UHF \cstar-algebras
with the same `supernatural' number, $2^\infty$.
Further, $K_0(\C_j)$ can be identified with 
$G = \{ k/2^n : k \in \bbZ, n \in \bbN \}$, with the usual order and scale
$G \cap [0,1]$.
With this identification, $i_{j*}$ is the usual trace from $K_0(\D_j)$ into 
$G \subset \bbR$.
Since $\alpha$ is the identity on $K_0(\D_j)$, we have 
$i_{2*} \circ \alpha_* = i_{1*}$.
Thus, we can take $h$ to be the identity map on $G$.

It remains to show that there is no $*$-isomorphism 
$\tau : \C_1 \to \C_2$ with $\tau |_{\D_1} = \alpha$.  
We argue by contradiction, so assume such a $\tau$ exists. 

We may build an intertwining diagram as follows.  For brevity, let
$C_i$ be $M_{2^i}$, $\rho_{i,j}$ denote $\rho_i \circ \rho_{i+1} \circ
\dots \circ \rho_{j-1}$, and define $\phi_{i,j}$ similarly.  By
~\cite[Theorem~2.7]{MR0312282}, there are sequences $(m_i)$ and
$(n_i)$ and $*$-monomorphisms $(\psi_i)$ and $(\eta_i)$ so that the following
diagram commutes:
$$ \begin{CD}
C_1 @>{\rho_{1,m_1}}>> C_{m_1} @>{\rho_{m_1,m_2}}>> C_{m_2}
                                   @>{\rho_{m_2,m_3}}>> 
	C_{m_3} @>>> \cdots & \C_1 \\
& \symbse{\psi_1} & \symbup{\eta_1} & \symbse{\psi_2} & \symbup{\eta_2}
      & \symbse{\psi_3} & \symbup{\eta_3} & \symbse{\psi_4}
      & & \symbdown{\tau}\\
C_1 @>{\phi_{1,n_1}}>> C_{n_1} @>{\phi_{n_1,n_2}}>> C_{n_2} 
      @>{\phi_{n_2,n_3}}>> C_{n_3} @>>>  \cdots & \C_2
\end{CD} $$
Since $\tau$ maps $\D_1$ onto $\D_2$, we can use this diagram to show that
each $\psi_k$ and $\eta_k$ are restrictions of $\alpha$ and $\alpha^{-1}$,
respectively, and so are the identity map at the level of matrix algebras.

To obtain the contradiction, first fix $C_k=M_{2^k}$ and observe that
if $e$ is the $(1,1)$ matrix unit and $f$ the $(2^k,2^k)$ matrix unit
in $C_k$, then for any $l < k$, $e \rho_{l,k}(C_l) f = 0$ while $e
\phi_{l,k}(C_l) f \ne 0$.  Now consider the two maps $\lambda =
\phi_{n_1,n_3}$ and $\mu = \psi_3 \circ \rho_{m_1,m_2} \circ \eta_1$
from $C_{n_1}$ into $C_{n_3}$.  Letting $e$ and $f$ be the $(1,1)$ and
$(2^{n_3},2^{n_3})$ matrix units in $C_{n_3}$, the observation implies
that $e \lambda(C_{n_1}) f \ne 0$. 
To see that $e \mu(C_{n_1}) f = 0$, let $e'$ and $f'$ be the $(1,1)$
and $(2^{m_2},2^{m_2})$ matrix units in $C_{m_2}$ and observe that
$e' \rho_{m_1,m_2}(\eta_1(C_{n_1})) f' = 0$ by the observation.
Applying $\phi_3$ and noting that $e$,$f$ are subprojections of 
$\phi_3(e')$,$\phi_3(f')$ respectively completes the argument.

Thus, no such diagram exists, and hence no such $\tau$ exists.
\end{remark}

\providecommand{\bysame}{\leavevmode\hbox to3em{\hrulefill}\thinspace}
\providecommand{\MR}{\relax\ifhmode\unskip\space\fi MR }
\providecommand{\MRhref}[2]{%
  \href{http://www.ams.org/mathscinet-getitem?mr=#1}{#2}
}
\providecommand{\href}[2]{#2}

\end{document}